\documentclass[11pt,leqno]{amsart}

\usepackage{color,graphicx,times}

\makeatother

\let \ttorg \tt \def \tt{\ttorg \obeyspaces}

\makeatletter
 
\newcommand{\cross}[1]{%
\mbox{\vbox{\kern 1pt\hbox{\vbox{\hrule
\kern 2pt\hbox{\ensuremath{\vphantom{b}#1}\kern 2pt}}\vrule\kern 1pt}}}\,}
\newcommand{\rcross}[1]{%
\mbox{\vbox{\kern 1pt\hbox{\vbox{\hrule
\kern 2pt\hbox{\ensuremath{\vphantom{b}#1}\kern 2pt}}\vrule\kern 1pt}}}\,}
\newcommand{\lcross}[1]{%
\mbox{\vbox{\kern 2pt\hbox{\vrule\kern 0pt\vbox{\hrule
\kern 1pt\hbox{\ensuremath{\vphantom{b}#1}\kern 1pt}}   }}}   }

 \newcommand{\Fcross}{\raisebox{-0.25\height}{\includegraphics[width=0.5cm]{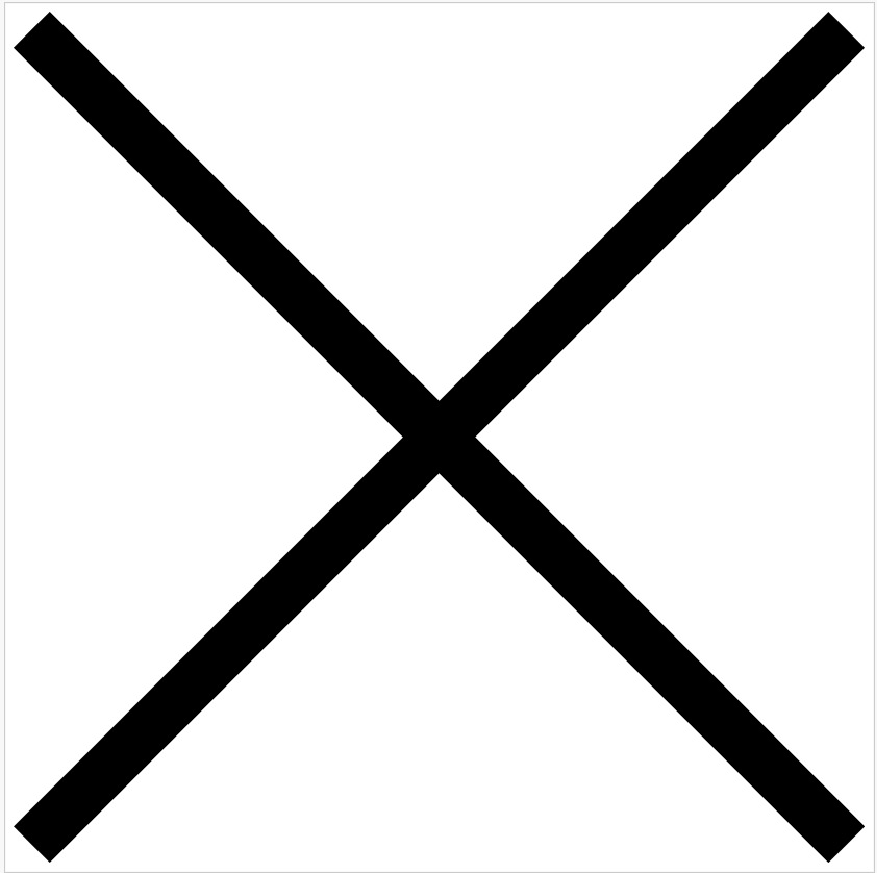}}}
\newcommand{\Across}{\raisebox{-0.25\height}{\includegraphics[width=0.5cm]{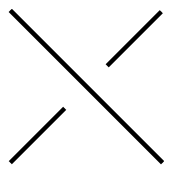}}}
\newcommand{\Bcross}{\raisebox{-0.25\height}{\includegraphics[width=0.5cm]{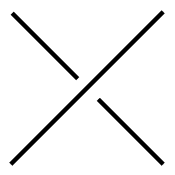}}}
\newcommand{\Asmooth}{\raisebox{-0.25\height}{\includegraphics[width=0.5cm]{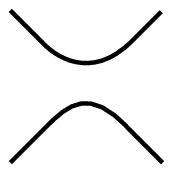}}}
\newcommand{\Bsmooth}{\raisebox{-0.25\height}{\includegraphics[width=0.5cm]{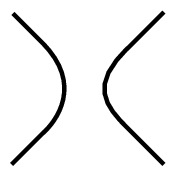}}}
\newcommand{\Rcurl}{\raisebox{-0.25\height}{\includegraphics[width=0.5cm]{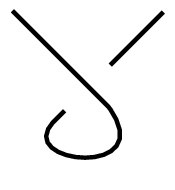}}}
\newcommand{\Lcurl}{\raisebox{-0.25\height}{\includegraphics[width=0.5cm]{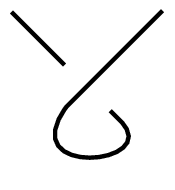}}}
\newcommand{\Arc}{\raisebox{-0.25\height}{\includegraphics[width=0.5cm]{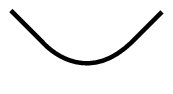}}}

\begin{document}
 
\date{}

\title{\bf Knot Logic and Arborescent Links}
\author{Louis
H. Kauffman\\ Department of Mathematics, Statistics \\ and Computer Science (m/c
249)    \\ 851 South Morgan Street   \\ University of Illinois at Chicago\\
Chicago, Illinois 60607-7045\\ $loukau@gmail.com$\\
and\\
International Institute for Sustainability with Knotted Chiral Meta Matter ($WPI-SKCM^{2}$)\\
Hiroshima University, 1-3-1 Kagamiyama, Higashi-Hiroshima, Hiroshima 739-8526, Japan}

\maketitle

\noindent{\bf Abstract.} This paper introduces a new algebra, the {\it crossing algebra}, that is applied to count the number of components for arborescent knots, links, tangles or states
(of a state polynomial expansion such as the Kauffman bracket). This algebra is elementary and foundational, and it is related to generalisations of boolean logic and to aspects of foundations based in diagrams and distinctions.
Applications are given to rational knots, links and tangles and to the structure of the bracket polynomial and the beginnings of Khovanov homology.\\



\section{Introduction}

Rational and aborescent links \cite{Alex,Conway,KL,KL1,KL2,KL3,SB,S} can be described by algebraic expressions that generalise continued fractions.
In this paper we point out that there is a simple algebraic-combinatorial way to determine whether a given aborescent link has one or many components, and to count the number of components. In the case of rational knots, there will be either one or two components and the 
algebraic method can be made even simpler. We call this algebraic method the {\it crossing algebra} and explain both how to use it to find the component count (by hand and by program), and how the crossing algebra is related to the indicational calculus of G. Spencer-Brown \cite{GSB} and to boolean and multiple valued logics. Precursors to the constructions in this paper can be found in \cite{VL,KnotLogic,GoldKauff}. \\

Finding the number of components of a link is a very simple question, but when the link is presented as an algebraic expression, such as a continued fraction, or as a diagram, it is advantageous to be have a
method other than tracing the edges of diagram to determine the component count. In the case of rational tangles and rational knots and links \cite{S}, the structure is determined by a continued fraction and in this case by the rational number that corresponds to this continued fraction. Let $[a_1,a_2,\cdots , a_n] = P/Q$ denote the rational value of the continued fraction $a_1 + 1/a_2 + 1/a_3 + \cdots +1/a_n .$ It is assumed that $P/Q$ is a fraction in reduced form so that $P$ and $Q$ are relatively prime.
Let $K= K(P/Q)$ denote the rational knot or link obtained by taking the numerator closure of the rational tangle associated with this continued fraction (The terminology is defined in Section 1 below.). In Section 3 we prove the {\it Fraction Theorem} that states that $K(P/Q)$ has two components if and only if $P$ is even, and that when $P$ is odd the question of the parity of $Q$ is decided by a calculation in the crossing algebra.\\

We give background on rational tangles and rational knots and links in Section 2 and continue this with the development of the crossing algebra in Sections 2 and 3. Section 3 gives many examples of the use of the crossing algebra. Along with the determination of component count, one can find those crossings in an arborescent link that are self-crossings of single components and those crossings that involve two components. The strategy for such 
determination involves seeing the crossing algebra expression for the link as a function of its local crossings. If the value of this function cannot change when a given variable is changed  (corresponding to smoothing the crossing in a diagram), then that crossing is a self-crossing. Thus the crossing algebra's opacity or transparency to transmission from its variables reflects on the topology of the corresponding link.\\

Other examples are pursued in Section 3. We point out that rational knots are determined by ordered partitions of positive integers such that the ends of the orderings have values greater than one. This is a well-known reformulation of the Schubert Theorem classifying rational knots and links. We combine this with the ease with which the crossing algebra tells whether a given partition yields a knot or a link. Then it is easy to generalte all rational knots and links with $n$ crossings in the minimal diagram (and know which ones are knots and which are links). We do the calculation here for $n$ less than or equal to 8. See Figures~\ref{part1},~\ref{part2},~\ref{twosix},~\ref{seven},~\ref{eight}. \\

At the end of Section 3 and continuing to Section 4 we discuss how to use the crossing algebra to compute the Kauffman bracket polynomial \cite{KNew,JO,JO1,Witten} and the details of $Mathematica^{TM}$ programs for this purpose.
The Appendix contains further details about the programming using string replacements.  We point out that for arborescent links, the crossing algebra yields a way to construct the Khovanov complex for a link from its crossing algebra expression. However, the full explication of the Khovanov complex is explained here by using a mapping of the crossing algebra to a an algebra of abstract tensors that are in principle directly codeable in a computer and directly interpretable as diagrams of the bracket states of the link. We make these constructions to explore the possiblity of purely algebraic approaches to the construction of the Khovanov complex for its own sake and for quantum computation of Khovanov homology.\\

Section 5 discusses component count in relation to medial graphs of plane graphs and so-called checkerboard graphs of link diagrams. The medial graph of a plane graph is a {\it flat} knot or link diagram, meaning that crossing types (over or under) have not been chosen. A flat link diagram has a component count equal to the component count of any given choice of knot or link obtained by making crossing choices. 
The problem of finding the number of components of the medial graph of a plane graph is a generalisation of our problem of component counts for arborescent links. We explain the details of this correspondence and point out how it is that the component count of the medial graph can be seen as the nullity of a graphical laplacian matrix.\\

Section 6 recalls the work of Claude Shannon who showed how boolean algebra applies to the structure of switching circuits. The classical problem for switching circuits is a connectivity problem: Is there a path in the network from one node to another. We explain how to design switching circuits and consider the problem to control one light  with $n$ switches. We show how a {\it crossing switch} can be used to solve this problem where a crossing switch has two input lines and two output lines and two states. In one state, the lines are parallel in a given instantiation as in $\Asmooth$. In the other state, two lines cross over one another as in $\Fcross.$
These are key elements in our crossing algebra. The details of the crossing algebra are given below, but here is a sample. Let $$O = \Fcross.$$Then $$OO = \Fcross \Fcross = \Asmooth = E$$ (The two crosses are joined at their middle tops and bottoms, as in tangle addition. The result is the same connectivity as in the single smoothing $E$.)The signs $O$ and $E$ are
the names of these iconic local possibilities in a network. We see from this that if $O$ is regarded as a switch then a simple linear connection of a row of $O$'s will suffice to control the light, via the parity of the connectivity. See Section 6 for more details.\\

Another word about the iconics of the crossing algebra: We use the notation $\cross{A}$ (`` cross A ")  for $1/A$ in ordinary algebra and for the $\pi/2$ turn of a tangle. Thus we have 
$$\cross{\cross{A}} = A$$ for any A and any tangle that is equivalent to itself under a $\pi$ rotation. The notation $\cross{A}$ is the analog of the negation of $A$ in logic.
$$\cross{\Asmooth} = \Bsmooth,$$ and $$\cross{\Bsmooth} = \Asmooth.$$ Since $E= \Asmooth$ acts as an identity element in our algebra, we can often replace $E$ by the empty word and write
$$\cross{E} = \cross{\,\,} = \Bsmooth.$$
Note also that if $O = \Fcross,$ then $$\cross{O} = \cross{\Fcross} = \Fcross = O.$$
$$\cross{O} = O.$$
This shows the parallel between the crossing algebra and mutiple valued logic at the point of the inclusion of fixed points for negation. In the well-known Lukasiewicz three valued logic we would have a {\it tertium non datur}
$@$ such that $@@= @$ and $\cross{@} = @.$ Here we have $O$ with $OO = ~~~~ $ (the empty word), and $\cross{O} = O.$ The crossing algebra excels at determining the connectivity of networks, including knots, links and tangles. Thus the crossing algebra can be seen as a development in the same direction as Shannon's switching network theory.\\

In this way the crossing algebra can be seen as a departure from boolean algebra, similar to but quite distinct from multiple valued logics and other constructions, as a method for analysing certain switching circuits and as a method for understanding connectivity and component count in topology.
While component count is very simple, the consequences of counting components is quite serious for the theory of knots and links. The Kauffman bracket polynomial and the Khovanov homology are based on the numbers and relations of loops related to the smoothing states of links.  \\

The present paper is foundational and gives new points of view for these invariants. This paper is part of a line of papers by the author in this theme
\cite{KnotLogic,CSR,FD,KC,FKT,KOK,KD,TL,KS,SRRF,K,KNew,KP,GoldKauff,VL,Eigenform,Reflexivity,Reflexphys,KaufDiag,MLogic} and we hope it will stimulate further thought along these lines.\\

\noindent{\bf Remark.} While we have not invoked virtual knot theory \cite{VKT1,VKT2,KaufDiag} in this paper, all the methods given here apply to the clear generalisation of classical arborescent knots links and tangles to the
virtual category. Such relationships will be explored in a subsequent paper.\\

\section{Recalling Rational and Arborescent Knots and Links}
\begin{figure}[htb]
     \begin{center}
     \begin{tabular}{c}
     \includegraphics[width=12cm]{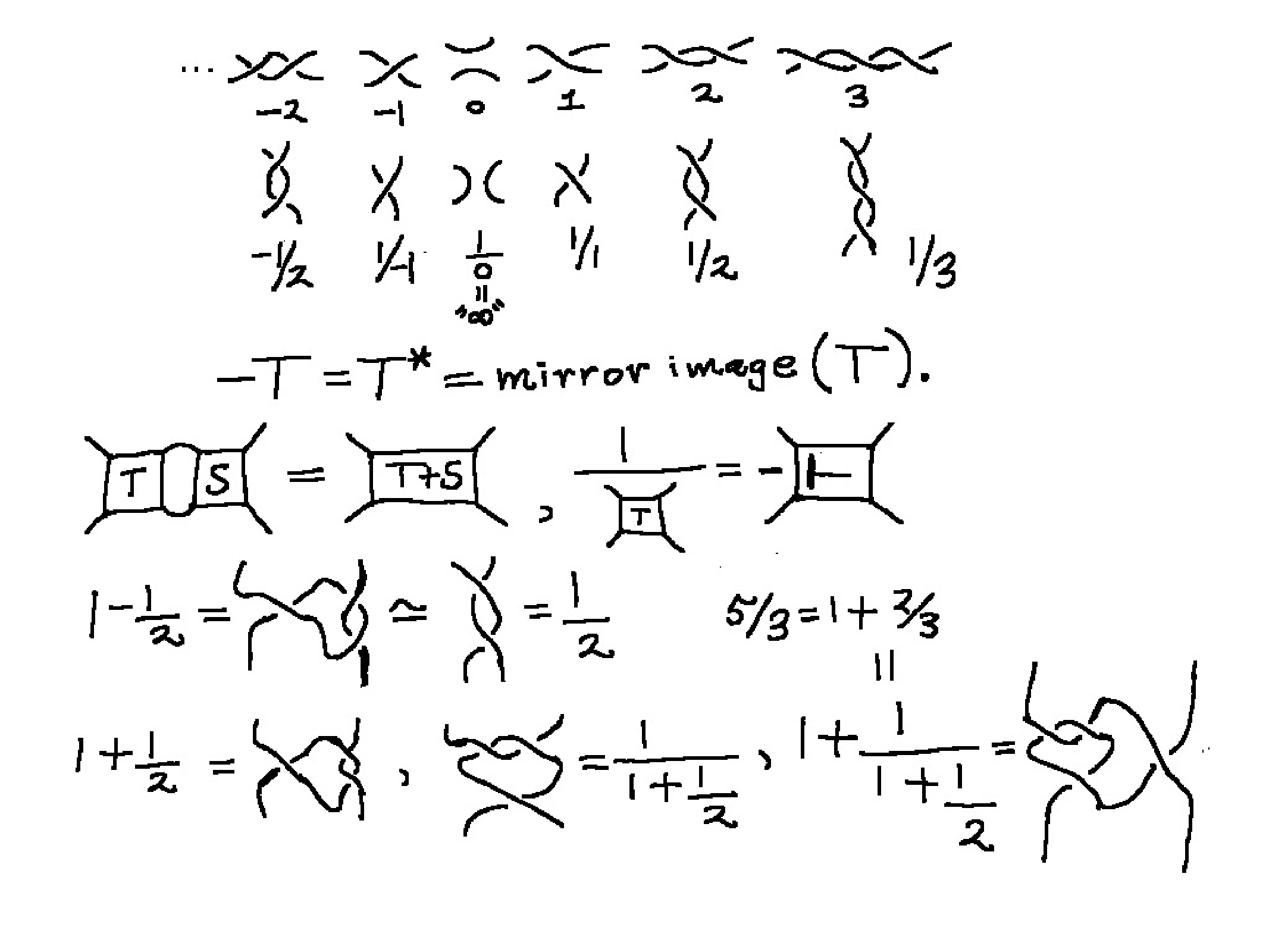}
     \end{tabular}
     \caption{\bf Tangles, Tangle Operations and Rational Tangles}
     \label{ratl1}
\end{center}
\end{figure}

We will use the concept of a {\it tangle} \cite{Conway}. A tangle is (represented by) a knot diagram with four free ends entering into a rectangular plane region wherein there is further diagram with no free ends. The four free ends are in the single outer region of the rectangle. One also can interpret tangles as embeddings in three dimensional space where the four free ends extrude from the boundary of a three-ball and there are embeddings of arcs and circles within the three-ball without free ends. Tangles sometimes are generalised to have a different number of free ends than four, but all tangles in this paper will have four free ends. Tangles are indicated as shown in Figure~\ref{ratl1}.
In that figure we illustrate that each tangle is shown with free ends in the four positions $\{ nw, ne,sw,se \}$ with $nw$ and $ne$ the upper left and right positions respectively, and $sw$ and $se$ the lower left and right positions. With this convention we can define the {\it sum}  $T+S$ of tangles $T$ and $S$ by attaching strands $ne(T)$ to $nw(S)$ and attaching strands
$se(T)$ to $sw(S),$ producing a new tangle whose strands are $\{ nw(T), sw(T), ne(S), se(S) \}.$\\

\begin{figure}[htb]
     \begin{center}
     \begin{tabular}{c}
     \includegraphics[width=12cm]{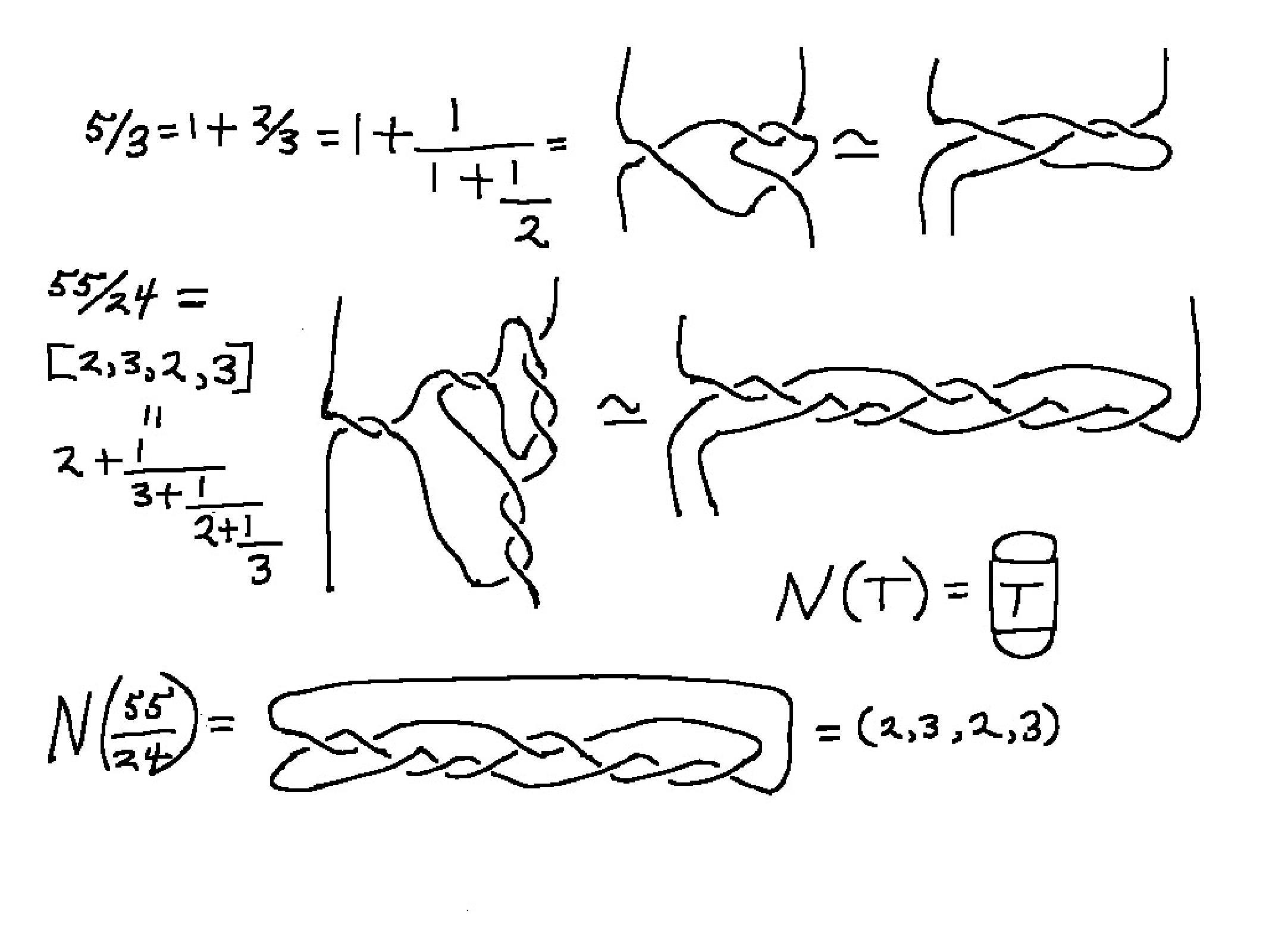}
     \end{tabular}
     \caption{\bf Forming Numerators}
     \label{ratl2}
\end{center}
\end{figure}

Two tangles $T$ and $S$ are {\it topologically equivalent} is there is an ambient isotopy fixing the tangle ends and restricted to the tangle box that makes one of them identcal to the other.
Equivalently, they are equivalent if there is a series of Reidemeister moves \cite{KOK,KNew,KP} taking one tangle to the other. No Reidemeister move is allowed to occur outside the tangle box.\\

\begin{figure}[htb]
     \begin{center}
     \begin{tabular}{c}
     \includegraphics[width=12cm]{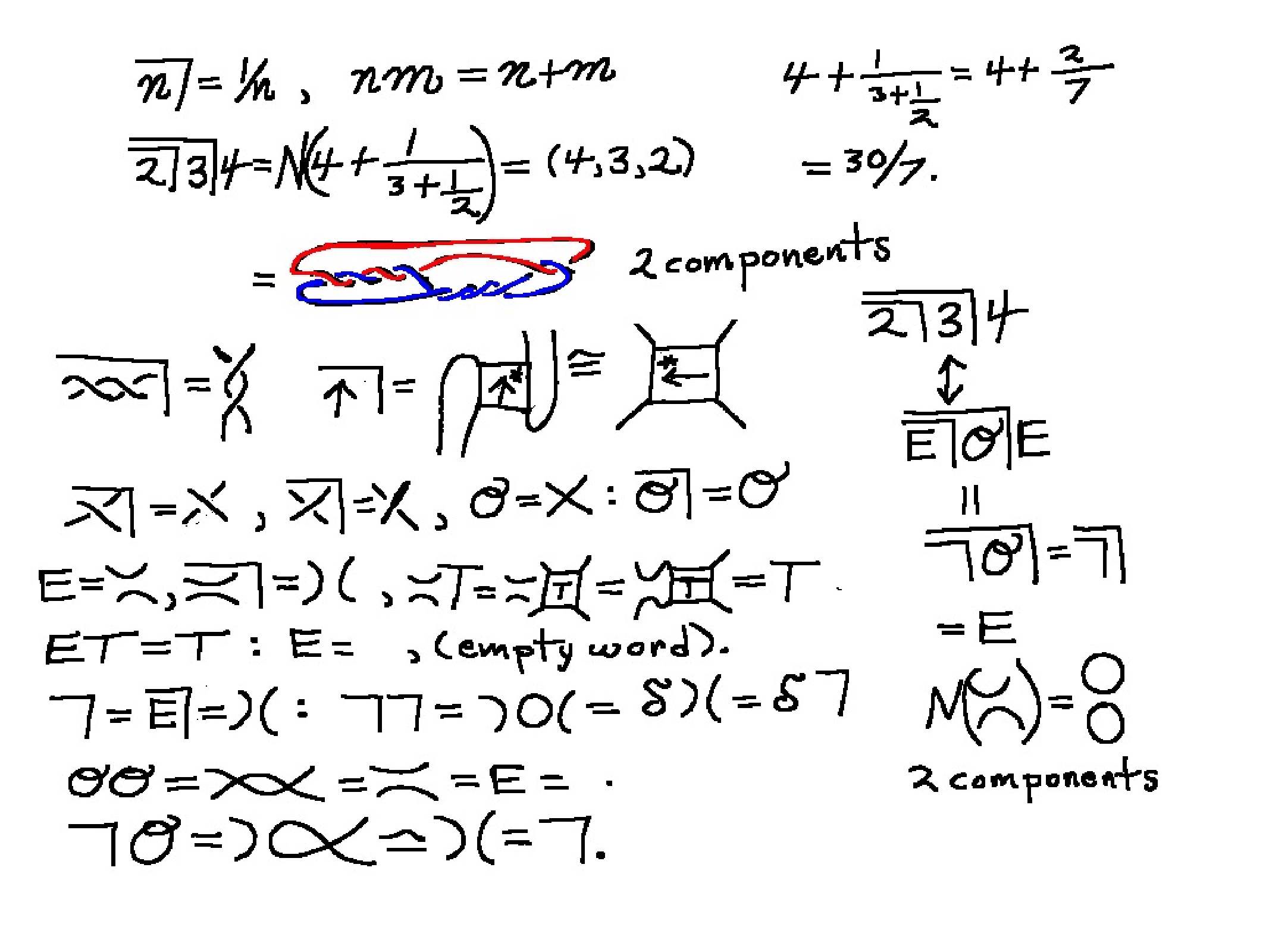}
     \end{tabular}
     \caption{\bf Mirror Operations and Crossing Algebra}
     \label{ratl3}
\end{center}
\end{figure}

\begin{figure}[htb]
     \begin{center}
     \begin{tabular}{c}
     \includegraphics[width=12cm]{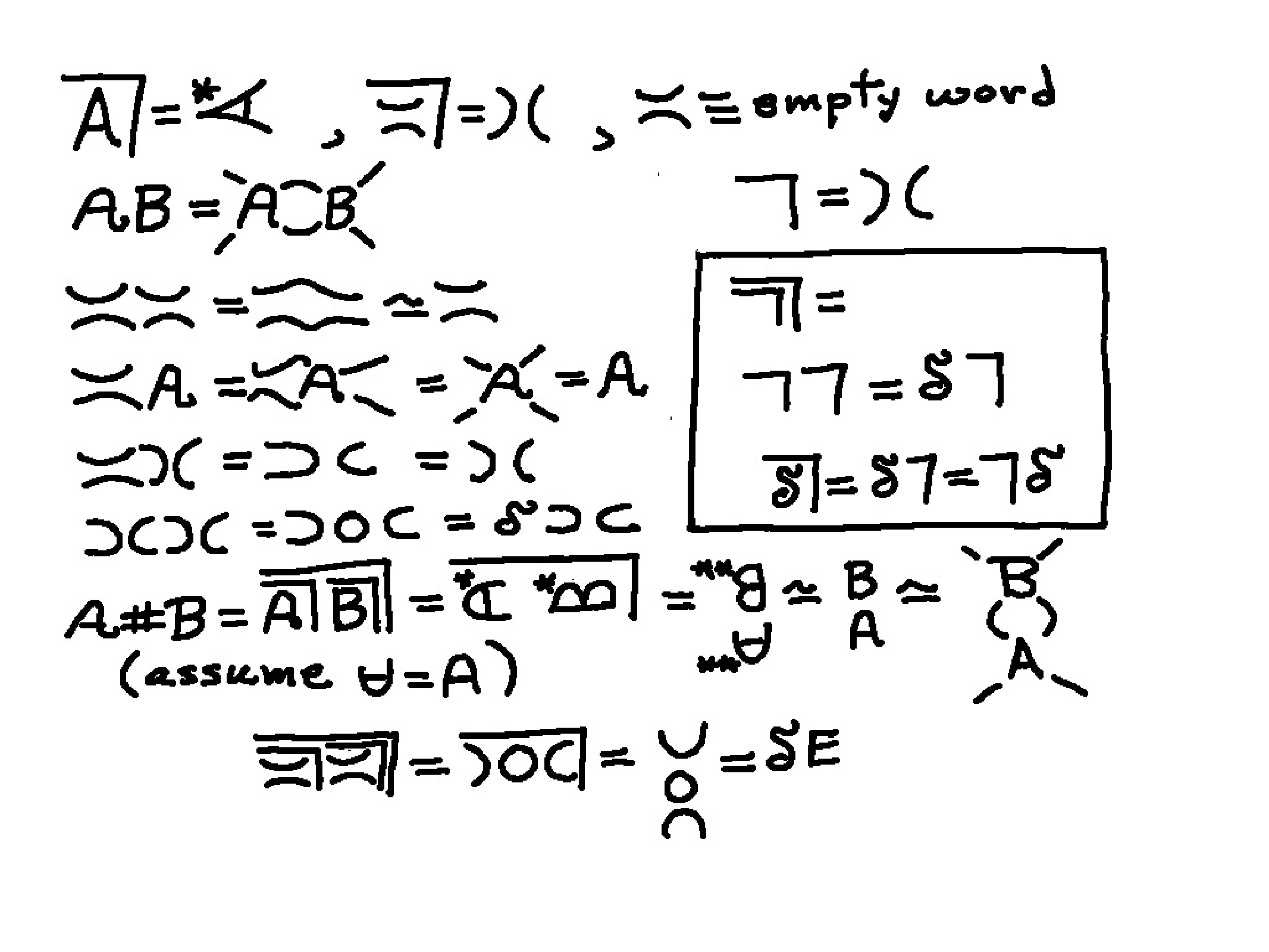}
     \end{tabular}
     \caption{\bf General Operations}
     \label{general}
\end{center}
\end{figure}

Along with the concept of addition, we have the notion of the {\it mirror rotation} of a tangle $T$ which consists in rotating the tangle by $\pi /2$ counterclockwise around a vertical axis through the center of the tangle box and perpendicular to the page of the diagram, and taking the mirror image of the result. See Figure~\ref{ratl1} and Figure~\ref{ratl3} for  illustrations of the mirror rotation operation. We shall denote the mirror rotation of a tangle $T$ as $\cross{T}.$ Thus
$\cross{\,\,}$ is our symbol for the mirror rotation operator. In the case of rational tangles $\cross{\,\,}$ has order two. In general, the mirror rotation operator has order four. In Figure~\ref{general} we review these operations once more. In particular we note that $ET = TE$ for any tangle $T$ and that 
$E$ can often be replaced by an empty word. In particular, we will write $\cross{\,\,} = \cross{E}$ and consequently we have $\cross{\,\,} \cross{\,\,} = \delta \cross{\,\,}$ where $\delta$ counts the extra loop that occurs in this tangle sum. We further note that for tangles $A,B$ where the operation $T \longrightarrow \cross{T}$ is of order two, the defined binary operation
$$A \sharp B =\cross{\cross{A}\cross{B}}$$ is realized by the vertical connection of the tangles $A$ and $B$ (where $AB$ is the horizontal connection of the tangles). These operations are tangle analogues of the dual operations of "or" and "and" in logic. We shall discuss this analogy further in Sections 6 and 7.\\

\begin{figure}[htb]
     \begin{center}
     \begin{tabular}{c}
     \includegraphics[width=12cm]{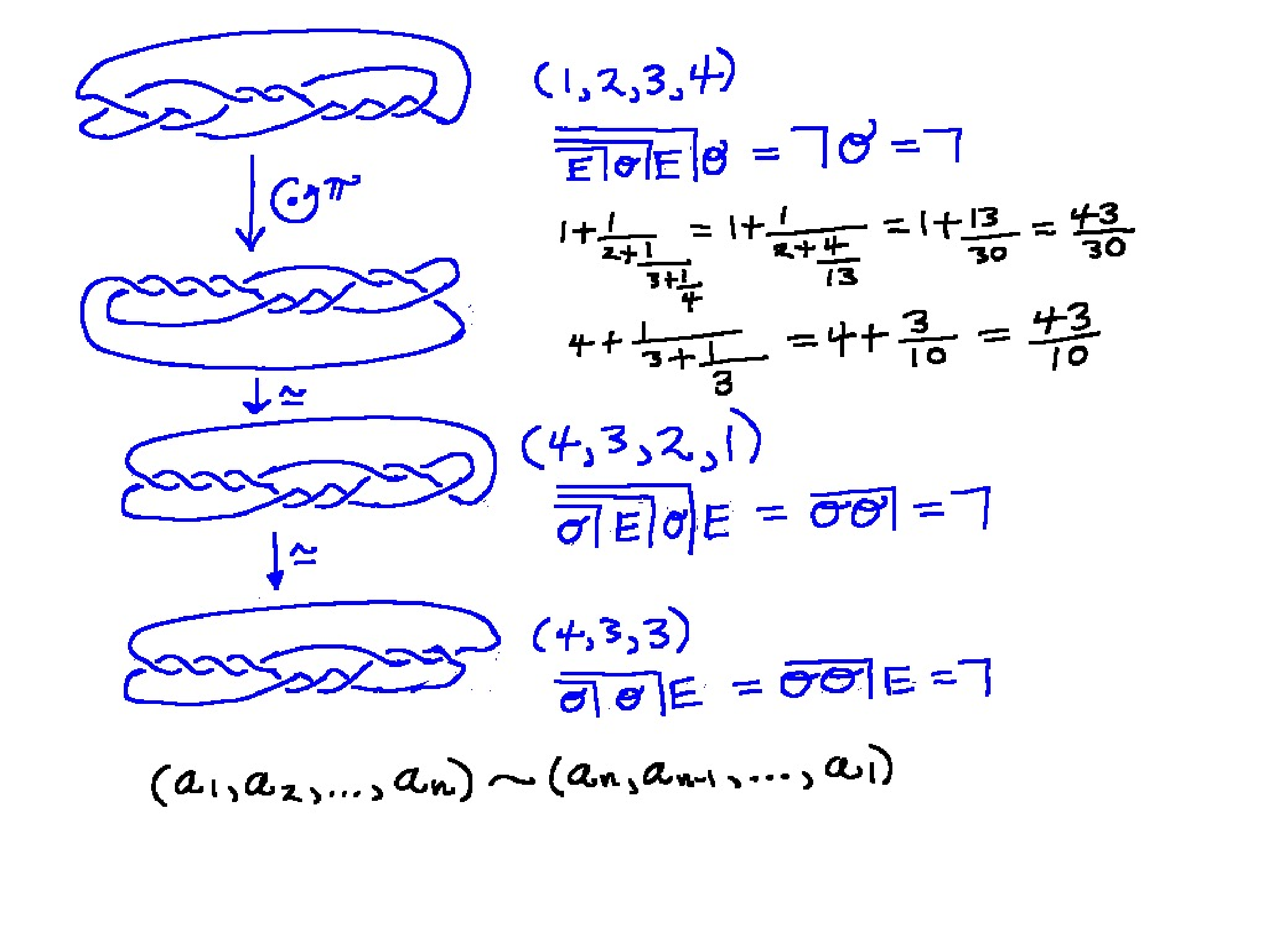}
     \end{tabular}
     \caption{\bf Order Reversal}
     \label{order}
\end{center}
\end{figure}

Given a tangle $T$, we construct a knot or link from the tangle that is called the {\it numerator} of $T$,  $N(T).$ The numerator of $T$ is constructed by identifying $ne(T)$ with $nw(T)$ and identifying $se(T)$ with $sw(T).$ See Figure~\ref{ratl2} for examples of this construction.  This figure illustrates ways to form tangles and numerators in the form of continued fractions that correspond also to braids. The numerators that are made by closing braids can be indicated by a sequence of integers as in $a = (a_1,a_2,\cdots , a_n).$ While we shall discuss these constructions below, the reader can apprehend them directly from the figure.
In Figure~\ref{order} we show that the numerator closure of $(a_1,a_2,\cdots , a_n)$ is the same (topologically and geometrically) as the closure of its reversal $(a_n,a_{n-1},\cdots , a_1).$ This is a key fact leading to the classification of the so-called {\it rational knots and links}, and we will use it in the discussion below.\\

\noindent {\bf Remark.} The notation $\cross{A}$ can be regarded as a shorthand for a box placed around $A$ as in $\fbox{A}.$ Thus we can write
$$ \cross{\cross{A}\cross{B}} = \fbox{\fbox{A} \fbox{B}}.$$ Some readers may find the ``box" notation more intutive and they are encouraged to write boxes when using this new algebra.\\

Rational and arborsecent links are composed from elementary integral tangles (See Figure~\ref{ratl1}), and in fact, these integral tangles can all be made from
the tangles shown below.
$$[0]=\Asmooth,  \,\, [\infty] = \Bsmooth,  \,\, [1]=\Across,  \,\, [-1]=\Bcross$$
Note that
$$\cross{[0]} = \cross{\Asmooth} = \Bsmooth = [\infty],$$
and
$$\cross{[1]} = \cross{\Across} = (\Bcross)^{\star} = \Across = [1].$$
$$\cross{[-1]} = \cross{\Bcross} = (\Across)^{\star} = \Bcross = [-1].$$
where $T^{\star}$ denotes the mirror image of the tangle $T$ via the diagram plane as the mirror.
The mirror rotation operator has order two on these basic tangles and leaves $[1]$ and $[-1]$ fixed.\\

The integral tangles $[n]$ are defined inductively by the equation $[n+1] = [n] + [1]$ and $[n-1] = [n] + [-1]$ so that, starting from $[0]$ one has tangles corresponding to each integer.
It is easy to see that topologically there is one integral tangle for each integer. For example the tangle $[1] + [-1]$ is topologically equivalent to $[0].$ Since an integral tangle has the appearance of a horizontal twist, its rotate has the appearance of a vertical twist as shown in Figure~\ref{ratl1}. Rational tangles are a special class of tangles that are obtained from the 
$[0]$ tangle by alternating horizontal and vertical twists. Thus we can imagine first creating all integral (horizontal twist) tangles $[n]$ and then creating all vertical twist tangles $\cross{[n]}.$
Then we can create tangles of the form $\cross{[b]} + [a].$ Tangles of this form can be interpreted as an integral tangle $[a]$ with a vertical twist of size $b$ made at the bottom. And one can make tangles of the form $\cross{\cross{[c]} + [b]} + [a]$ and so on in this pattern. In the illustration below and from now on, we let $a$ denote the integral tangle $[a].$
$$a$$
$$\cross{b} + a$$
$$\cross{\cross{c} + b} + a$$
$$\cross{\cross{\cross{d} + c} + b} + a$$
$$\cross{\cross{\cross{\cross{e} + d} + c} + b} + a$$
$$\cdots$$
The {\it rational tangles} are the tangles produced from integer tangles in this pattern.\\

\noindent {\bf Notation.} We can further abbreviate tangle operations by writing $TS$ instead of $T+S.$ With this, our chart of possible rational tangles takes the form:
$$a$$
$$\cross{b}a$$
$$\cross{\cross{c}b}a$$
$$\cross{\cross{\cross{d}c}b}a$$
$$\cross{\cross{\cross{\cross{e}d}c}b}a$$
$$\cdots$$
The inductive definition of rational tangles is 
\begin{enumerate}
\item Each integral tangle is a rational tangle.
\item If $T$ is a rational tangle and $a$ is an integral tangle, then $\cross{T}a = \cross{T} + a$ is a rational tangle.
\end{enumerate}

To each rational tangle $T$ there is an associated fraction $F(T)$ defined inductively by 
\begin{enumerate}
\item $F([n]) = \frac{n}{1}$
\item $F(\cross{[n]} =  \frac{1}{n}$ and $F(\cross{[0]} = F([\infty]) = \infty$ where $\infty$ is regarded as a formal infinite number whose arithmetic rules \cite{Conway} we shall discuss below using the iconic $\Bsmooth$ for $\infty = \cross{[0]} = \cross{\Asmooth} = \Bsmooth.$
\item If $T$ is a rational tangle for which $F(T)$ is defined, then $F(\cross{T} + [n]) = n + \frac{1}{F(T)}.$
We will write $$F(\cross{T} + a) = a+ \frac{1}{F(T)}.$$
\end{enumerate}

It follows from this definition that the fraction of a standard form of rational tangle is a continued fraction in the integral tangles that make it up. For example,

$$F(\cross{\cross{\cross{\cross{\cross{f}e}d}c}b}a)=F(\cross{\cross{\cross{\cross{\cross{f} + e} + d} + c} + b}+a) = a + \frac{1}{b + \frac{1}{c + \frac{1}{d + \frac{1}{e + \frac{1}{f}}}}}.$$\\

\noindent {\bf Definition.} Let $T(a)=[a_{1}, a_{2}, \cdots , a_{n}]$ denote the rational tangle with continued fraction
 $$a_{1} + \frac{1}{a_{2} + \frac{1}{a_{3} + \cdots + \frac{1}{a_{n}}}}.$$  For our purposes in this paper, this formula can be taken as the definition of a representative for any rational tangle.\\
  
 Conway's Theorem \cite{Conway,KL} states that {\it two rational tangles are topologically equivalent if and only if they have the same fraction.}\\
 
\noindent {\bf Definition.}. A {\it rational} knot or link is a link of the form $N(T)$  (numerator of $T$) where $T$ is a rational tangle. Using $T(a)$ as above, let
$K(a) = (a_{1}, a_{2}, \cdots , a_{n})$ denote the rational knot or link obtained as $N(T(a)).$ \\

See Figure~\ref{ratl2} for an illustration for the fraction
$55/24 = 2 + 1/(3 + 1/(2+ 1/3)).$ In that figure we denote the continued fraction by the notation $[2,3,2,3]$ and the numerator closure (the corresponding rational knot) by the notation
using curved brackets $(2,3,2,3).$\\

 Rational knots and links are classified by their continued fraction forms up to the following equivalences.
 \begin{enumerate}
 \item We have the following equality of rational links: $$(a_{1}, a_{2}, \cdots , a_{n}) = (a_{n}, a_{n-1}, \cdots , a_{1}).$$ 
 \item The endpoints $a_{1}$ and $a_{n}$ are greater than $1$.
 \item Every rational number in reduced form $P/Q$ with $P > Q > 0$ has a continued fraction of the
 form $[a_{1}, a_{2}, \cdots , a_{n}] $ with all positive terms, and $a_n > 1.$ If $a_1 = 1$ then 
 $(1, a_{2}, \cdots , a_{n}) = (a_{n}, a_{n-1},\cdots, a_{3}, a_{2}+1).$ Thus every rational number of this type corresponds to a unique rational knot or link $K(P/Q).$
 \item It is a fact of continued fractions that if $P/Q = [a_{1}, a_{2}, \cdots , a_{n}] $ then $[a_{n}, a_{n-1}, \cdots , a_{1}]  = P/Q'$ where
 $QQ' \equiv (-1)^{n+1} mod(P).$ This provides the connection between this description of the classification and the classical Theorem of Schubert \cite{KL}.
 \end{enumerate}

\noindent {\bf Remark.} The Schubert Theorem states that two reduced fractions $P/Q$ and $P'/Q'$ represent the same rational link if and only if $P=P'$ and either $QQ' \equiv \pm 1 mod(P)$ or $Q \equiv Q' mod(P).$
See \cite{KL}. \\

It follows from these remarks that rational knots are in 1-1 correspondence with sequences $(a_{1}, a_{2}, \cdots , a_{n})$ with positive entries and with $a_1$ and $a_n$ greater than $1$ where we identify $(a_{1}, a_{2}, \cdots , a_{n})$ with  $(a_{n}, a_{n-1}, \cdots , a_{1}).$ This means that we can enumerate rational knots and links by enumerating such sequences.
In the next section we show how to determine, without tracing diagrams, whether such a sequence is a knot or a link.\\

\section{Crossing Algebra}
Given a rational link or, more generally, an arborescent link, we wish to determine the number of components in the link. The case of rational links is special in that a rational link has either one component or two components. We will give an algebraic method for determining the number of components. As we have explained in Section 2, an arborescent link can be encoded as an arbitrary expression $F[a_1,a_2, \cdots , a_n]$ involving the cross operator $\cross{\,\,}$ and the indicated variables. For example, we may take 
$F = \cross{2} \cross{3} \cross{5}$ corresponding to a pretzel knot of type $(2,3,5)$. An example of the type of problem we should like to solve is to find the number of components of a 
generalised pretzel link of type $\cross{a_1}\cross{a_2}\cross{a_3} \cdots \cross{a_n}.$\\

It is clear that the component count depends only upon the parities of the integers $a_{k}.$ Accordingly, let $O$ denote ``odd" and let $E$ denote ``even".\\
Then we have the following rules for combinations of $O$ and $E.$\\
\\
\\
\\

\noindent {\bf Crossing Algebra Rules.}
$$EE=E$$
$$EO=OE=O$$
$$OO = \,\,\,$$
$$\cross{\,\,} O = \cross{\,\,}$$
$$\cross{\,\,} E = \cross{\,\,}$$
$$\cross{E} = \cross{\,\,}$$
$$\cross{O} = O$$
$$\cross{\cross{\,\,}} = \,\,\,$$
$$\cross{\,\,}\cross{\,\,} = \delta \cross{\,\,}$$\\

\begin{figure}[htb]
     \begin{center}
     \begin{tabular}{c}
     \includegraphics[width=12cm]{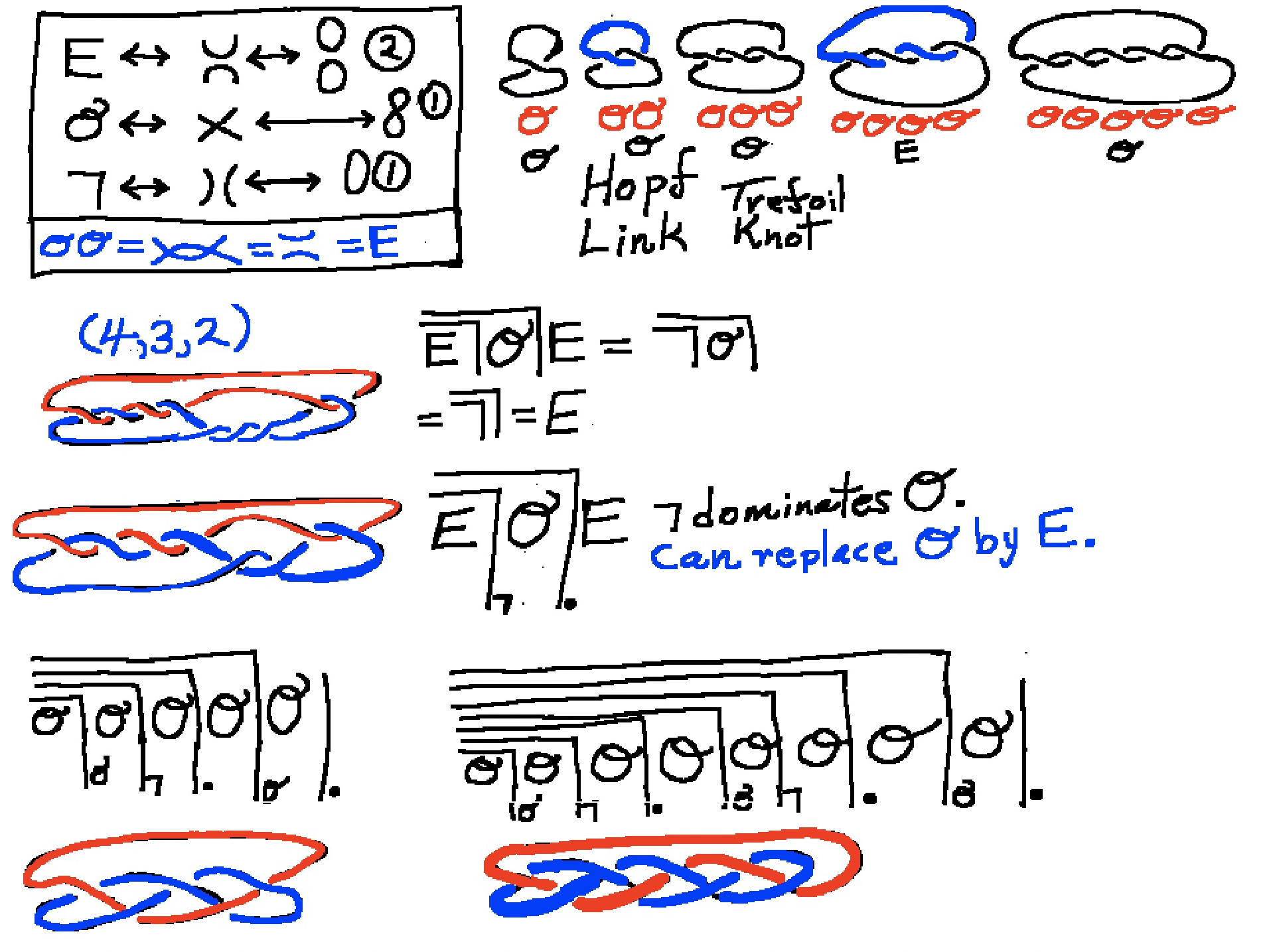}
     \end{tabular}
     \caption{\bf Even (E) and Odd (O) determine parity of components in the numerator.}
     \label{ratl4}
\end{center}
\end{figure}

See Figure~\ref{ratl3} and Figure~\ref{ratl4} for the diagrammatics for these identities. We describe this corresondence in detail below.\\

In order to see how these rules arise let $$O= \Fcross$$ $$E= \Asmooth$$ $$\cross{E} =\Bsmooth$$
These glyphs are representative odd, even and inverted even tangles. We can then consider their combinations.
$$OO = \Fcross \Fcross = E =\Asmooth$$
since the tangle sum of two odd integer tangles is an even integer tangle.
Similarly, we have
$$EE = \Asmooth \Asmooth = \Asmooth = E$$
since the sum of two integer even tangles is even, and
$$EO= \Asmooth \Fcross = \Fcross = O$$ since the addition of an even tangle to an odd tangle yields and odd tangle.
Note that $E$ behaves as an identity in this algebra. Thus we can use the empty word for $E$ as in
$$\cross{\,\,} = \cross{E} = \cross{\Asmooth} = \Bsmooth.$$
$\cross{\,\,}$ stands for the rotate of a horizontal smoothing.
We have
$$\cross{\,\,}E = \cross{\Asmooth}\Asmooth = \Bsmooth \Asmooth = \Bsmooth = \cross{\,\,},$$
$$\cross{\,\,}O =  \cross{\Asmooth}\Fcross = \Bsmooth \Fcross = \Bsmooth =  \cross{\,\,},$$
$$\cross{O} = \cross{\Fcross} = \Fcross$$
and
$$\cross{\cross{\,\,}} = \cross{\cross{E}}= \cross{\Bsmooth} = \Asmooth = E = \,\,.$$
Thus, using the empty word on the right, we have
$$\cross{\cross{\,\,}} = \,\,.$$
Finally, note that the concatentation of the rotated identity produces a loop in the middle.
Letting $\delta$ denote this loop, we have
$$\cross{\,\,}\cross{\,\,} =\cross{\Asmooth}\cross{\Asmooth} = \Bsmooth \Bsmooth = \delta \Bsmooth = \delta \cross{\,\,}.$$
This completes our iconic verification of the crossing algebra identities. We can now see, by using the algebra, that one can determine the number of components in an
arborescent link, from its structural specification in terms of tangle concatenations.\\

 In tangle calculus, $\cross{T}$ represents the ninety degree turn of the tangle combined with taking its mirror image, so that for a
 tangle fraction $P/Q$ we have $\cross{P/Q} = \frac{1}{P/Q} = Q/P$ and generally for a number or ordinary algebraic variable x,
 $\cross{x} = 1/x.$ In the crossing algebra we take AB to mean the analog of A+B. Thus $\cross{y}x=x\cross{y} = x + \frac{1}{y}.$  
  Note that $Odd + Odd = Even$ corresponds to the equation $OO=\,\,$ and $Odd + Even = Odd$ and $Even + Even = Even$
 correspond to $OE= O$ and $EE = E =\,\,.$\\
 
 We note that if the final evaluation of a tangle is $E$, than its numerator closure has two loops and so the numerator closure will be a link. If the final evaluation is $O$ (for example $\Across$ or $\Bcross$) or if it is  $\cross{\,\,} = \Bsmooth$ then the numerator closure will be a single component. And so the numerator closure will be a knot in these two cases.\\
 
\noindent{\bf Example.} The expression $$R = \cross{\cross{\cross{\cross{E}D}C}B}A$$ represents the continued fraction $$A + \frac{1}{B + \frac{1}{C+ \frac{1}{D+\frac{1}{E}}}},$$ seen either as a numerical fraction with $A,B,C,D$ integers or as the corresponding rational tangle. We further consider the numerator closure of that rational tangle and regard these expressions as representatives of the numerator closure. With this we see how to find the component count using the crossing algebra. For example, suppose that $A,B,C,D,E$ are odd. Then consider the expression
\\
\\

$$\cross{\cross{\cross{\cross{O}O}O}O}O$$ 
$$=\cross{\cross{\cross{OO}O}O}O$$
$$=\cross{\cross{\cross{\,\,}O}O}O$$
$$=\cross{\cross{\cross{\,\,}}O}O$$
 $$=\cross{O}O$$
 $$=OO$$
 $$=E$$
 
 Hence, since the numerator of an even twist has two components, we conclude that $R$ has two components.\\

More generally, let $K(a)=[a_{1}, a_{2}, \cdots , a_{n}]$ denote the numerator closure of the rational tangle with continued fraction
 $$a_{1} + \frac{1}{a_{2} + \frac{1}{a_{3} + \cdots + \frac{1}{a_{n}}}}.$$ Then we can determine whether $K(a)$ is a knot or a link
 by computing the the parities of the integers $a_{k}$ in the crossing algebra. Specifically, let $e_{i} = O$ if $a_{i}$ is odd, and let
 $e_{i} = E$ if $a_{i}$ is even.  Then the crossing algebra expression
 $$Cross(K(a))=\cross{\cross{\cross{\cross{e_{n}}e_{n-1} \cdots }e_{3}}e_{2}}e_{1}$$ will evaluate to either $O$ or $\cross{\,\,}$ for a knot, and $E$ for a link, determining the connectivity of the rational knot or link.\\

\subsection {\bf Rational Counting.} Rational knots and links are classified by their continued fraction forms as we explained in the previous section.
 Every rational number in reduced form $P/Q$ with $P > Q > 0$ has a continued fraction of the
 form $[a_{1}, a_{2}, \cdots , a_{n}] $ with all positive terms, and $a_n > 1.$ If $a_1 = 1$ then 
 $[1, a_{2}, \cdots , a_{n}] \sim [a_{n}, a_{n-1},\cdots, a_{3}, a_{2}+1].$ Thus every rational number of this type corresponds to a unique rational knot or link $K(P/Q).$\\
 
 \noindent {\bf Fraction Theorem.} {\it The rational link $K= K(P/Q)$ with $P/Q$ as above is a link of two components if and only if $P$ is even.
When $P$ is odd, $K(P/Q)$ has one component. If $P/Q$ has continued fraction expansion $[a_1,\cdots, a_n ],$ let $$C = \cross{\cross{\cross{e_n} e_{n-1}} \cdots e_1 }$$ be its parity expression in the crossing algebra, where $e_k = E$ when $a_k$ is even and $e_k=O$ is odd when $a_k$ is odd. In the crossing algebra we have:
\begin{enumerate}
\item $C = E$ if and only if $K$ has two components.
\item $C= O$ if and only if $K$ has one component and the denominator $Q$ is odd. 
\item $C=\cross{\,\,}$ if and only if $K$ has one component and the denominator $Q$ is even.
\end{enumerate}}
 
 \noindent {\bf Proof.} We will prove this result by induction. Fractions are of the type $e/o, o/o, o/e$ where here we use the symbols $e$ and $o$ as shorthand for {\it even} and {\it odd}.
 The induction hypothesis is: \\
 
 {\it The continued fraction cross-algebra evaluation for $e/o$ is $E$, for $o/o$ is $O$, for $o/e$ is $\cross{\,\,}.$} \\
 
 Base cases are easy to check:
 $2/1$ has crossing algebra expression $E.$ $1/1$ has crossing algebra expression $O.$ $1/2 = 0 + 1/2$ has crossing algebra expression $\cross{E}O = \cross{\,\,}O = \cross{\,\,}.$ 
 Thus the induction hypothesis is satisfied at the base. We now verify the induction step by taking each fraction type in turn and checking that the induction hypothesis remains satisfied in each case after we add either an even or an odd integer to a given fraction.\\
 
 \begin{enumerate}
 \item $e + 1/(e/o) = e + o/e = o/e$ and $\cross{E}E = \cross{\,\,}.$ 
 
 \item $o + 1/(e/o) = o + o/e = o/e$ and $\cross{E}O = \cross{\,\,}O = \cross{\,\,}.$
 
 \item$e + 1/(o/o) = e + o/o = o/o$ and $\cross{O}E = OE = O.$
 
 \item $o + 1/(o/o) = o + o/o = e/o$ and $\cross{O}O = OO = E.$
 
 \item $e + 1/(o/e) = e + e/o = e/o$ and $\cross{\cross{\,\,}}E = E.$
 
 \item $o + 1/(o/e) = o + e/o = o/o$ and $\cross{\cross{\,\,}}O = O.$\\
 \end{enumerate}

 Since each new induced fraction continues to satisfy the induction hypothesis, this completes the proof of the Theorem.
 QED.\\
 
  \noindent{\bf Example.}   As we stated in the previous section, if $p/q= [a_{1}, a_{2}, \cdots , a_{n}] $ and $[a_{n}, a_{n-1}, \cdots , a_{1}]  = p'/q'$ then $p=p'$ and 
 $qq' \equiv (-1)^{n+1} mod(p).$  Thus if we know the continued fraction for $p/q$ then we can determine both the parity for $q$ and the parity for $q'$ from the crossing algebra.
 For example, if we have $p/q= 355/113 = [3,7,16]$ then we have the corresponding parity expression in crossing algebra  $\cross{\cross{E}O}O =  EO  = O.$  The continued fraction $[16,7,3]$ has the parity expression
 $\cross{\cross{O}O}E = \cross{\,\,}E = \cross{\,\,}. $ This shows that $q$ and $q'$ have different parity.  Without further calculation we then know that $[3,7,16] = 355/q'$ with $q'$ even. This agrees with the direct calculation $ [16,7,3]  = 16 + 1/(7+ 1/3) = 355/22.$\\  
 
\begin{figure}[htb]
     \begin{center}
     \begin{tabular}{c}
     \includegraphics[width=10cm]{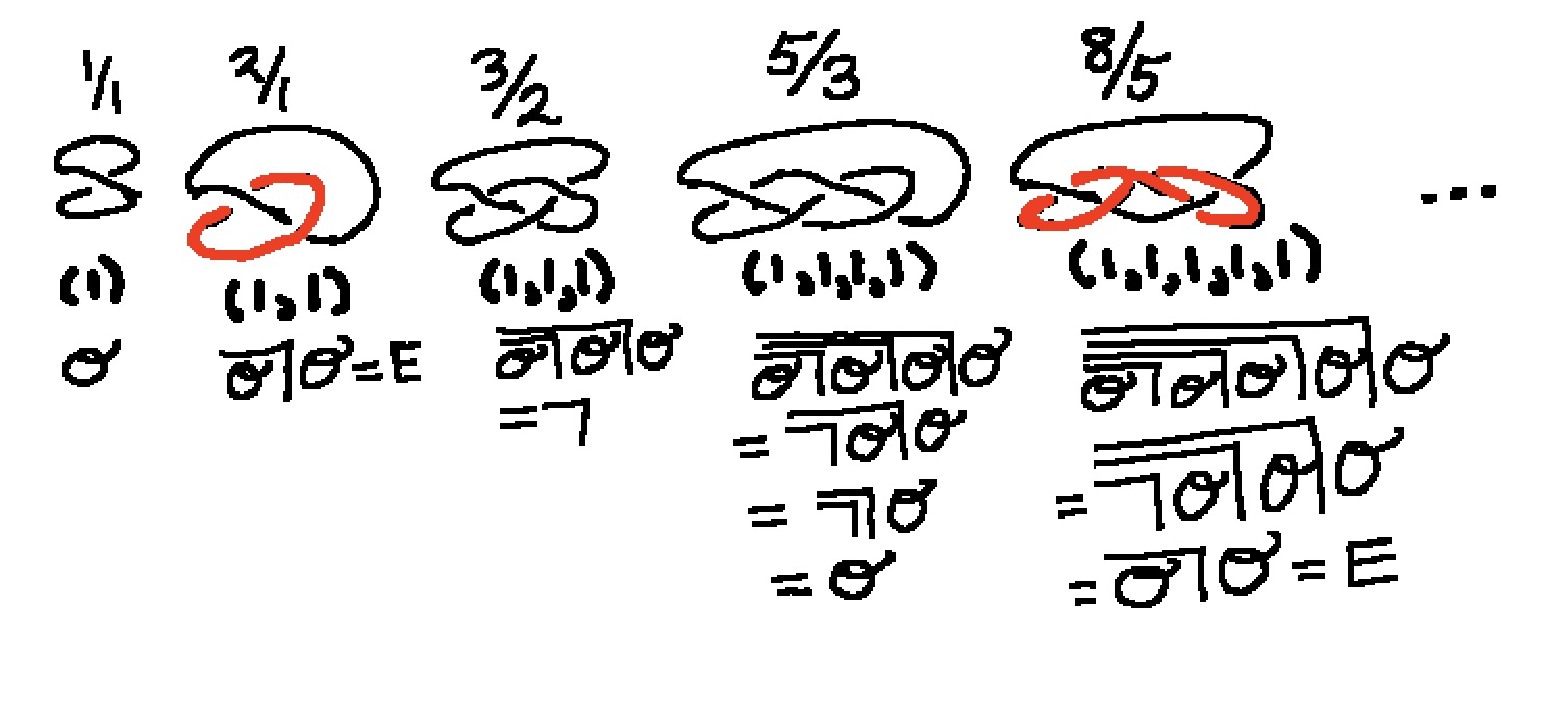}
     \end{tabular}
     \caption{\bf The Fibonacci Rational Knots and Links}
     \label{fibo}
     \end{center}
     \end{figure}

 \noindent{\bf Example.}
Consider the knots and links in Figure~\ref{fibo}. These are rational knots and links corresponding to the fractions $1/1, 2/1, 3/2, 5/3, 8/5$ and so can be called ``Fibonacci" rationals, as the Fibonacci sequence $\{f_{n}\}$ is
$1,1,2,3,5,8,13,21,34,\cdots$ where $f_{1}= f_{2}=1$ and the equation $f_{n+1} = f_{n} + f_{n-1}$ defines the sequence inductively. The Fibonacci rational knots and links are $\{ K_{n}\}$ where $K_{n}$ corresponds to the fraction $f_{n+1}/f_{n}.$ Note also that every third Fibonacci number is even. Thus every third $K_{n}$ will be a link of two components. It is of interest to see how the component count works out in the crossing algebra.
The continued fraction representations of the Fibonacci knots and links are uniform: 
$K_{1} = (1), K_2 = (1,1), K_3 = (1,1,1)$ and $K_n$ is represented by the continued fraction $(1,1,\cdots, 1)$ with $n$ $1$'s. This means that in crossing algebra $K_{n}$ is represented by
$\cross{\cross{{\cross{O}O}}O\cdots O}O$ with $n$ appearances of $O$ in the expression.
Here are the first few of them.

\begin{enumerate}
\item $K_{1}: O$
 \item $K_{2} : \cross{O}O = OO =E$
 \item $ K_{3} : \cross{\cross{O}O}O = \cross{OO}O = \cross{\,\,}O = \cross{\,\,}$
 \item $K_{4} :  \cross{\cross{\cross{O}O}O}O = \cross{\cross{\,\,}}O = O$
  \item $ K_{5} : \cross{\cross{\cross{\cross{O}O}O}O}O = \cross{O}O = OO =E$
  \end{enumerate}
 
 This shows that since all the terms in $K_{n}$ are odd, then $K_{n}$ will be a link when $n = 2,5,8,11,14,17...$. That is,
  $n = 3k+2$ for $k=0,1,2,3,\cdots.$\\

  It is natural to ask, for an arbitrary continued fraction $a = (a_1,a_2,\cdots , a_n)$, when $K(a)$ is a link of two components, which $a_{i}$ correspond to self-crossings of one of the  link components of $K(a).$
  This information is encoded in the crossing algebra.  To see how this works, consider
  $$K = \cross{\cross{\cross{\cross{O}O}O}O\,}\, O.$$
  The direct evaluation of $K$ proceeds as follows:
  $$K = \cross{\cross{\cross{\cross{O}O}O}O\,}\, O$$
  $$= \cross{\cross{\cross{OO}O}O\,}\, O$$
$$= \cross{\cross{\cross{\,\,}O}O\,}\, O$$
$$= \cross{\cross{\cross{\,\,}}O\,}\, O$$
$$= \cross{ O\,}\, O$$
$$= OO = E = u.$$

Here $u$ denotes unmarked (empty word, even) or equivalently, $E.$ We can use $m$ for the marked value ($\cross{\,\,}$) and $o$ for
$O$ as a value (odd).  Then the computation above can be summarized succinctly on the original expression by subscripting each mark with the value that emerges from its inside. See the expression below.

 $$K = \cross{\cross{\cross{\cross{O}_{o}O}_{m}O}_{u}O\,}_{o}\, O \, |_{u}$$
  
Here we use abbreviations at the subscripts with $o=O, m=\cross{\,\,}, u = E =\,\, .$
We write $$\cross{\cross{O}O} = \cross{\cross{O}_{o}O} = \cross{\cross{O}_{o}O}_{m}$$
Because $\cross{O}$ has value $o$ and $OO$ has value $u$ so that $\cross{OO} = \cross{\,\,}$ has value $m.$
In this way, each expression can be regarded as a tree structure that processes values from inside to outside, from the leaves of the tree to
its root.  With this we see that the assignment of $O$ or $E$ to  certain spaces in the expression do not affect the final value because a marked value emerges in that space in the course of this process evaluation. 
Thus in $K = \cross{\cross{\cross{\cross{O}O}O}O\,}\, O.$ we see from the above calculation that we can change the third $O$ from the left to an $E$ (change a twist from odd to even) without changing the component count of the numerator closure.\\

It is useful to write the computation in this succinct form and we can save even more notational noise by compressing it further by the rules
 $$a\}b = \cross{a}b,$$
 $$a\}b\}c = \cross{\cross{a}b}c,$$
 $$a\}b\}c\}d = \cross{\cross{\cross{a}b}c}d.$$
 In this compressed notation we can replace a nest of crossings by the right operator symbol $\}$ and then write the compressed calculation in the form below.\\

 $$O\}_{o}O\}_{m}O\}_{u}O\}_{o}O |_{u}$$
 
 In some cases it is convenient to notate the above just using a right-angle bracket, with the understanding that it is being used outside the usual conventions of the crossing algebra.
 Then we can write
 $$\cross{O}_{o}\cross{O}_{m}\cross{O}_{u}\cross{O}_{o}O |_{u}.$$

\begin{figure}[htb]
     \begin{center}
     \begin{tabular}{c}
     \includegraphics[width=8cm]{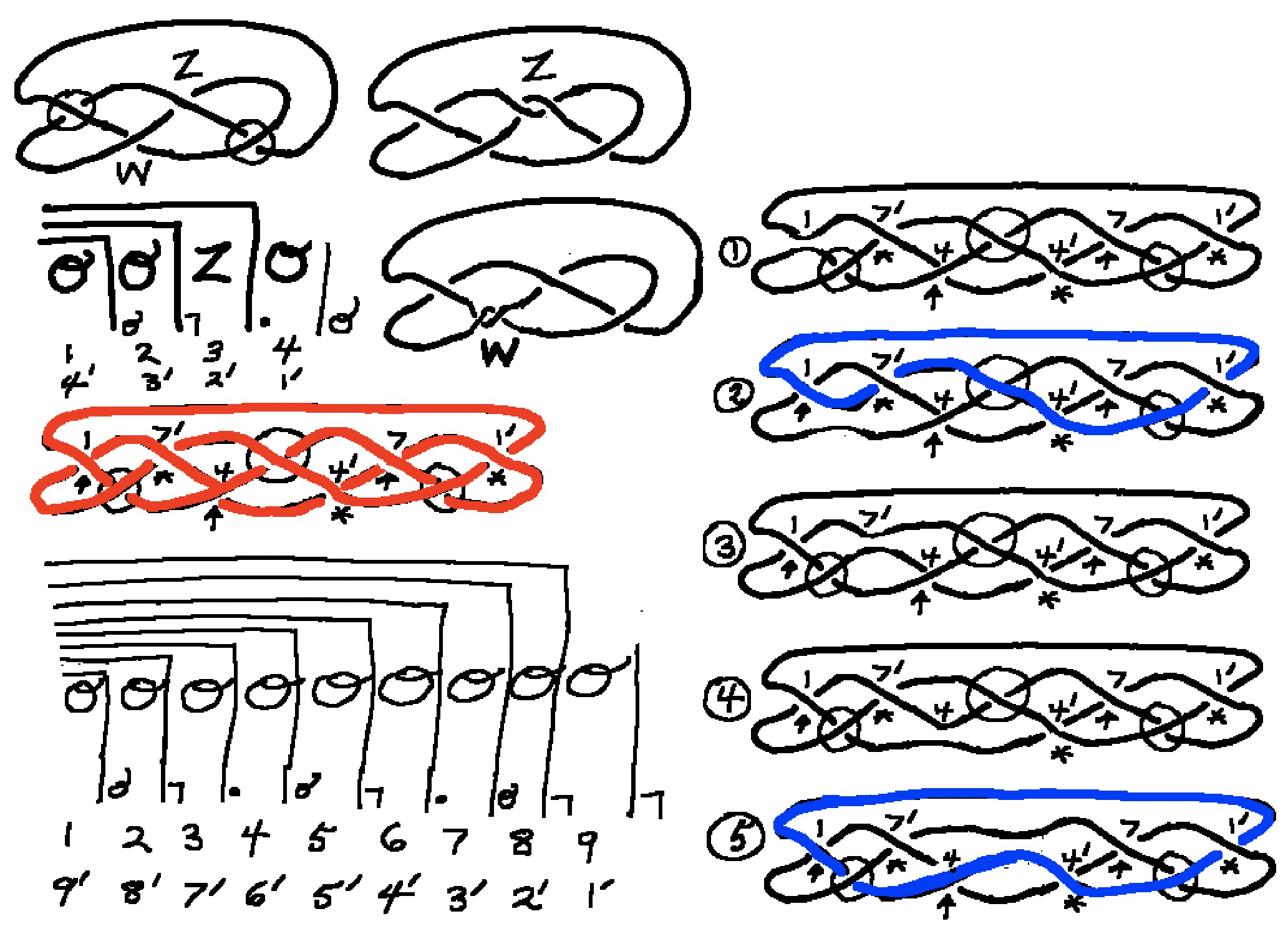}
     \end{tabular}
     \caption{\bf Example showing opacity and transparency.}
     \label{twist}
     \end{center}
     \end{figure}

     \begin{figure}[htb]
     \begin{center}
     \begin{tabular}{c}
     \includegraphics[width=6cm]{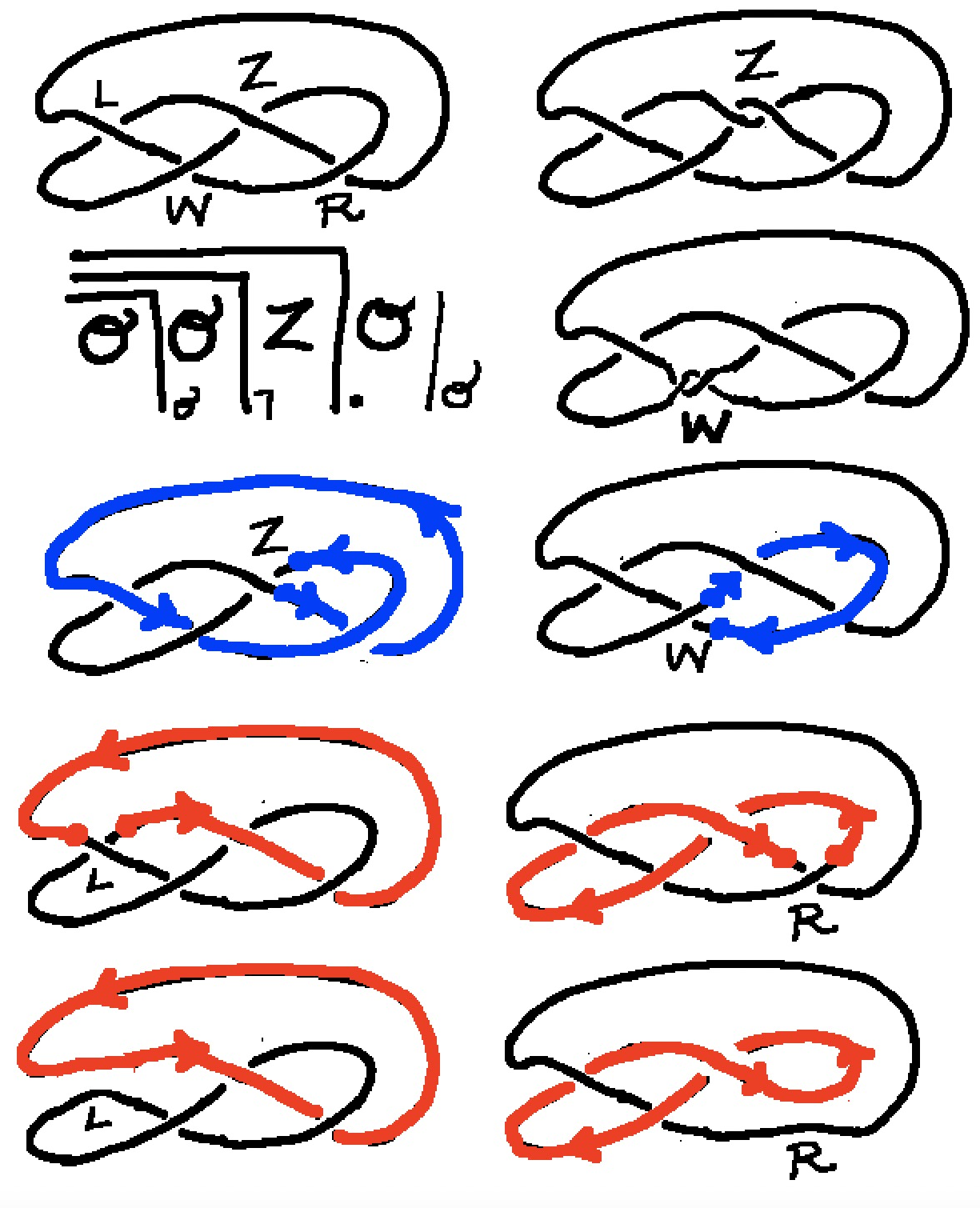}
     \end{tabular}
     \caption{\bf Opacity, transparency and loop structure at a crossing.}
     \label{optrans}
     \end{center}
     \end{figure}

     \begin{figure}[htb]
     \begin{center}
     \begin{tabular}{c}
     \includegraphics[width=6cm]{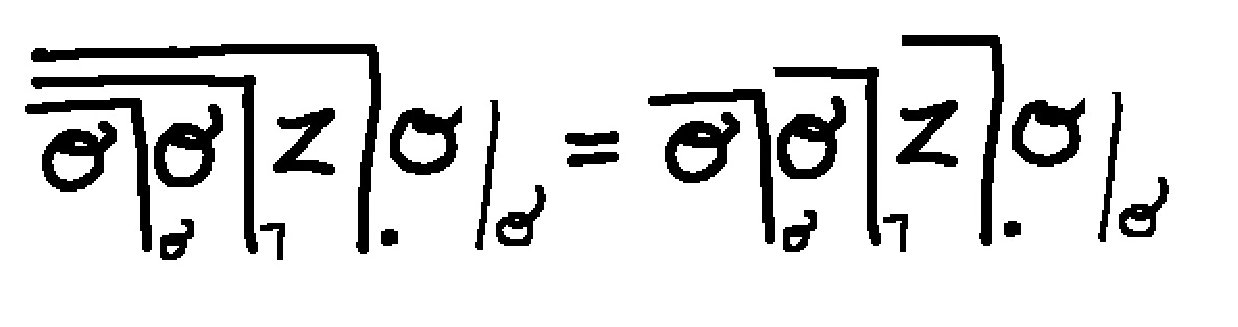}
     \end{tabular}
     \caption{\bf Abbreviated Expression for Sequential Calculation.}
     \label{abbrev}
     \end{center}
     \end{figure}

\noindent {\bf Example.} Consider 
 $$L = \cross{\cross{\cross{\cross{O}_{o}O}_{m}Z}_{u}O\,}_{o}\, O \, |_{u}$$
 Where $Z$ is either $O$ or $E.$ The value of $L$ is independent of the value of $Z$ since $\cross{\,\,}Z = \cross{\,\,}$ in either case.
 Thus $$\cross{\cross{\cross{\cross{O}O}Z}O}O$$ will have two components whether $Z$ is even or odd.
 We say that the value of $L$ is {\it opaque to transmission from $Z$}. We see that {\it opacity to transmission from $Z$ means that Z is a self-crossing!}.
 See Figure~\ref{twist} for an illustration of the weave corresponding to this expression.\\
 
The twist Z is a twist of one component with itself. Removing or adding a crossing to such a twist does not change the 
component count. Conversely, if the value of the expression can be changed by replacing an O by an E or an E by an O, then the corresponding twist (in the case of two components) must be
a twist between two distinct components.  Such changes will change the component count. Thus in the case of rational tangles we can locate, from the crossing algebra, the self-crossing twists as well as the number of components. Opacity to transmission in the value of an expression means that the
 twist at that location in the corresponding knot or link can be changed from even to odd or from odd to even without affecting the component count of the link.  \\
  
 \noindent {\bf Example.} Here is an opacity - transparency example for a knot. Let $K = \cross{\cross{\cross{O}O} Z }O,$ where $Z$ denotes a twist that is either even or odd.  Thus we have
 $$\cross{\cross{\cross{O}_{o}O}_{m}Z }_{u}O|_{o},$$ showing that $K$ is a knot and that its component count is opaque to transmission from $Z.$ If $Z$ is either even or odd, the resulting diagram will have one component. See Figure~\ref{twist} for an illustration for this example with specific choices for the twists. But there is more to say. The rational knot $(1,1,Z,1)$ depicted in Figure~\ref{twist}  can be regarded just as well as 
 $(1,Z,1,1)$ since $(a,b,c,d)$ and $(d,c,b,a)$ represent identical closures as we illustrated in Figure~\ref{order}. So if we start with $(1,W,Z,1)$ corresponding to $\cross{\cross{\cross{O}W} Z }O,$ we conclude that the resulting
 closure will be a knot when {\it either} $W$ or $Z$ is even. Reading the crossings $1,2,3,4$ from left to right, if we start with all odd crossings, then the weave will be a knot. If we change either $3$ or $4$ to even parity, then the weave remains a knot. Links will occur if we change the parity of $1$ or $4.$ \\
 
 In Figure~\ref{optrans} we illustrate the meaning of opacity and transparency in terms of topology of the link diagram in relation to a given crossing. If one chooses an edge at a crossing in a link diagram and walks along the diagram until one returns to that crossing (without going through the crossing during the walk), then a second edge of the crossing is chosen in relation to the initial edge. If the crossing is a part of a twist in a rational link diagram, then opacity means that smoothing the crossing in the twist direction (i.e. changing the parity of the twist) will not change the connectivity of the link. As the figure shows, this is what we detect algebraically when we detect opacity or transparency. It is of interest that the algebra can see this aspect of the diagram topology.\\
 
 The next example in  Figure~\ref{twist} shows a Fibonacci knot $K_9 $with $9$ crossings. The same argument with opacity shows that 
 crossings $\{ 3,6,9 \}$ can be made even without changing to a link by counting from the left. See the cross algebra expression in the figure showing opacity at these positions. Then by symmetry, we can count from right to left
 and conclude that crossings $\{3',6',9'\} = \{7,4,1\}$  are also opaque to transmission. Thus altogether we conclude that positions $\{1,3,4,6,7,9\}$ are opaque  and so can be individually switched to even and retain a single component weave. The places  $\{2,5,8\}$ when changes from odd to even produce links, as illustrated in the figure. We call these positions {\it transparent}.\\
 
  We can generalise this statement by considering Fibonacci knots $K_{n}.$ By our previous discussion $K_{n}$ is a Fibonacci knot when $n = 3k$ or $n=3k+1.$ In Figure~\ref{twist} we analyzed $k=3.$ Call a twist location $i$ for a twist $a_{i}$  of a rational knot  $K = (a_1,a_2,\cdots , a_n)$ {\it opaque} if changing it from odd to even or from even to odd does not change $K$ from being a knot to being a link.
  The reader will have no difficulty using our technique to show that the {\it transparent} crossings of $K_{3k}$  are $2,5,8,\cdots, (3k-1)$ and that the {\it transparent} crossings of $K_{3k+1}$ are $1,4,7,\cdots, (3k+1).$\\

  \noindent{\bf Example.} We are now in position to enumerate rational knots and links with $N$ crossings and to discriminate which are knots and which are links, by pure algebra and combinatorics. For a given value of $N,$ list all ordered partitions of $N$ with no $1$ at either end of the list and all positive entries in the partition. Call two such partitions equivalent if
 one is the reversed order of the the other. The equivalence classes are in $1-1$ correspondence with the rational knots and links with $N$ crossings by our remarks in the previous section about the classification of rational knots and links. Make a list of the equivalence classes. For each element in the list, use the crossing algebra to test whether it has one component or two components. Separate the list into a list of knots and links.\\
 
\begin{figure}[htb]
     \begin{center}
     \begin{tabular}{c}
     \includegraphics[width=8cm]{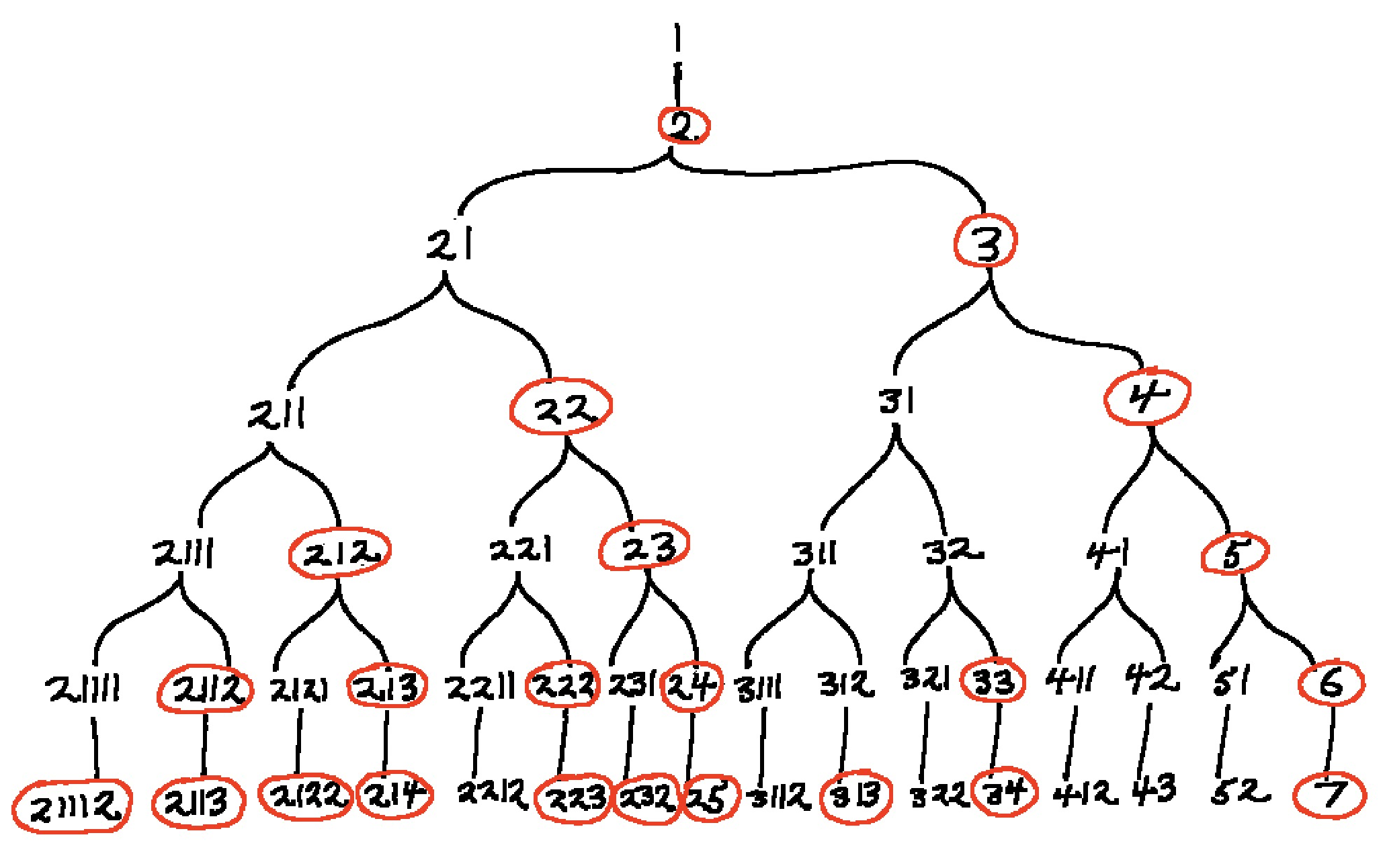}
     \end{tabular}
     \caption{\bf Ordered Partitions for Rational Knots and Links Up to Seven Crossings.}
     \label{part1}
     \end{center}
     \end{figure}

\begin{figure}[htb]
     \begin{center}
     \begin{tabular}{c}
     \includegraphics[width=8cm]{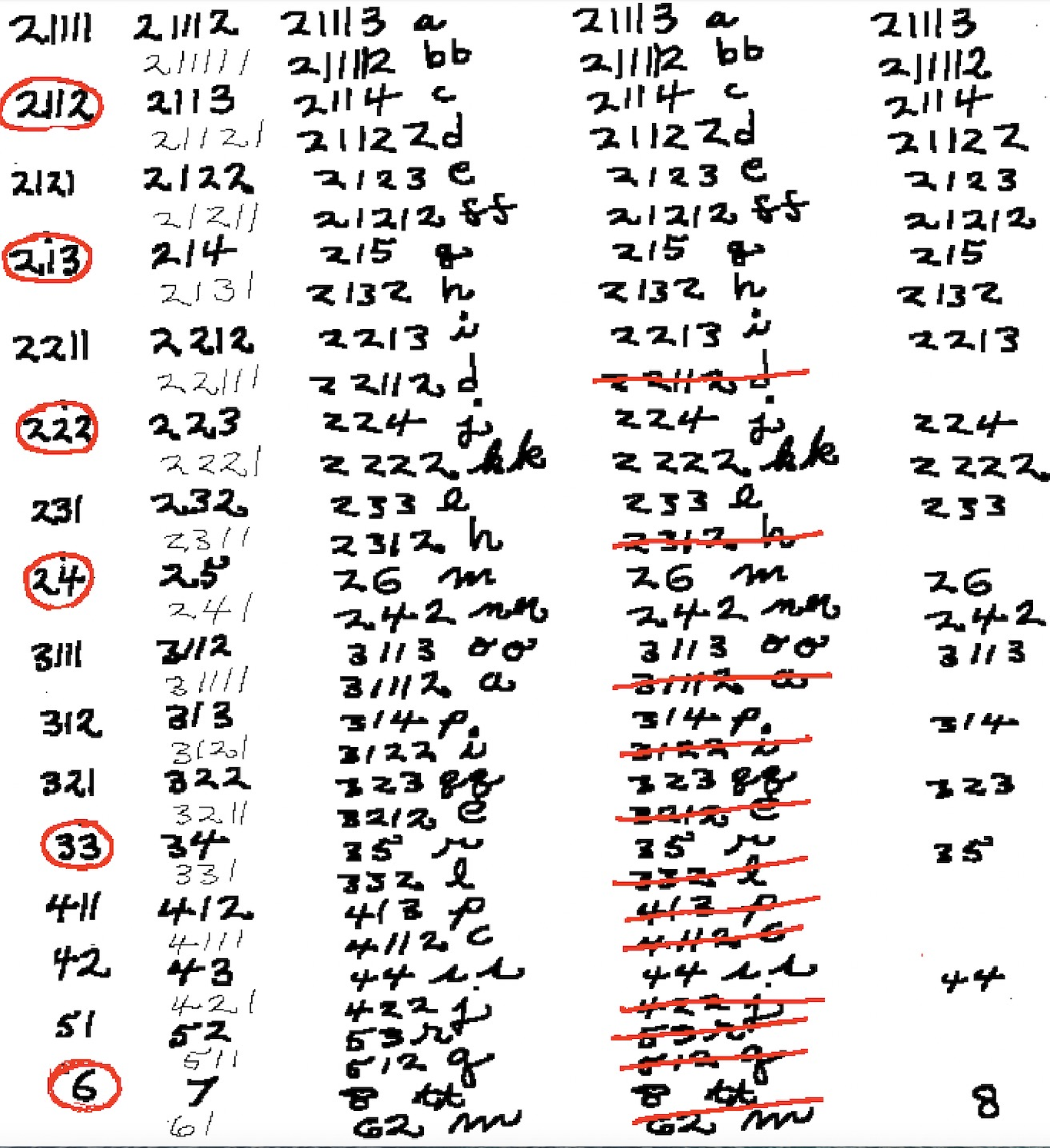}
     \end{tabular}
     \caption{\bf Extension for 8 Crossings.}
     \label{part2}
     \end{center}
     \end{figure}

\begin{figure}[htb]
     \begin{center}
     \begin{tabular}{c}
     \includegraphics[width=8cm]{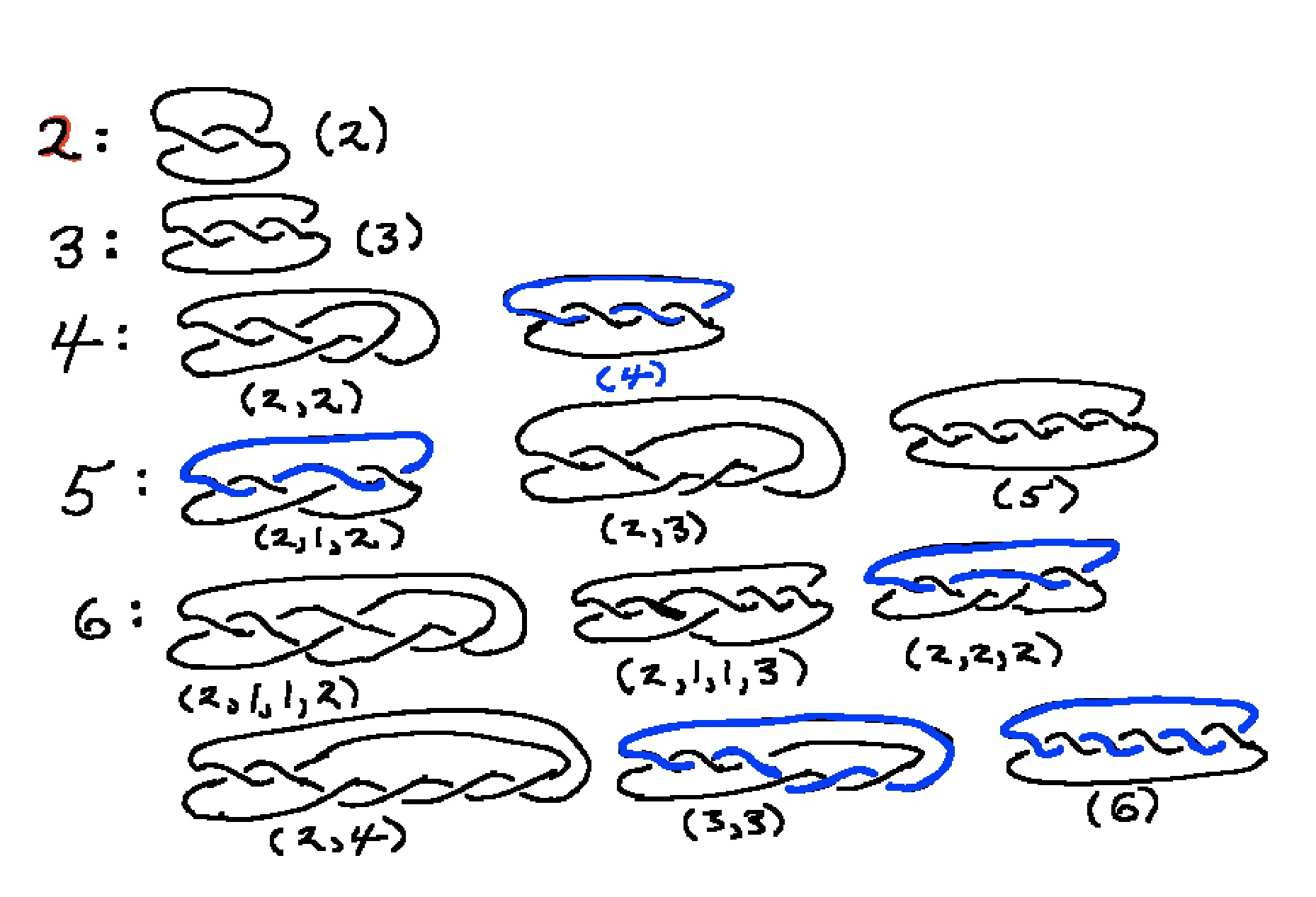}
     \end{tabular}
     \caption{\bf Rational knots and links with two to six crossings, from partitions .}
     \label{twosix}
     \end{center}
     \end{figure}

\begin{figure}[htb]
     \begin{center}
     \begin{tabular}{c}
     \includegraphics[width=8cm]{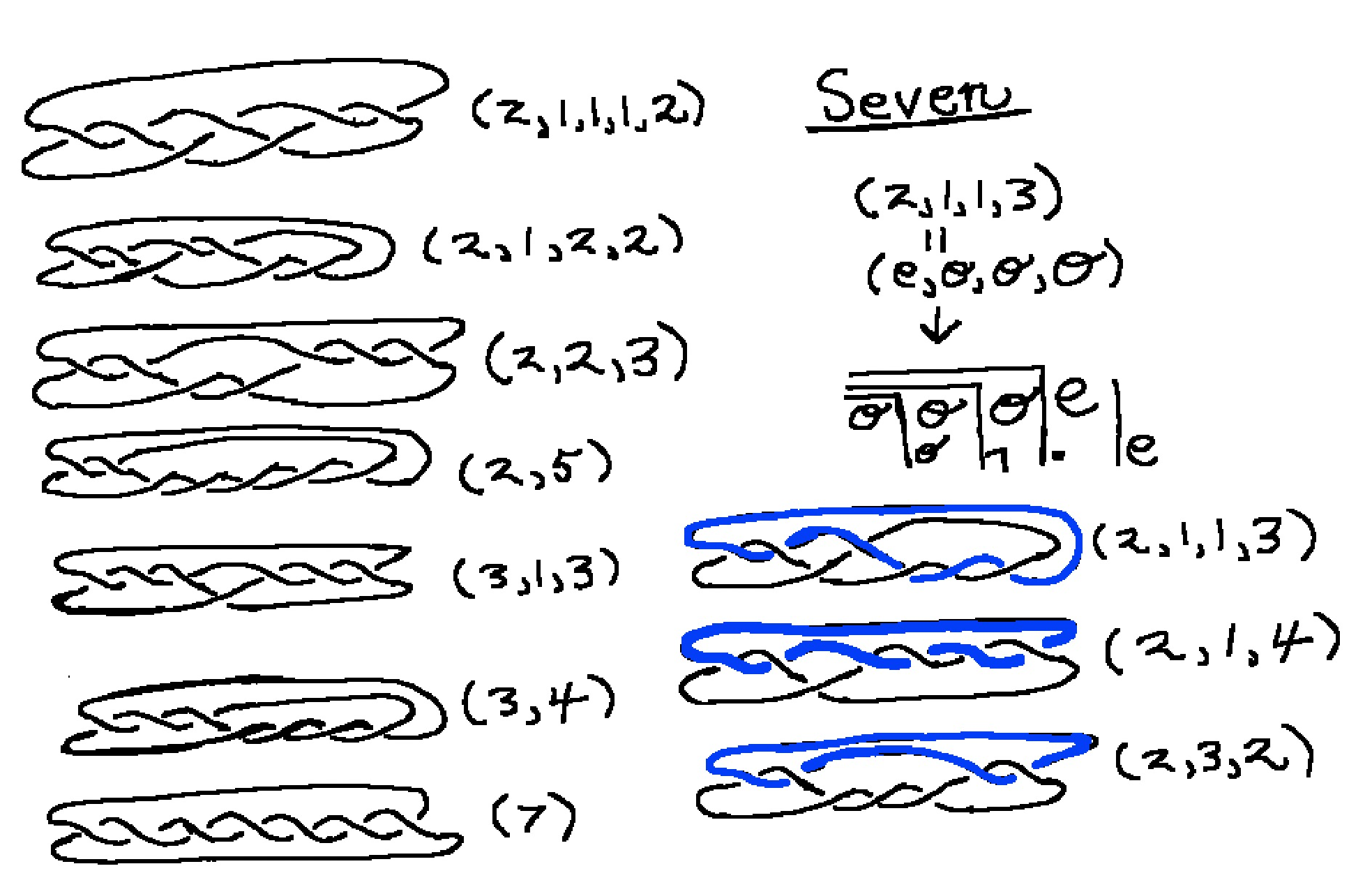}
     \end{tabular}
     \caption{\bf Rational links with seven crossings, from partitions.}
     \label{seven}
     \end{center}
     \end{figure}

\begin{figure}[htb]
     \begin{center}
     \begin{tabular}{c}
     \includegraphics[width=10cm]{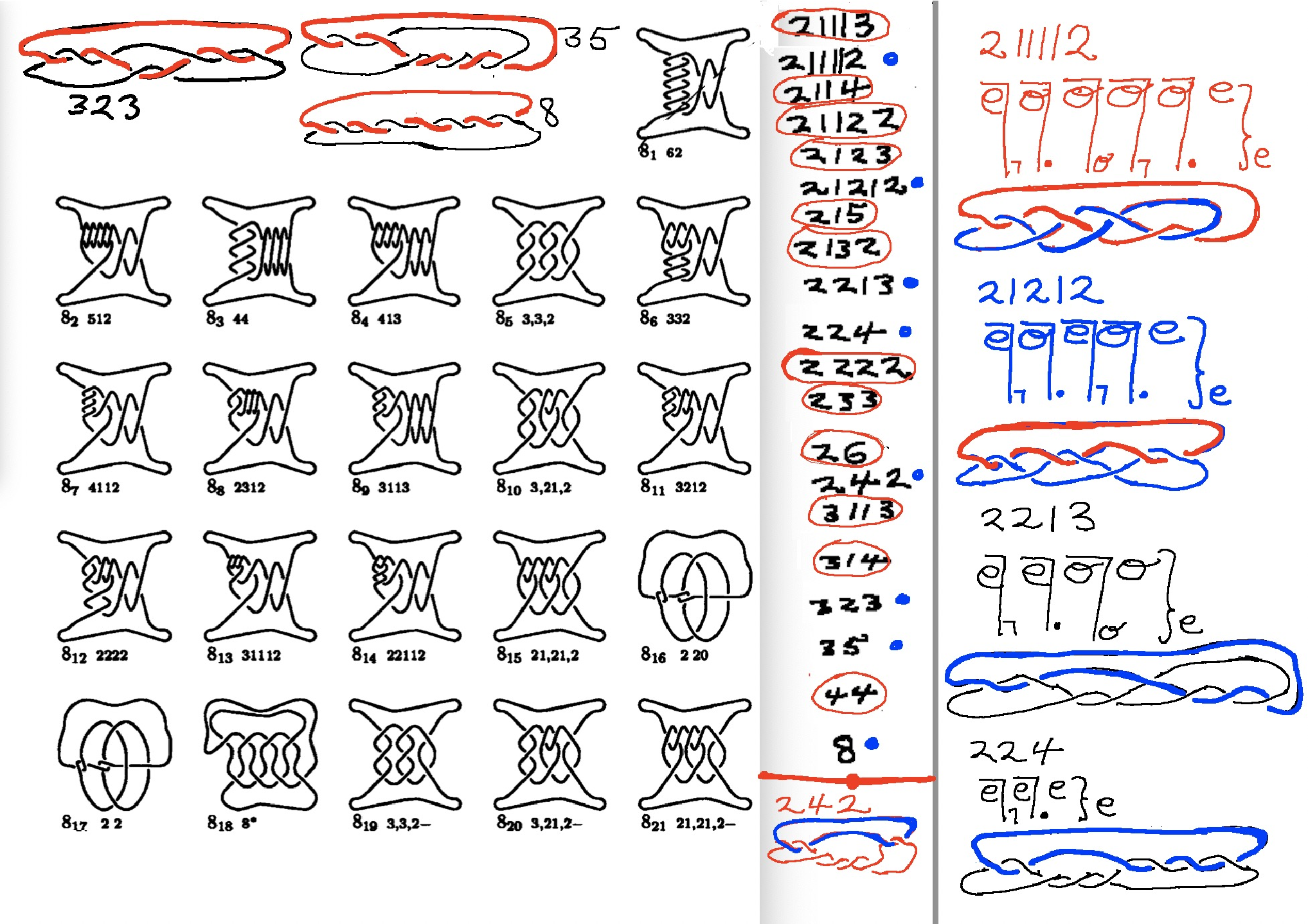}
     \end{tabular}
     \caption{\bf Table of 8 crossing knots and rational links.}
     \label{eight}
     \end{center}
     \end{figure}

 Here is a specific example for this procedure. Let $N=7.$ The knot list is $$(7), (2,5),(3,4),(2,2,3),(3,1,3),(2,1,2,2),(2,1,1,1,2)$$ and the links are the list
 $$(2,3,2),(2,1,4),(2,1,1,3).$$ This example can be done by hand. In fact, we have illustrated a method in Figure~\ref{part1}, Figure~\ref{part2}, Figure~\ref{twosix}, Figure~\ref{seven}, Figure~\ref{eight}. The first figure illustrates a method to produce all the needed ordered partitions ( by either adding one to end of the list, or adding one to the last member of the list). In the first figure we do not need the branch of this process that always has a 1 as the left-most element, and when we arrive at the row for seven crossings, we do not write the partitions with a 1 at the right end. Thus the last row contains all the knots and links for seven crossings, with repetitions of the sort $(4,3) \sim (3,4)$ since order reversal gives the same link. The crossing calculus can be used to discriminate knots from links.  In Figure~\ref{part2} we contine the process from six to eight by including all the descendants of row six, making the single end additions that create row eight from them and then culling out just one representive for each rational link with eight crossings. In Figure~\ref{eight} we use the culled list of rational links and make an extended table using the table of 8 crossing knots from the book by Kawauchi \cite{Kawauchi}. Kawauchi's table lists all the knots with eight crossings. Rational knots are indicated and drawn in continued fraction form. We expand the table to include all the rational links of eight crossings. The figure includes drawings of all the eight crossing knots and all the eight crossing rational links. It also includes some sample crossing algebra computations. Note that we have adopted an abbreviated notation for calculating in the form shown below. The right-hand side is meant to be short hand for the left-hand side, and writing this way saves the nesting of the boxes.
 $$\cross{\cross{\cross{E}_{m}E}_{u}O }_{o}O\,|_{u} = \cross{E}_{m}\cross{E}_{u}\cross{O }_{o}O\,|_{u} $$

 The reader will appreciate our methods if she tries enumerating all the rational knots and links with (say)  $9$ or more crossings, possibly using a computer to produce the partitions and a crossing algebra program to determine which are knots and which are links. See the next section of the present paper for a discussion of programming the crossing algebra.\\

 \subsection {\bf Arborescent Links.}
 The formalism of the crossing algebra allows us to determine connectivity count for the larger class of {\it aborescent links}.
 From the point of view of our formalism, an aborescent link corresponds to any expression in positive integers that is generated by 
 addition of tangles ($x + y$) and inversion of tangles ($1/x$). Thus $K = Num[1/3 + 1/2]$ corresonds to the numerator of the non-rational tangle sum of two vertical twists of type $3$ and type $2.$ In the cross formalism we have
 $$\cross{O} \cross{E} = O \cross{\,\,} = \cross{\,\,}$$ and hence (as we knew) $K$ is a knot.
 
 We can indicate an aborescent link as a generalised continued fraction. For example,
 
  $$L = a_{1} + \frac{1}{    \frac{1}{a_{2} + \frac{1}{a_{3}}  } + \frac{1}{a_{4} + \frac{1}{a_{5} + \frac{1}{a_{6}} }  }    }$$ 
  $$= \cross{ \cross{ \cross{a_{3}}a_{2} }   \cross{ \cross{ \cross{a_{6}} a_{5} } a_{4} }  } a_{1}.$$

 On the other hand, note that $$\cross{\,\,}\cross{\,\,} = \cross{\Asmooth} \cross{\Asmooth} = \Bsmooth \Bsmooth$$
 and the sum of two $\Bsmooth$ tangles results in a circle component. Thus we should write
 $$\Bsmooth \Bsmooth = \delta \Bsmooth$$ where $\delta$ is an algebraic variable corresponding to a loop.
 Then correspondingly we would have $$\cross{\,\,}\cross{\,\,}  = \delta \cross{\,\,}$$ connoting a link of two components.
 
 From this we can use our calculus to determine the number of components in an arborescent link, given as such a tree structure.
 For example
 $$S = \cross{ \cross{ \cross{O}O }   \cross{ \cross{ \cross{O}O}O}  } O$$
 $$=\cross{ \cross{OO }   \cross{ \cross{ OO}O}  } O$$
 $$=\cross{ \cross{\,\, }   \cross{ \cross{ \,\,}O}  } O$$
 $$=\cross{ \cross{\,\, }   \cross{ \cross{ \,\,}}  } O$$
 $$=\cross{ \cross{\,\,} } O$$
 $$=O$$
 
 Hence if all the twists in $L$ above are odd, then the structure is a knot of one component.
 We see from this calculation that it will be a link if the rightmost $O$ is changed to $E.$ And, in fact we can do the transmission analysis on 
 the structure:
 
$$S = \cross{ \cross{ \cross{a_{3}}a_{2} }   \cross{ \cross{ \cross{a_{6}} a_{5} } a_{4} }  } a_{1} = \cross{ \cross{ \cross{O}_{o}O }_{m}   \cross{ \cross{ \cross{O}_{o}O}_{m}O}_{u}  }_{u}O$$

From this we deduce that if $a_{3}$ and $a_{2}$ are odd, then the connectivity of $S$ will be unaffected by any assignments to
$a_{4},a_{5},a_{6}$ but changing $a_{1}$ will toggle $S$ back and forth from knot to link.\\

\noindent {\bf Completing the Calculus.} To complete this component count calculus we can note that a single appearance of $E$ can be replaced by $\delta^{2}$ and a single appearance of $\cross{\,\,}$ can be replaced by 
a single $\delta$ since $\cross{\,\,} = \cross{\Asmooth} = \Bsmooth$ is a tangle whose numerator closure is one loop. Similarly, we have that $O = \Fcross$ with closure a single loop. Thus, we add to the crossing calculus these rules to be applied to reduced expressions.
$$E \longrightarrow \delta^{2},$$
$$O \longrightarrow \delta,$$
$$\cross{\,\,} \longrightarrow \delta.$$

For example $$B = \cross{ \cross{2} \cross{-2}}  \cross{2} \cross{-2}.$$ This arborescent link is the Borommean Rings. See Figure~\ref{boro} for an illustration of this formalism for the Borommean rings and two other arborescent examples.
To see that the rings have three components, we calculate
$$\cross{ \cross{E} \cross{E}}  \cross{E} \cross{E}$$
$$=\cross{ \cross{\,\,} \cross{\,\,}}  \cross{\,\,} \cross{\,\,}$$
$$=\cross{ \cross{\,\,} \delta}  \cross{\,\,} \delta$$
$$=\cross{ \cross{\,\,}}  \cross{\,\,} \delta^{2}$$
$$=\cross{\,\,} \delta^{2}$$
$$=\delta^{3}.$$

\begin{figure}[htb]
     \begin{center}
     \begin{tabular}{c}
     \includegraphics[width=12cm]{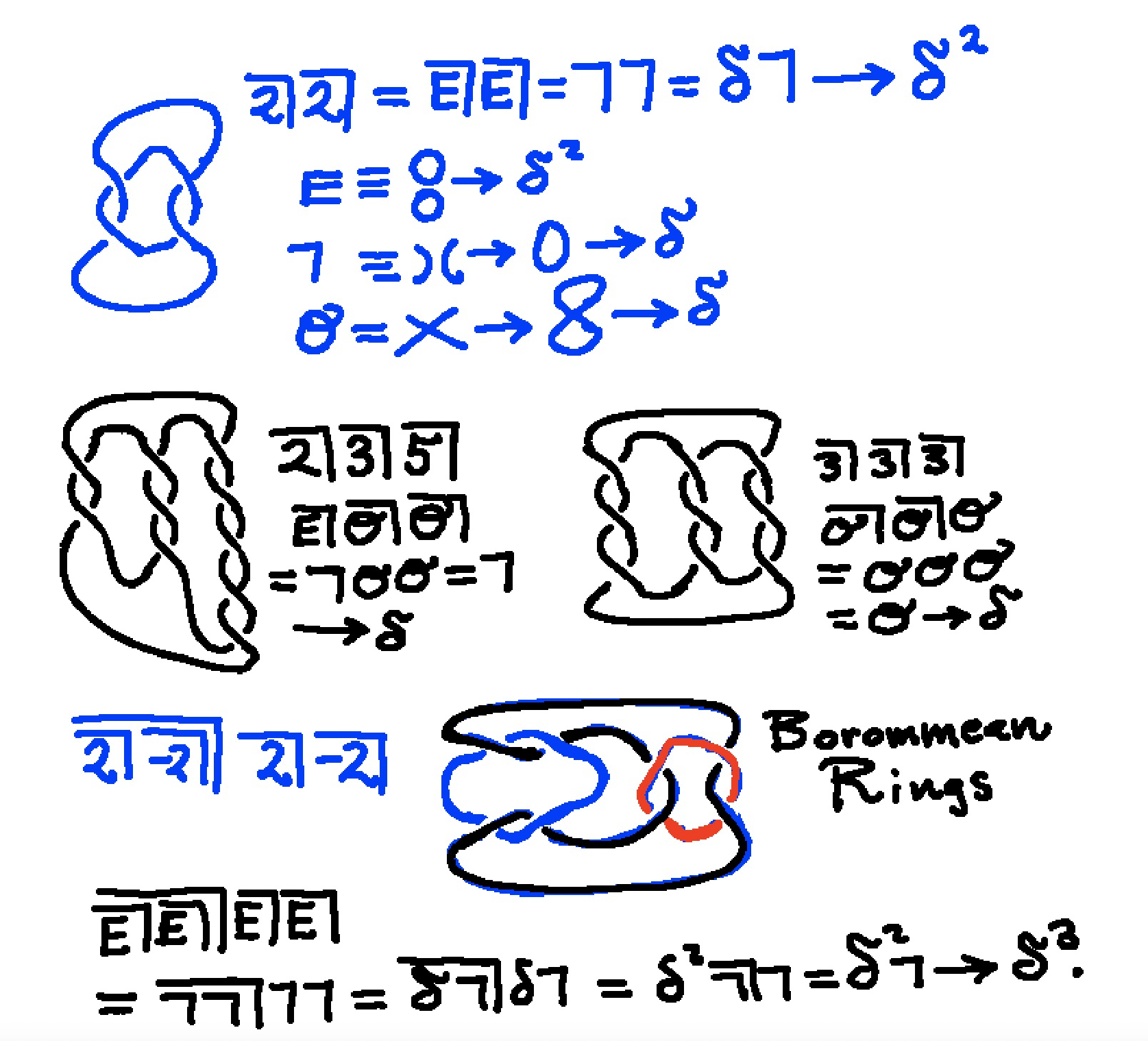}
     \end{tabular}
     \caption{\bf Pretzel Knots and Links and Borommean Rings}
     \label{boro}
\end{center}
\end{figure}

\begin{figure}[htb]
     \begin{center}
     \begin{tabular}{c}
     \includegraphics[width=12cm]{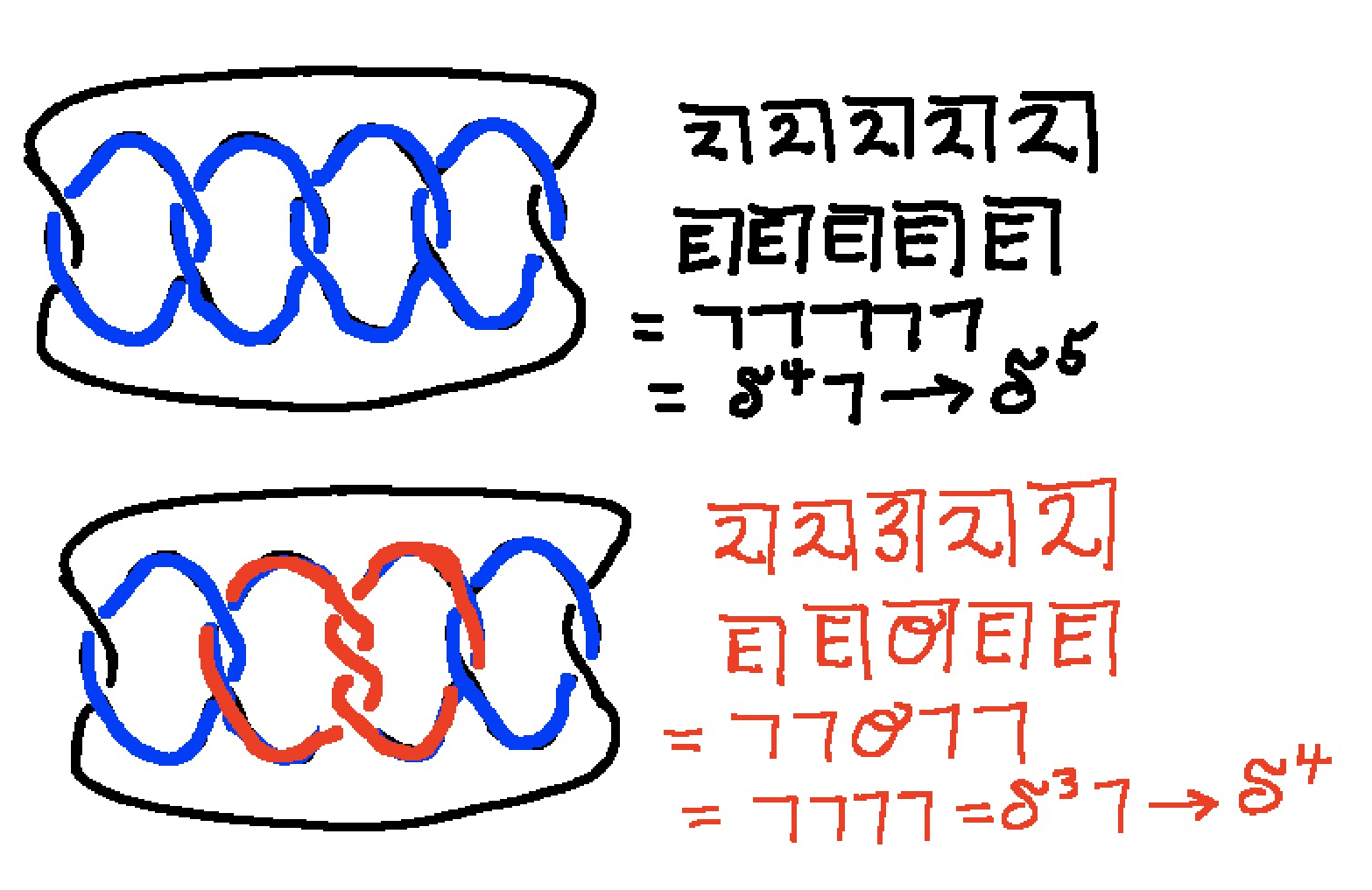}
     \end{tabular}
     \caption{\bf Pretzel knots and links}
     \label{pretzel}
\end{center}
\end{figure}

\begin{figure}[htb]
     \begin{center}
     \begin{tabular}{c}
     \includegraphics[width=12cm]{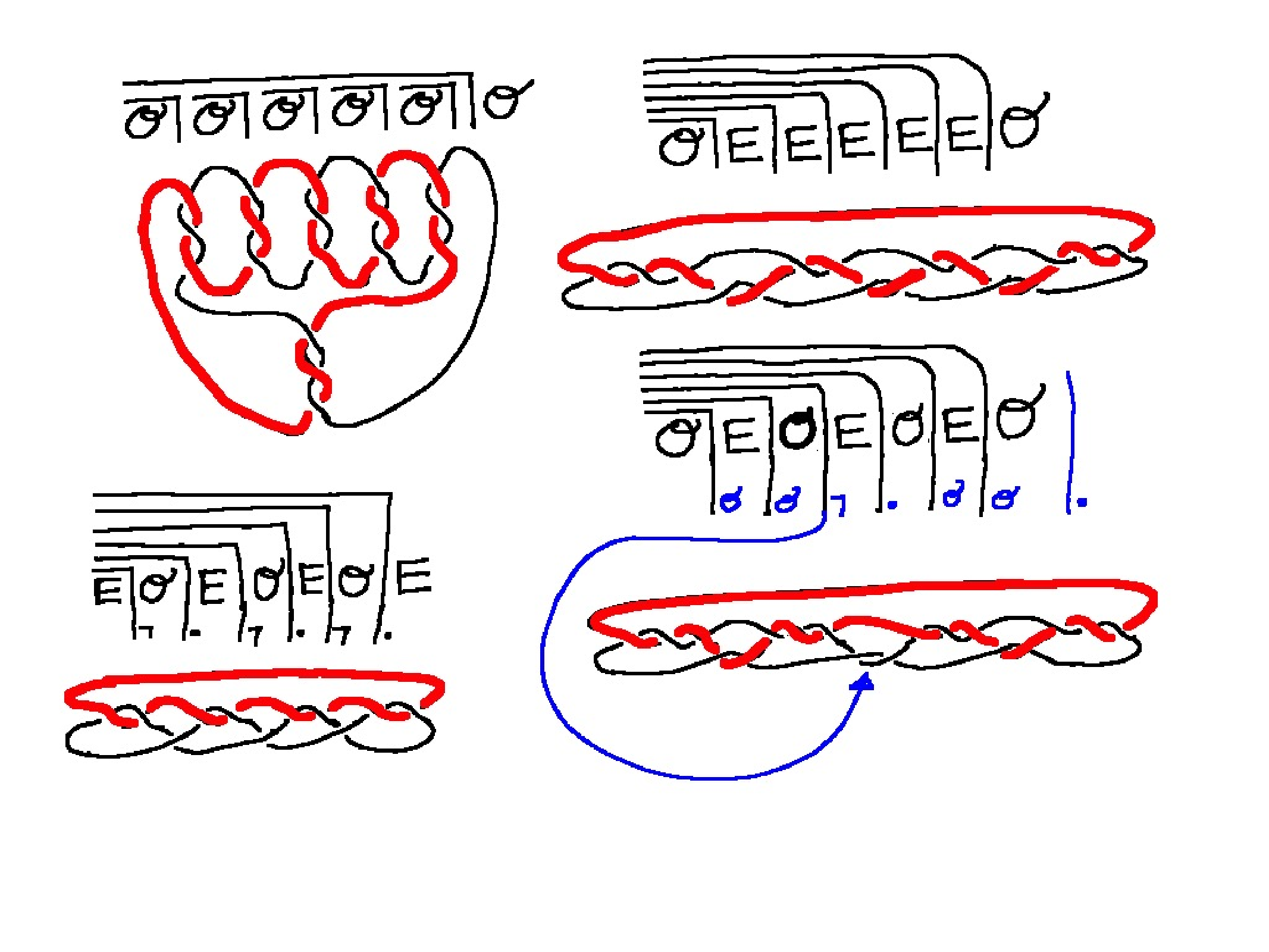}
     \end{tabular}
     \caption{Linking patterns}
     \label{linking}
\end{center}
\end{figure}

\begin{figure}[htb]
     \begin{center}
     \begin{tabular}{c}
     \includegraphics[width=12cm]{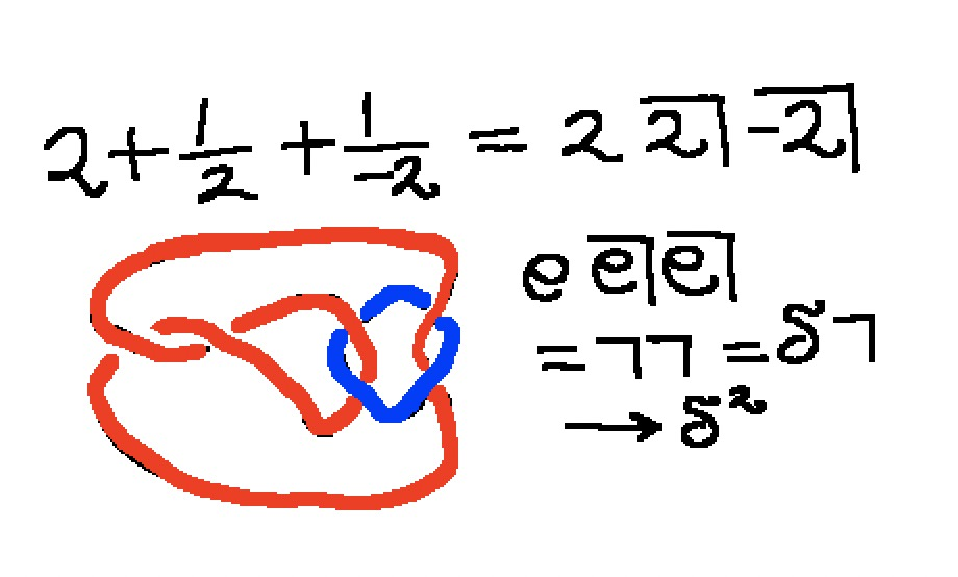}
     \end{tabular}
     \caption{\bf Whitehead Link}
     \label{arbor}
\end{center}
\end{figure}

\begin{figure}[htb]
     \begin{center}
     \begin{tabular}{c}
     \includegraphics[width=12cm]{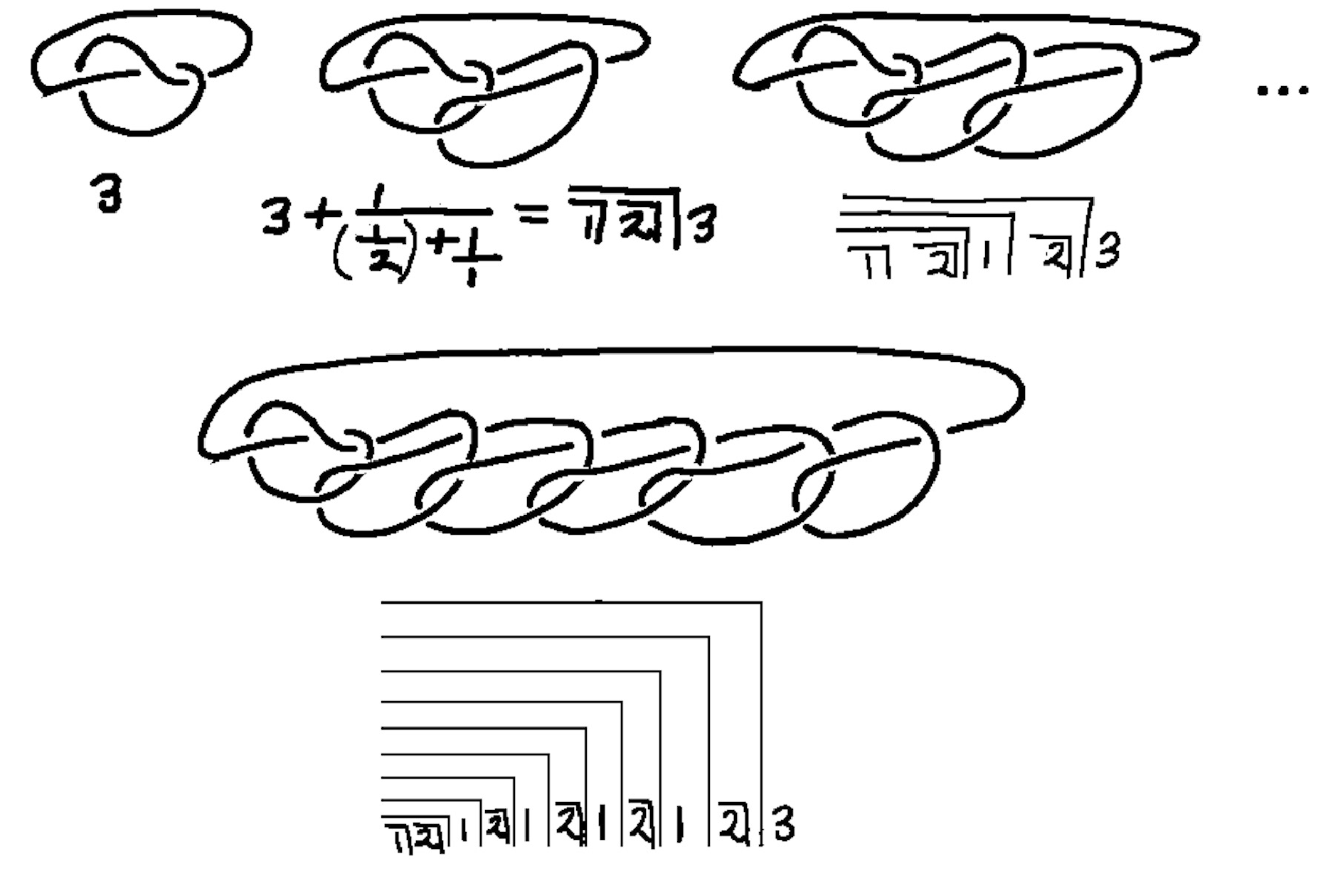}
     \end{tabular}
     \caption{\bf chain}
     \label{chain}
\end{center}
\end{figure}

\begin{figure}[htb]
     \begin{center}
     \begin{tabular}{c}
     \includegraphics[width=12cm]{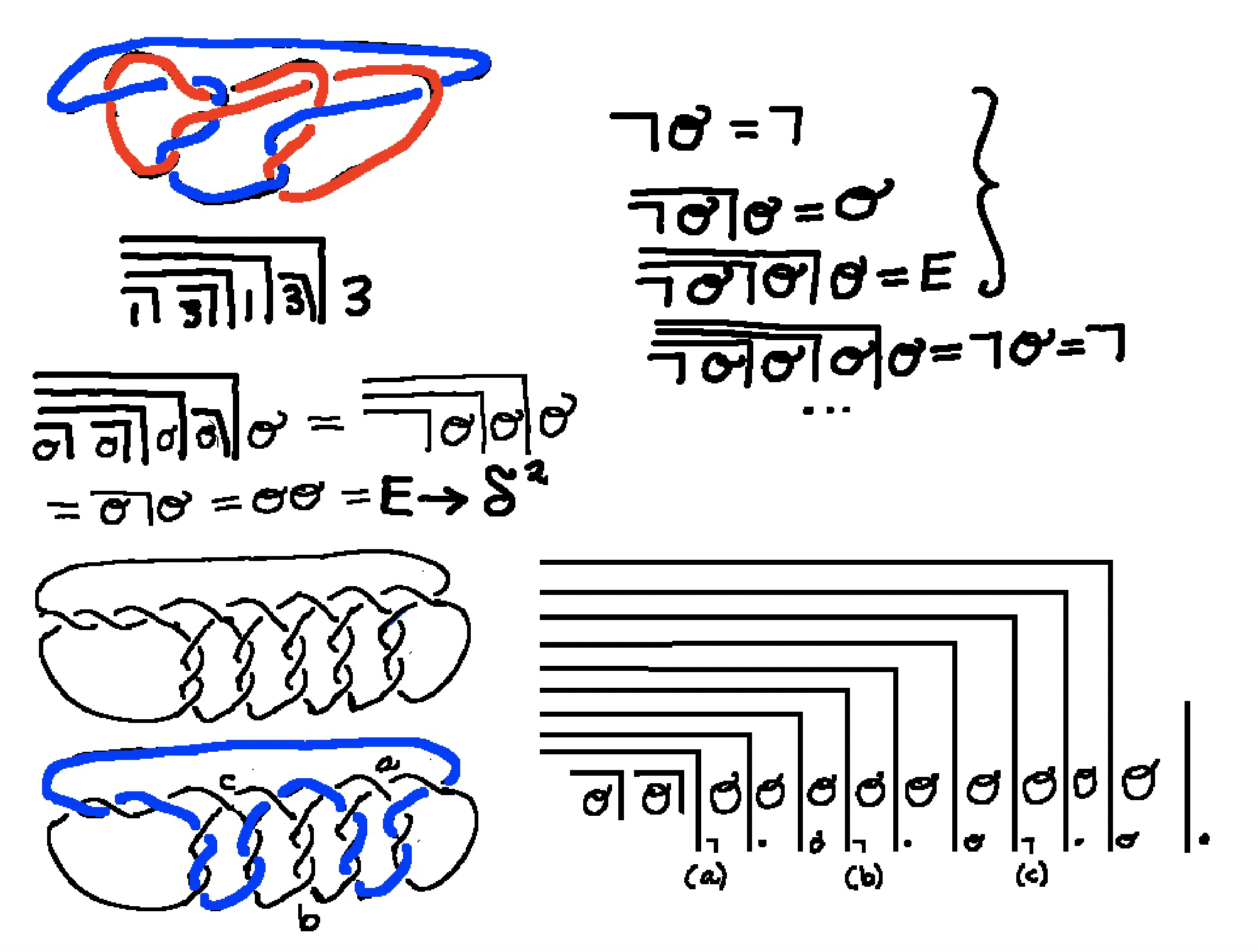}
     \end{tabular}
     \caption{\bf Stitching}
     \label{stitch}
\end{center}
\end{figure}

 \begin{figure}[htb]
     \begin{center}
     \begin{tabular}{c}
     \includegraphics[width=12cm]{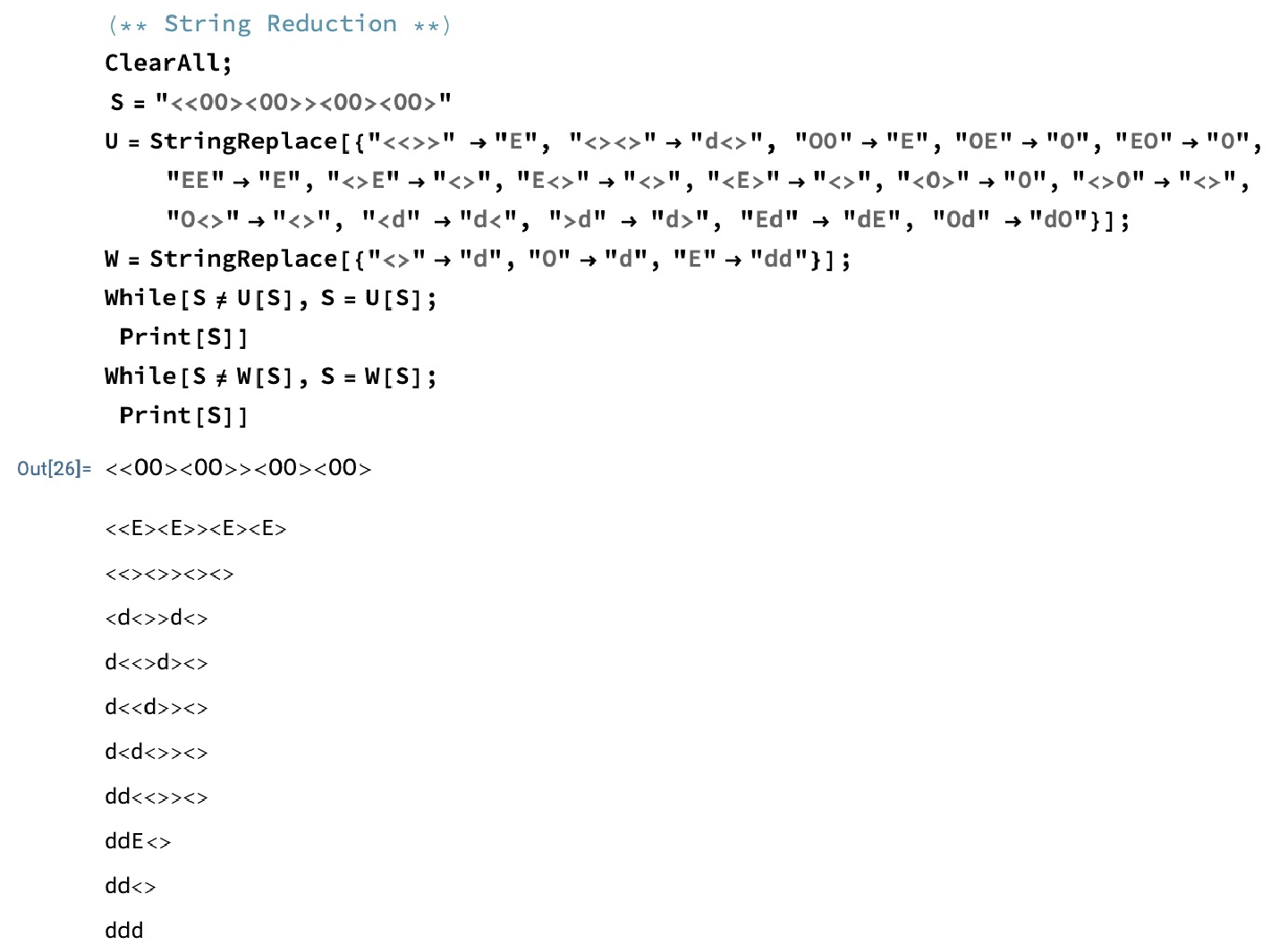}
     \end{tabular}
     \caption{\bf Crossing Algebra Reduction for Component Count by String Manipulation}
     \label{strings}
\end{center}
\end{figure}

\noindent{\bf Examples.} We end this section with a selection of examples that illustrate the principles and properties of the crossing algebra in relation to knots and links.
\begin{enumerate}
\item In Figure~\ref{pretzel} we illustrate more preztel knots and links.
\item In Figure~\ref{linking} we illustrate how linking patterns arise in relation to the algebra.
 \item In Figure~\ref{arbor} illustrates the Whitehead Link $W$(non trivial but with zero linking number) in the form $2 \cross{2}\cross{-2}.$ Note that 
$E \cross{E} \cross{E} = \cross{\,\,}\cross{\,\,} = \delta \cross{\,\,} \longrightarrow \delta^{2},$ showing that $W$ has two components.
\item In Figure~\ref{chain} we illustrate the crossing algebra expressions for a chain stitch.
\item In Figure~\ref{stitch} we illustrate properties of a chain stitch with different parity and multiple components.
\end{enumerate}

\noindent {\bf Programming.} 
This counting method can be written as a string reduction program in $Mathematica^{TM}$ as shown in Figure~\ref{strings}.
In this figure, the operator $U$ does the first crossing algebra reduction. Then the operator $W$ performs the completion that we have just discussed.
The computation in this figure is an automatic version of our hand computation for the Borommean Rings above. In the program we write $\langle A \rangle = \cross{A}$ and $d$ for $\delta,$ but otherwise use the 
same notation as in the present text.
In the next section, we explain how to use the 
crossing algebra and string reduction techniques to compute the Kauffman bracket polynomial.\\

\section{Bracketology}

The Kauffman bracket polynomial model for the Jones polynomial \cite{K,KNew,KOK,KD,TL,KS,KP,KD} is usually described by the expansion
$$
\langle \Across \rangle=A \langle \Asmooth \rangle + A^{-1}\langle
\Bsmooth \rangle \label{kabr10}.
$$

and we have

\begin{equation}
\langle K \, \bigcirc \rangle=(-A^{2} -A^{-2}) \langle K \rangle \label{kabr1}
\end{equation}

\begin{equation}
\langle \Rcurl \rangle=(-A^{3}) \langle \Arc \rangle \label{kabr2}
\end{equation}

\begin{equation}
\langle \Lcurl \rangle=(-A^{-3}) \langle \Arc \rangle \label{kabr3}
\end{equation}
\bigbreak

\begin{equation}
\langle O \rangle= 1 \label{kabr4}
\end{equation}
\bigbreak

The bracket is a Laurent polynomial in $A$ and $A^{-1}$ with integer coefficients and it is an invariant of regular isotopy of knots and links
(Reidemeister moves 2 and 3). By a normalization, we obtain a polynomial $f_{K}(A) = (-A^{3})^{-wr(K)}\langle K \rangle$ that is invariant under all three Reidemeister moves, for oriented links $K.$ The function $wr(K)$, the {\it writhe} of $K$, is the sum of the crossing signs in the diagram $K.$ See \cite{KNew} for more information about the conventions for the bracket polynonmial and basic theorems about the invariant.\\

The bracket polynomial can be expressed as a state summation in the form $$\langle K \rangle = \sum_{S}\langle K | S \rangle d^{||S||}$$ where this is a sum over {\it states} $S$ of the diagram $K.$ A {\it state} $S$ of a diagram $K$ is a specific choice of smoothing for each crossing together with a label $A$ or $B$ at the smoothing. The skein relation
$
\langle \Across \rangle=A \langle \Asmooth \rangle + B\langle
\Bsmooth \rangle \label{kabr11}
$
with $B$ instead of $A^{-1}$ shows the two smoothings at a crossing and the labels $A$ and $B$. We can form a three variable bracket in the form of the state sum by using these labels and the variable $d$ for the loop value. In the above formulas we have $d=-A^{2} - A^{-2}.$ In the state sum formula $d$ is raised to the power $||S||,$ designating the number of loops in the state $S$.
To obtain the topologically specialized bracket we set $B= A^{-1}, d = -A^{2}-A^{-2}$ and divide the raw polynomial by $d.$\\

For arborescent links we can use the crossing algebra to determine the number of loops in each state, and thus we have an algebraic method to compute the bracket polynomial for this class of 
links. We will give small examples, and indicate the strategy for a computer program that can perform this calculation. To this end, we shall let $O$ designate the specific odd crossing
$$O = \Across$$ and $O^{\star}$ denote its mirror image $$O^{\star} = \Bcross.$$ Then $L = OO$ designates the numerator of the tangle $\Across \Across,$ the Hopf link shown in Figure~\ref{ratl4}. Note that in terms of our crossing algebra notation, the bracket skein expansion has the form 
$$
\langle O \rangle= \langle \Across \rangle =A \langle E \rangle + B\langle
\cross{E} \rangle \label{kabr12}
$$
$$
\langle O^{\star} \rangle = \langle \Bcross \rangle=A \langle \cross{E} \rangle + B\langle
E \rangle \label{kabr13}
$$

We shall denote the states of $L$ in the form shown  below. We write the labels for the state, followed by the expression for the link diagram with the corresponding smoothings $E$ or $\cross{E}$ inserted in place of $O$ or $O^{\star}.$
\begin{enumerate}
\item $AAEE$
\item $ABE\cross{E}$
\item $BA\cross{E}E$
\item $BB\cross{E}\cross{E}$
\end{enumerate}
Each of these expressions reduces to the corresponding state evaluation via the crossing algebra.
\begin{enumerate}
\item $AAEE= A^2 E = A^2 d^2$
\item $ABE \cross{E} = AB \cross{\,\,} = ABd $
\item $BA \cross{E} E = BA \cross{\,\,} = BAd $
\item $BB \cross{E} \cross{E} = BB \cross{\,\,} \cross{\,\,} = BBd \cross{\,\,} = B^2d^2$
\end{enumerate}
Thus $$\sum_{S}\langle K | S \rangle d^{||S||} = A^2 d^2 + 2 ABd + B^2 d^2$$
Dividing by $d$ and taking the above values for $B$ and $d,$ we have
$$\langle L \rangle = A^2 d + 2 AB + B^2 d = (A^2 + A^{-2})(-A^2 - A^{-2}) + 2  = -A^4 - A^{-4},$$
the bracket evaluation of the Hopf link diagram.\\

\begin{figure}[htb]
     \begin{center}
     \begin{tabular}{c}
     \includegraphics[width=11cm]{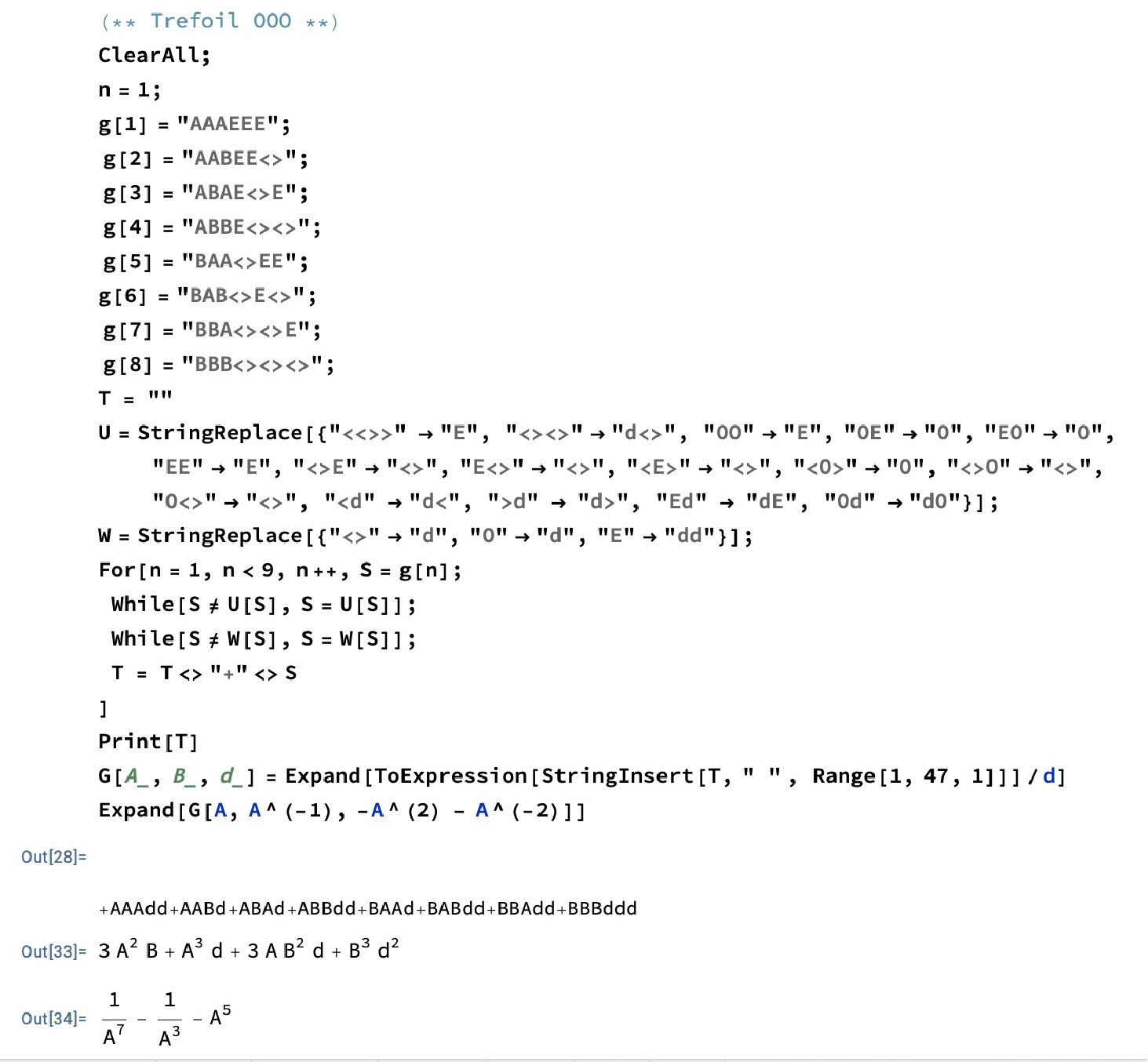}
     \end{tabular}
     \caption{\bf Crossing Algebra Bracket Calculation for the Trefoil Knot $OOO$}
     \label{progtrefoil}
\end{center}
\end{figure}

In Figure~\ref{progtrefoil} we illustrate the corresponding calculation for the trefoil knot, whose aborescent formula is $K=OOO.$
The eight labeled states are listed and then given to a computer program that does the crossing algebra reduction and
below that is the code that produces the bracket polynomial from the raw bracket polynomial expressed in terms of character strings.\\

\begin{figure}[htb]
     \begin{center}
     \begin{tabular}{c}
     \includegraphics[width=10cm]{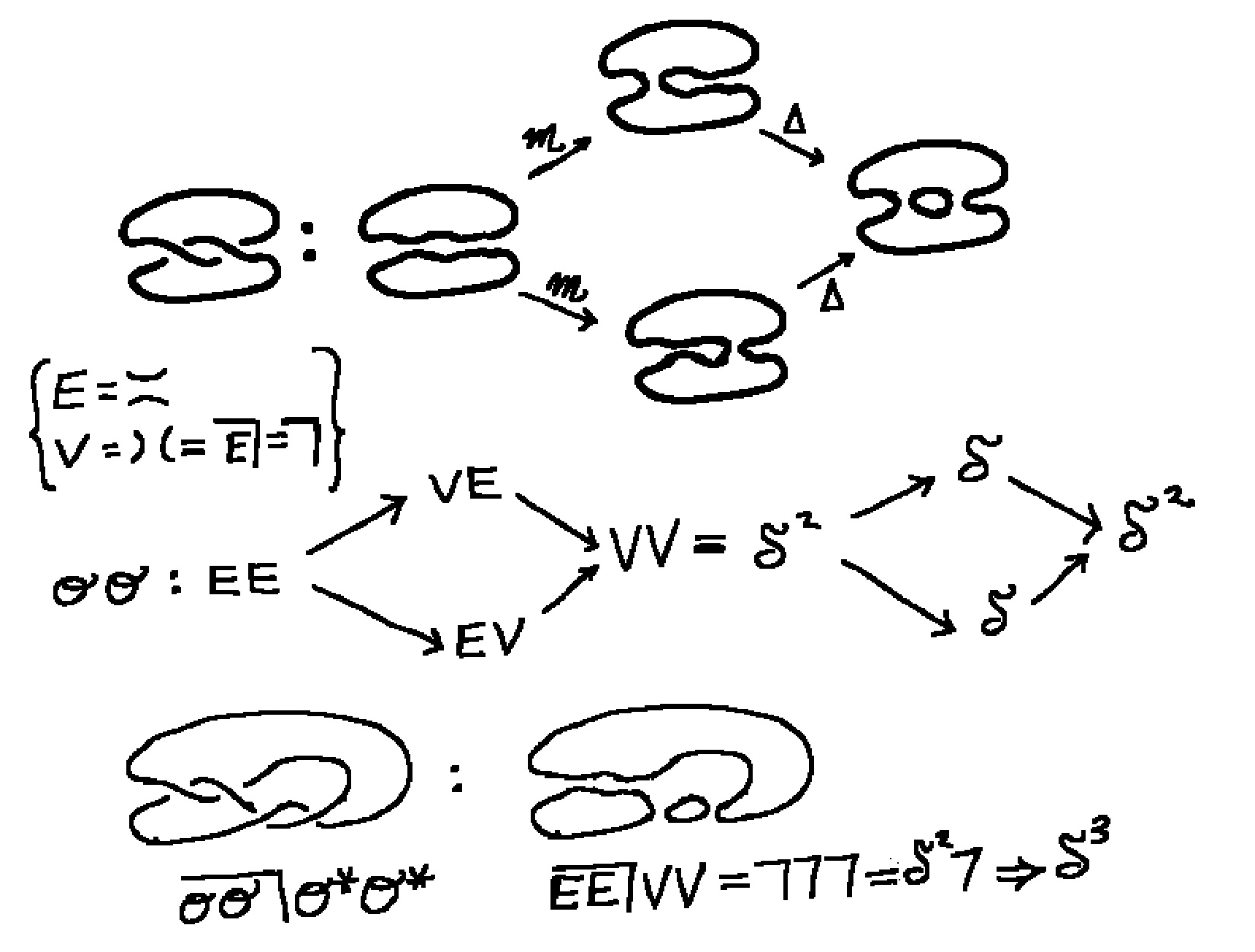}
     \end{tabular}
     \caption{\bf Khovanov Category for the Hopf Link $OO$ and one state of the Figure Eight Knot.}
     \label{hopf}
\end{center}
\end{figure}

\begin{figure}[htb]
     \begin{center}
     \begin{tabular}{c}
     \includegraphics[width=10cm]{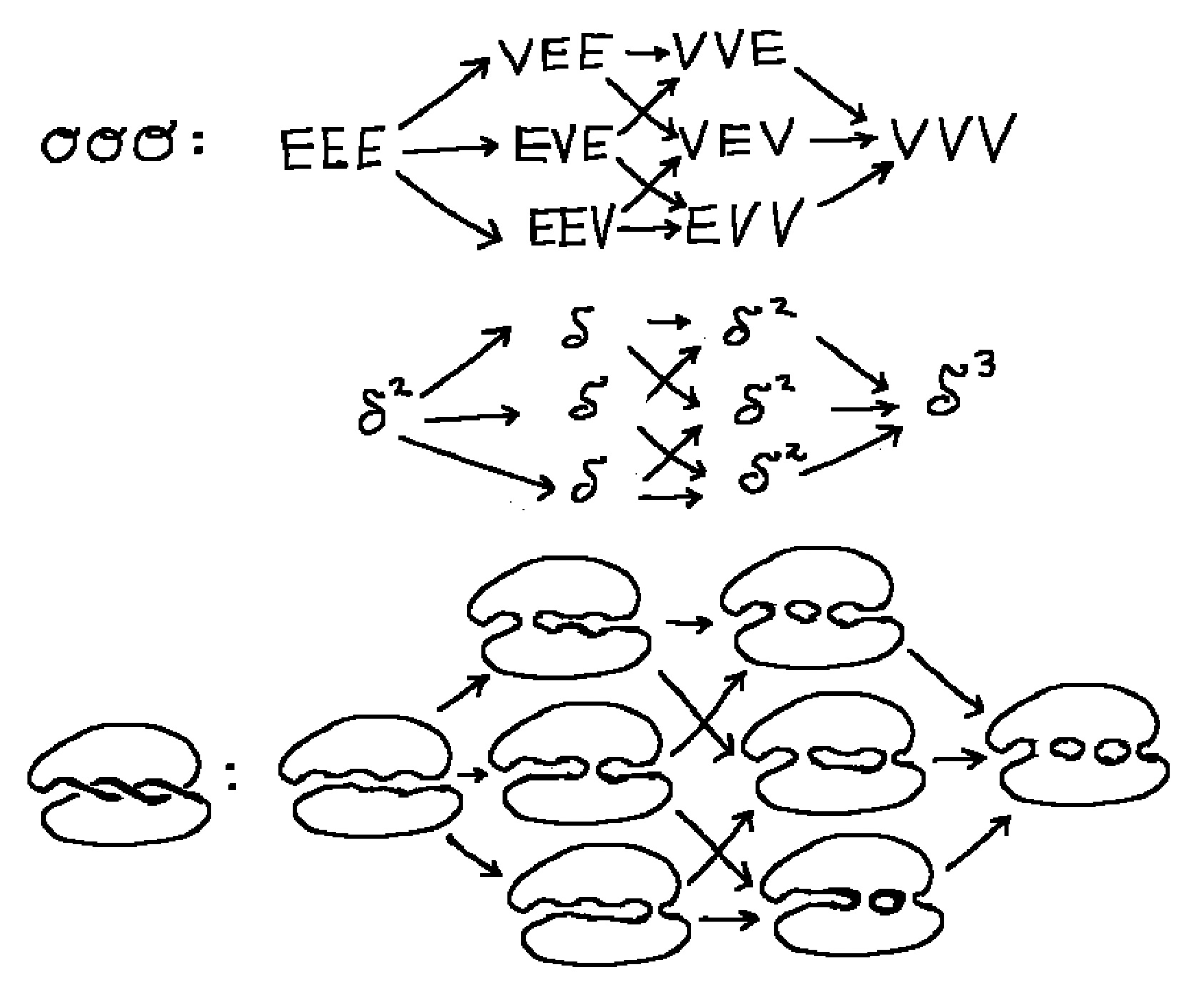}
     \end{tabular}
     \caption{\bf Khovanov Category for the Trefoil Knot $OOO$.}
     \label{trefoil}
\end{center}
\end{figure}

\noindent{\bf Remark on Khovanov Homology.}
In Figure~\ref{hopf} we illustrate the states for the Hopf link (closure of $OO$) and its bracket states in the form of its Khovanov Category \cite{Kho,KP,LKho1,LKho2} where the leftmost state is the $A$-smoothing
corresponding to $EE$ in the language of this paper. The other states proceed with arrows between them where an $A$-smoothing is replaced by a $B$-smoothing. All of this works at the level of writing the states in the form of replacing the $O$'s in $OO$ by $E$ or by $V = \cross{\,\,}$ in the notation of this figure. Then each state, seen in the crossing algebra produces $\delta$ raised to the number of loops in the state. Thus our representartion of the Khovanov Category begins to produce its structure. What is not directly derivable from our algebra is the exact nature of the smoothing sites of each state. We need to know which sites are between a loop and itself and which sites are between two loops. This information can be obtained graphically from the code $OO$ but it is not obvious how to obtain it from the algebra alone. The Figure~\ref{hopf} also illustrates one state of the figure eight knot $\cross{OO}O^{\star}O^{\star}.$  Figure~\ref{trefoil} shows the same type of analysis for the Khovanov Category of the trefoil knot $OOO.$ It would be of great interest to have a purely algebraic method for constructing all the details of the Khovanov Category for arborescent links, particularly for creating quantum algorithms for the bracket polynomial \cite{LKho1,LKho2} and for Khovanov homology.\\

Note that we can determine, by algebra, whether a given site in a state is an interaction of two separate loops or the interaction of a loop with itself. For we can change the sign at that site to $O$ (from $E$ or $V$): If the component count of the resulting expression remains the same, then the site is a self-interaction, while if the component count changes then the site is an interaction of two distinct loops. We also want to know the loop components of a given state and which sites are incident to a given loop component. This appears to require a different algebraic approach. \\

\begin{figure}[htb]
     \begin{center}
     \begin{tabular}{c}
     \includegraphics[width=10cm]{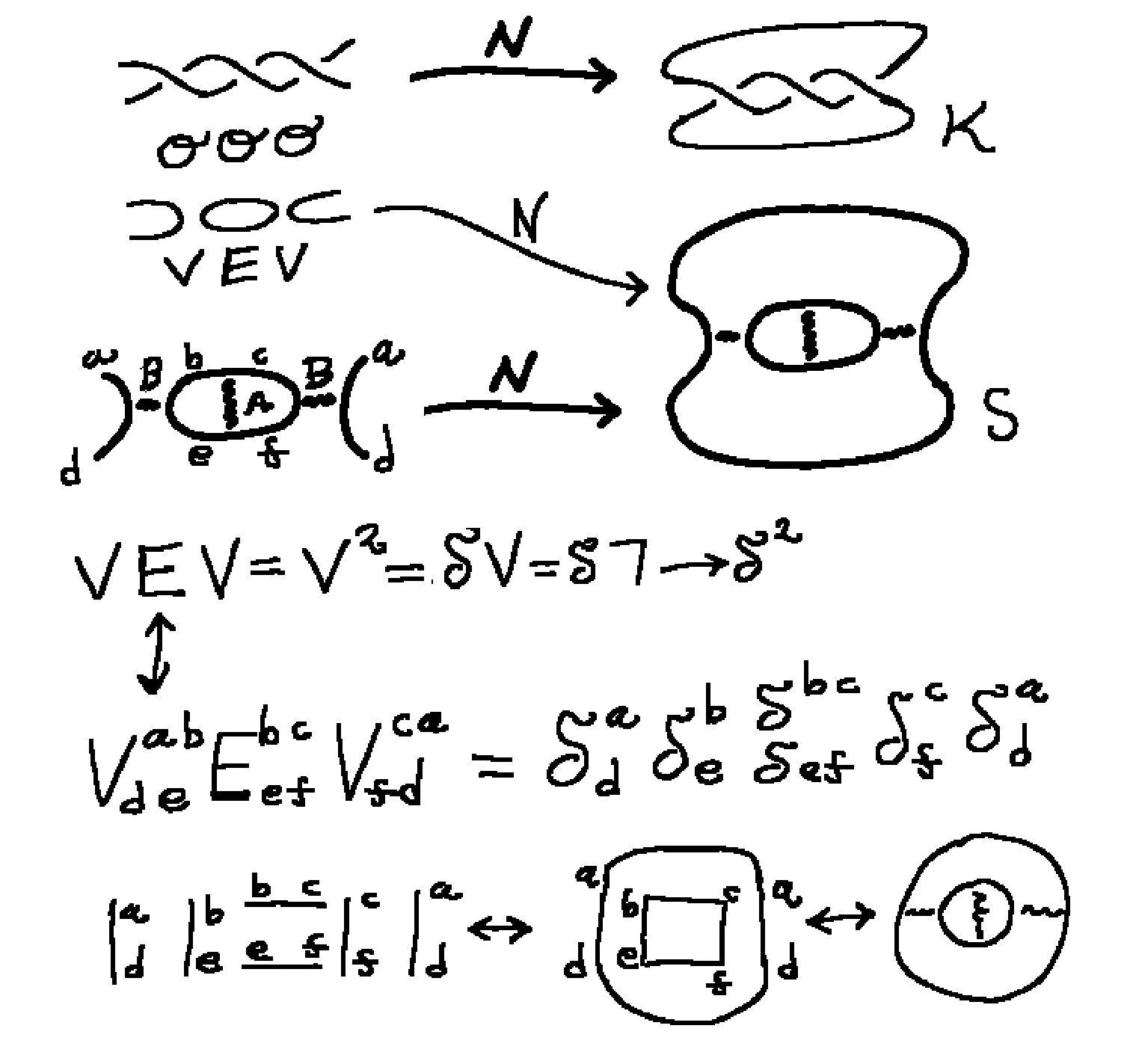}
     \end{tabular}
     \caption{\bf State Structure by Tensor Net for the Trefoil Knot $OOO$.}
     \label{tensor}
\end{center}
\end{figure}

There is an approach that accomplishes all these goals, but it involves a translation from crossing algebra to abstract tensor algebra (in the sense of \cite{P}). Regard each tangle as a box with four lines each labeled with a letter. Thus we replace 
$O$ by $O^{ab}_{cd}$ where the indices $a,b,c,d$ are here indicated to the right of the symbol to which they belong. We let $\delta^{a}_{b}, \delta^{ab}, \delta_{ab}$ (and similarly for the obvious variations) be formal Kronecker delta symbols so that 
$\delta_{ab} \delta_{bc} = \delta_{ac}.$ And we take $\delta_{aa} = \delta$ to indicate a loop. We define $$(TS)^{ab}_{cd} = T^{ai}_{cj}S^{ib}_{jc},$$ and we define
$$\cross{T}^{ab}_{cd} = T^{ca}_{db},$$
$$E^{ab}_{cd} = \delta^{ab}\delta_{cd},$$
$$V^{ab}_{cd} = \delta^{a}_{c}\delta^{b}_{c},$$
$$O^{ab}_{cd}= \delta^{a}_{d}\delta^{b}_{c}.$$
Using this tensor formalism, any product in the crossing algebra will resolve into an articulated set of loops whose sites are available from the corresponding indices. In this way, one can build the full details of the Khovanov complex from the image of the crossing algebra in the tensor algebra. Diagrams corresponding to these abstract tensors are produced by associating $\delta^{a}_{b}$ with a vertical line segement with end-labels $a$ and $b,$ and associating a similar horizontal line segment to $\delta_{ab}.$ In general, an abstract tensor is associated with a tangle drawing with labeled lines corresponding to the tensor indices. If two lines have the same index, they are joined in the tensor diagram. See Figure~\ref{tensor} for an example for the state $VEV$ of the trefoil knot $OOO.$ In that figure we take the tensor corresponding to $VEV$ to have matching indices from left top to right top and left bottom to right bottom so that the diagram corresponds to the loop closure for the state of the knot. Note that from the tensor decomposition, the product of deltas factors into closed loops and one can see from associated indices whether the sites in the state are internal to the loops or between one loop and another. Thus the full Khovanov category can be read fron the states of a knot or link that is expressed in the crossing algebra. It is possible that for arborescent links and tangles there may be a simpler route to the Khovanov complex other than through the abstract tensor algebra. That is a problem to be investigated beyond the present paper and to be compared with physical approaches such as \cite{Atiyah,Witten, Baxter,KP}\\

\begin{figure}[htb]
     \begin{center}
     \begin{tabular}{c}
     \includegraphics[width=11cm]{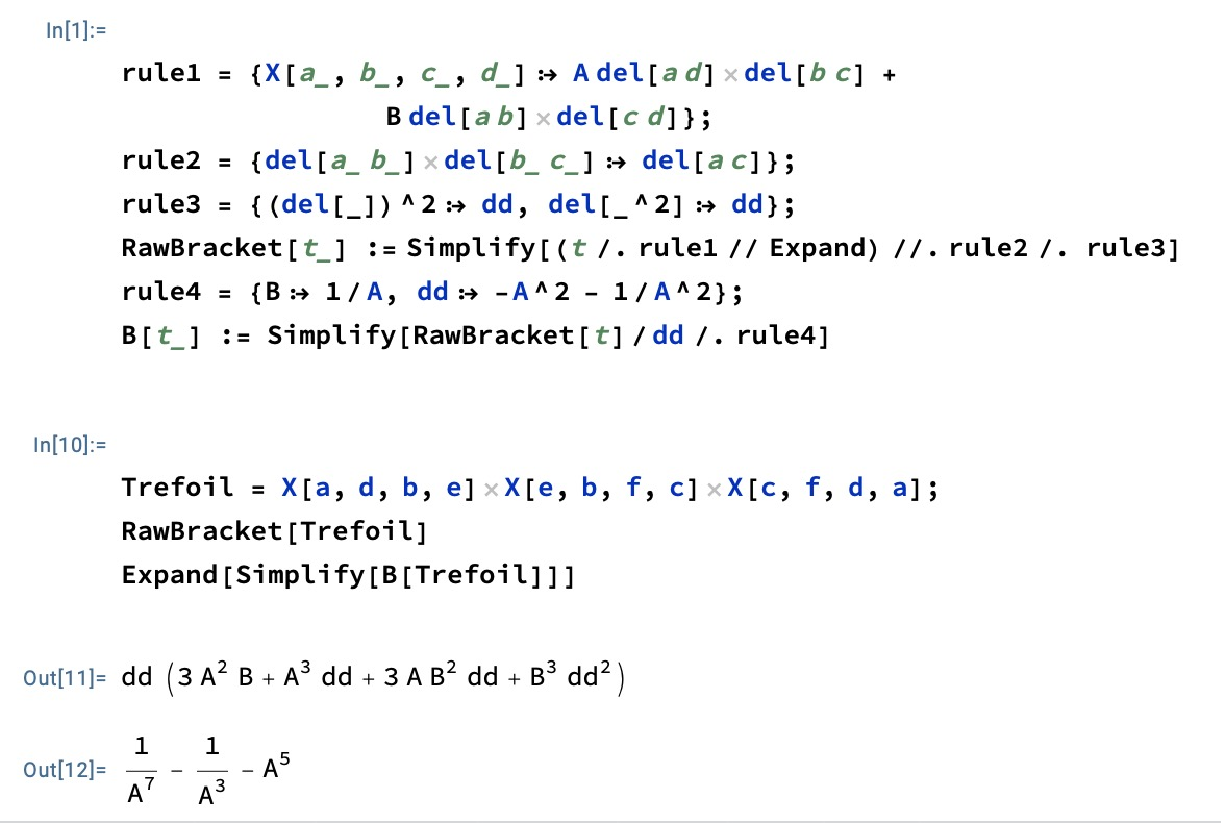}
     \end{tabular}
     \caption{\bf Abstract Tensor Expansion for Bracket Calculation of the Trefoil Knot $OOO$}
     \label{brackettrefoil}
\end{center}
\end{figure}

Returning to the bracket polynomial calculation, one can define a tensor expansion directly in terms of a crossing by the formula
$$\Across^{cb}_{da} = X[a,b,c,d] = A\delta^{cb} \delta_{da} + A^{-1}\delta^{c}_{d}\delta^{b}_{a}$$
corresponding the the bracket expansion
$$\Across = A \Asmooth + A^{-1} \Bsmooth.$$
This is instantiated in a general-purpose bracket calculation program as in Figure~\ref{brackettrefoil}. Note that in the notation $\Across^{cb}_{da}$ we can take the tensor labels on the crossing edges as proceeding counterclockwise from the lower right in the order $abcd$ and this is the mnemonic for the symbol $X[a,b,c,d]$ that is used for coding. The program illustrated in this figure is a highly efficient method for producing the bracket polynomial and is modified to handle Khovanov homology as well. (Note that $dd$ represents the loop value in this program.) All this occurs at the tensor network level for any knot or link diagrams.\\
 
\section{Medial Graphs}
This section reviews the reformulation of knot theory in terms of graphs - via the medial construction. We show how the arborescent forms can be viewed as plane graphs and 
how the component count is the nullity of a mod two Laplacian matrix for the graph. Our algebraic approach to the counting of components leads to problems of generalisation in the graph category.\\

\begin{figure}[htb]
     \begin{center}
     \begin{tabular}{c}
     \includegraphics[width=10cm]{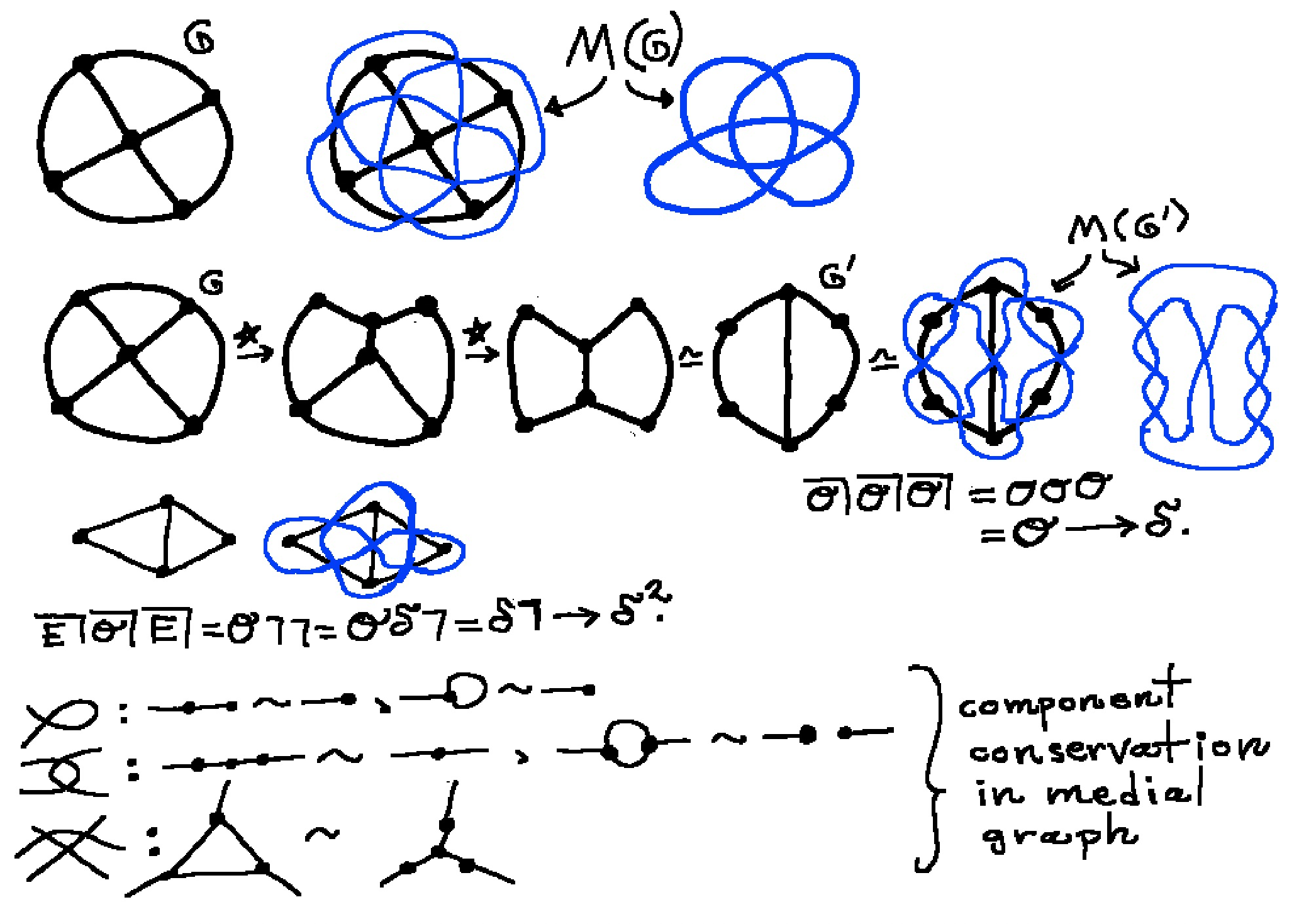}
     \end{tabular}
     \caption{\bf Medial graphs, component counts and graphical Moves preserving component count.}
     \label{medial}
\end{center}
\end{figure}

We have discussed the structure of knot and link diagrams that are coded in algebraic expressions that generalise continued fractions. There is a larger category of expressions
that encode knots and links and this is the category of all finite connected  plane graphs. Given a plane graph $G$ one can construct its {\it medial graph}  $M(G),$ a 4-valent plane graph. The medial is constructed by placing a flat crossing (local 4-valent node) at an interior point of each edge of $G,$ and then connecting these local edges along the boundaries of the regions of $G.$ The process is shown in Figure~\ref{medial}. \\

\begin{figure}[htb]
     \begin{center}
     \begin{tabular}{c}
     \includegraphics[width=10cm]{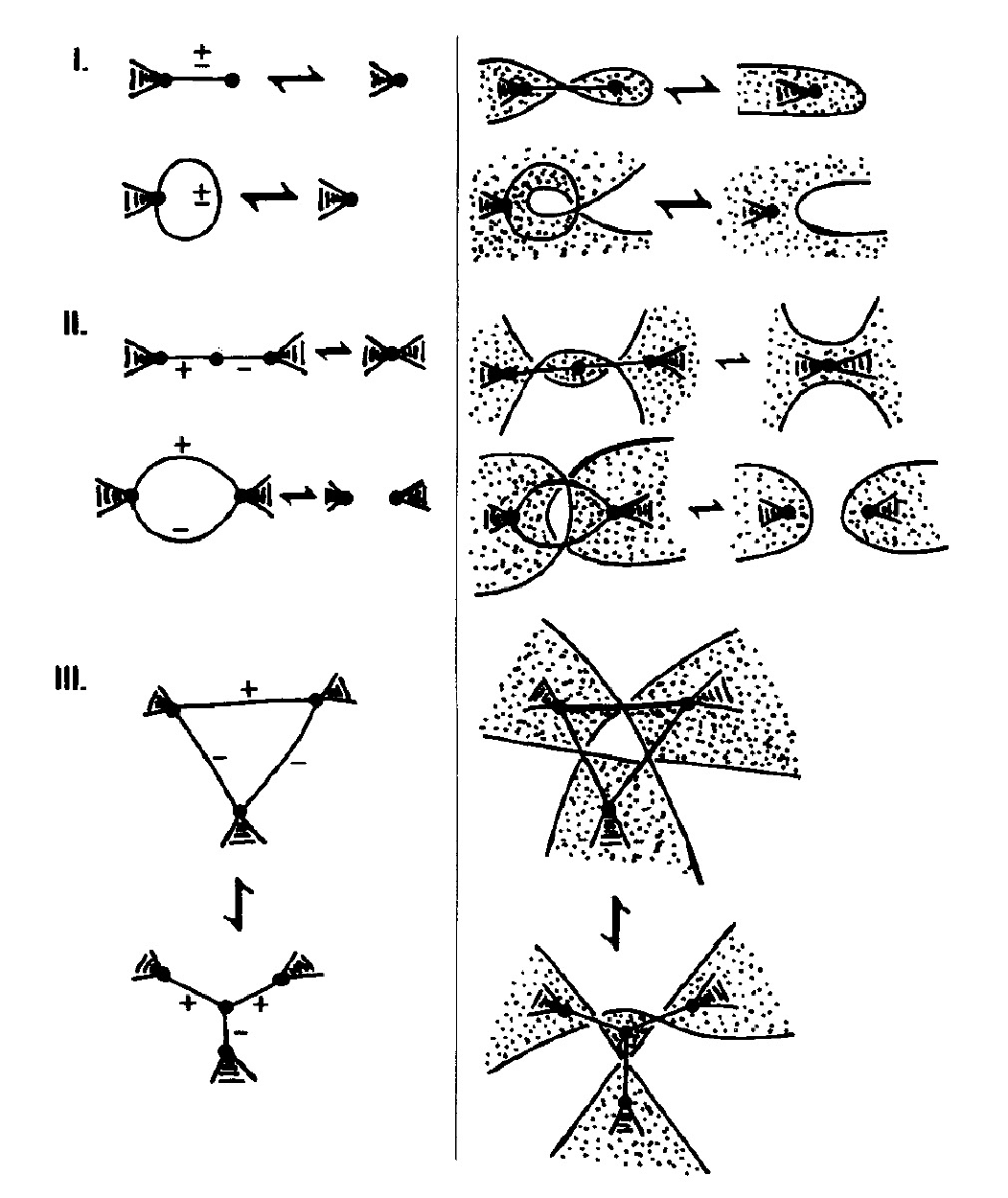}
     \end{tabular}
     \caption{\bf Graphical Reidemeister Moves}
     \label{grm}
\end{center}
\end{figure}

\begin{figure}[htb]
     \begin{center}
     \begin{tabular}{c}
     \includegraphics[width=10cm]{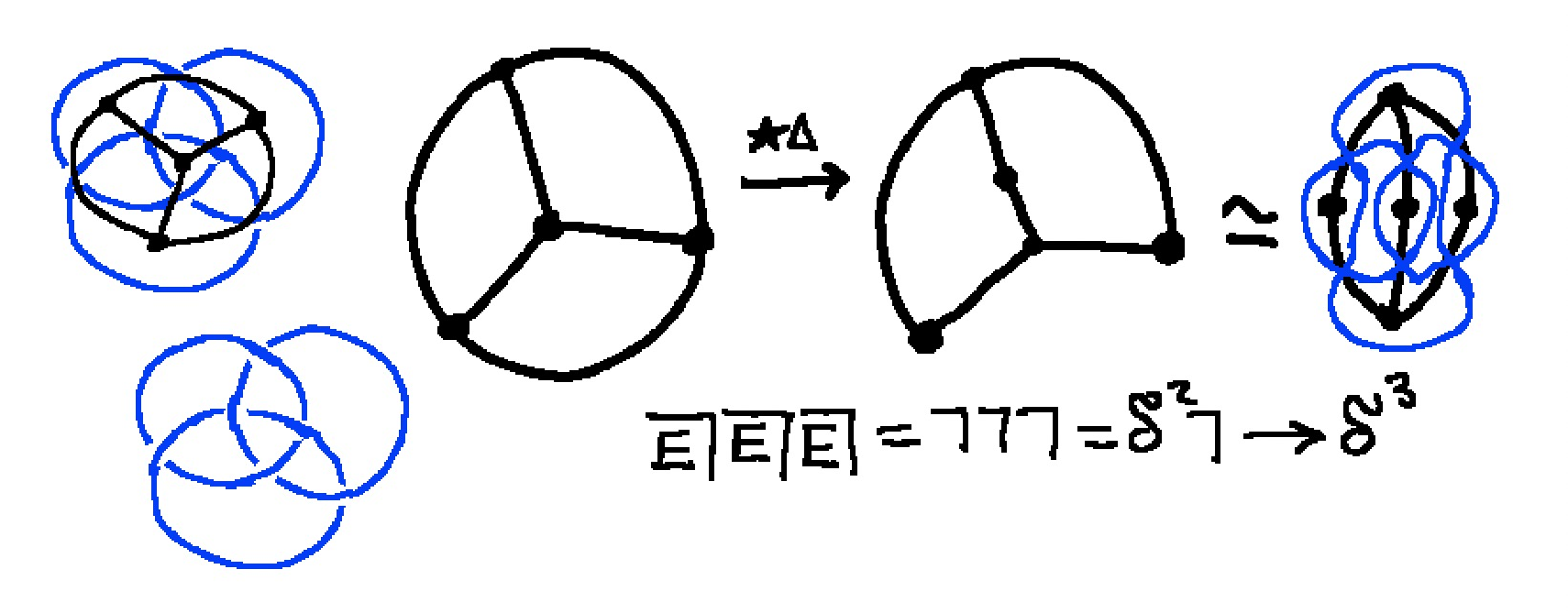}
     \end{tabular}
     \caption{\bf Graph for Borommean Rings.}
     \label{boro}
\end{center}
\end{figure}

The medial graph can be seen as a flat (no given over or under crossings) link diagram. See See Figure~\ref{grm} and Figure~\ref{medial}. Conversely, given a classical link diagram $K$, one can shade its regions with colors black and white so that adjacent regions have distinct colors. The shaded regions can each be assigned a graphical node and edges are constructed between two nodes if there is a crossing in the link diagram that is common to the two regions. The resulting graph, $G(K),$ is called the checkerboard graph of $K.$ The medial of $G(K)$ is equal to $Flat(K)$ where $Flat(K)$ is the 4-valent graph obtained from $K$ by ignoring the over and under crossing data in $K.$ One can also choose under and over crossings for the medial graph according to signs on the edges of the plane graph as shown in Figure~\ref{grm}.\\

In Figure~\ref{grm} we illustrate the translation of the Reidemeister moves for link diagrams to moves on their corresponding plane graphs. Here we show the general case where an arbitrary link diagram is translated to a {\it signed} plane graph (each edge has a plus or minus sign attached to it as a label). The signs correspond to the way the crossing interacts with the edge in the checkerboard graph. Smoothing the crossing in relation to the edge produces an $A$-smoothing with a plus edge and a $B$-smoothing with a minus edge. The key  point is that the Reidemeister moves translate into graphical moves of the following types:
\begin{enumerate}
\item Add or remove a pendant loop.
\item Add or remove a pendant edge.
\item Contract two edges in series when the edges have opposite signs.
\item  Delete two edges that are in parallel when the edges have opposite signs.
\item Change a triangle for a star or change a star for a triangle (with certain signed labels).
\end{enumerate}
See Figure~\ref{grm} for the details. This association amounts to a complete translation between a theory of moves on signed planar graphs and classical knot and link theory via diagrams.
However, given a plane graph it is not immediately obvious how many components there will be in its corresponding knot or link. This problem for plane graphs is exactly analogous to our 
problem of counting for arborescent links presented in algebraic forms. See Figure~\ref{medial} for examples of medial graphs and their component counts.\\

Just as we have analyzed the algebraic expressions of knots and links to determine the number of link components in their realizations, one can analyze the medial graphs of plane graphs $G$ to find out their number of components in the sense of link components. There are at least two distinct methods for determining this component count. One can construct the medial and then trace on it those cycles that correspond to the link components. This is, just as for link diagrams corresponding to arborescent links, tedious except for small examples. For the plane graphs a second method is to associate a matrix $Q_{2}[G] =(q_{ij}),$ the mod-2 graph Laplacian for the graph $G.$ Then the Nullity of the mod-2 Laplacian of $G$ is equal to the number of link components of the medial graph \cite{SW}. For the reader interested in examining this algebraic method, the definition of the Laplacian is as follows
\begin{enumerate}
\item $q_{ii}$ = the degree of the $i$-th node of $G$.
\item When $i \ne j$ then $q_{ij}$ = the number of edges between the $I$-th and $j$-th nodes of G.
\item Let $\mu(G)$ denote the dimension of the null space of $(q_{ij})$ taken modulo two as a linear map of vector spaces over the field with two elements.
\item Then $\mu(G)$ is equal to the number of link components of the medial graph $M(G).$
We will sketch a proof of this statement and discuss how this general component count is related to the counting properties of the crossing algebra.
\end{enumerate}

In working with the question of component counts it is not necessary to keep track of the crossing signs. Then we have a simplified set of moves exactly as in Figure~\ref{grm} except that one can ignore the signs. These are moves on unsigned plane graphs and are shown in Figure~\ref{medial}. They correspond to flat Reidemeister moves on flat knot and link diagrams. We know \cite{KP} that any flat link diagram can be transformed by flat Reidemeister moves to a disjoint collection of circles. The number of circles is equal to the number of components of the original link. Indeed each Reidemeister move preserves the component number. Thus for unlabeled plane graphs, each graphical move preserves the number of components in its medial. Every plane graph can be transformed by graphical moves to a disjoint union of isolated nodes. The number of nodes is equal to the number of components in the medial graph of the given graph.
With this understanding, we can sketch the proof that the nullity of the mod-2 Laplacian matrix gives the component count: {\it One verifies that the nullity is not changed by any of the graphical moves.} This completes the proof.\\

In Figure~\ref{medial} we illustrate the conversion of a plane graph to its medial graph, and we illustrate how graphical moves transform the graph to a graph that is the medial of an arborescent knot. In Figure~\ref{boro} we illustrate how  a plane graph whose medial has three components (and corresponds, with appropriate crossing choice, to the Borommean rings) is transformed by graphical moves to a graph with corresponding arborescent code $\cross{E}\cross{E}\cross{E}.$ It is an immediate consequence that the medial link has three components. A calculation of the nullity will confirm this result. \\

We see from these examples that the plane graph category can be seen as an extension of the algebraic coding of knots and links. It would be of interest to see a generalisation of the crossing algebra that can be directly applied to all plane graphs.\\

\section{Circuit Logic}
An independent motivation for the algebraic constructions in this paper comes from a problem in the design of switching circuits.
Recall that Claude Shannon \cite{Shannon} discovered a direct correspondence between Boolean algebra and the properties of networks of binary switches.
A switch has two states (denoted 0 and 1 or open and closed) and the most elementary switch acts on a single line to leave it open to passing a signal or breaking the line so that it cannot pass a signal.
Switches $p$ and $q$  wired in series correspond to a logical and, while switches connected in parallel correspond to logical or. Letting $a \vee b$ denote `a or b" and $a \wedge b$ denote `a and b" and 
$\sim a = \bar{a}$ denote ``not a", we can design circuits whose signal passing behaviour is the exact structure of any given logical expression.  \\

\begin{figure}[htb]
     \begin{center}
     \begin{tabular}{c}
     \includegraphics[width=10cm]{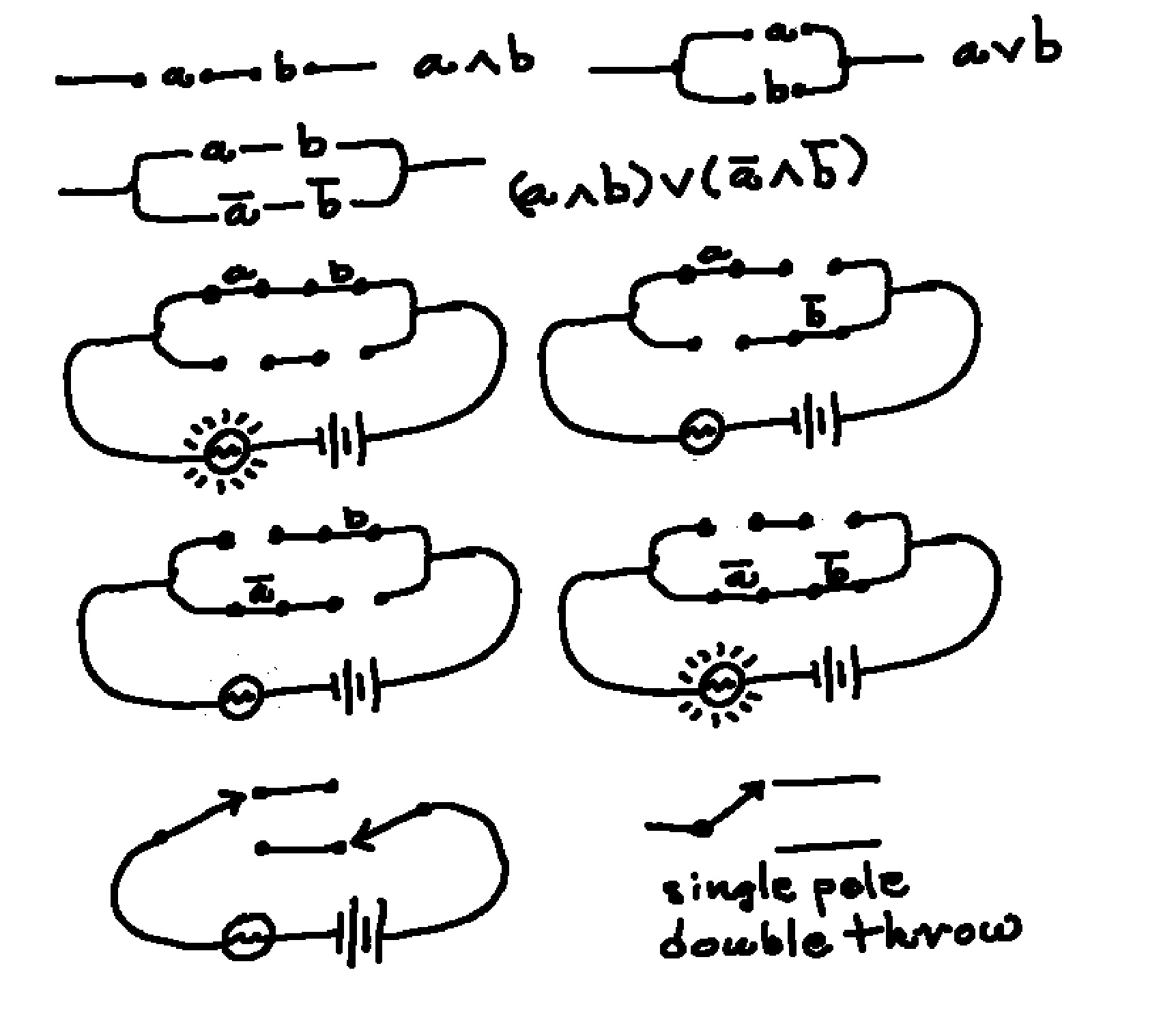}
     \end{tabular}
     \caption{\bf Logic and Switching Circuits.}
     \label{logoswitch}
\end{center}
\end{figure}

In Figure~\ref{logoswitch}  illustrates this translation between
symbolic logic and switching circuits. Specifically, the figure shows how the logical expression $ (a \wedge b) \vee (\bar{a} \wedge \bar{b})$ translates into a circuit with two switches, so that each switch independently controls a single light bulb. If the bulb is lit, then either switch can extinguish it. If the bulb is unlit, then either switch can bring the bulb to life. The logical expression indicates the two switch states that successfully light the bulb.
In the first case the switches are in the states $a$ and $b$ and in the second case the switches are both in their opposite states. \\

\noindent {\bf Remark.} Switching circuits generalise to electrical circuits where the switches are replaced by conductances. Then the variables $a,b,\cdots$ can take values in the real numbers with a formal value of $\infty$ added for the closed switch (infinite conductance) , retaining $0$ for the open switch. Let $C(N)$ denote the conductance of a network that has one input line and one output line.We have the  following two formulas for conductance using our logical notation for series and parallel connections: $$C(a \wedge b) = a + b,$$ $$C(a\vee b) = 1/((1/a) + 1/b)).$$  Note also that $C(a\vee b) = 1/((1/a) + (1/b)) = \cross{\cross{a}\cross{b}}$ (using our notation in this paper) indicates that the De Morgan Law holds in this extended boolean logic. With this we see that we should define $$\bar{a} = 1/a = \cross{a}$$ where it is understood that
$1/0 = \infty$ and $1/\infty=0.$ In this way $0$ and $\infty$ form the boolean algebra inside an $\it electrical$ extension of boolean logic with real values using these equations. \\

Note that in this extension $1$ and $-1$ are invariant under negation. This point of view connects with our work on conductivity and link invariants \cite{GoldKauff}. In that paper we associate a signed graph $G(T)$ to any knot or tangle $T$ by the checkerboard/medial construction described in the previous section of the present paper. The conductivity of $G(T)$ between two chosen nodes in the graph is a topological invariant of the knot or tangle (restricted to not move across the selected nodes). For example, the checkerboard graph of the Borommean rings (see Figure~\ref{boro}) has all crossings of the same type since the rings are alternating and so all conductances will be non-zero, indicating the linkedness of the rings.\\

\begin{figure}[htb]
     \begin{center}
     \begin{tabular}{c}
     \includegraphics[width=10cm]{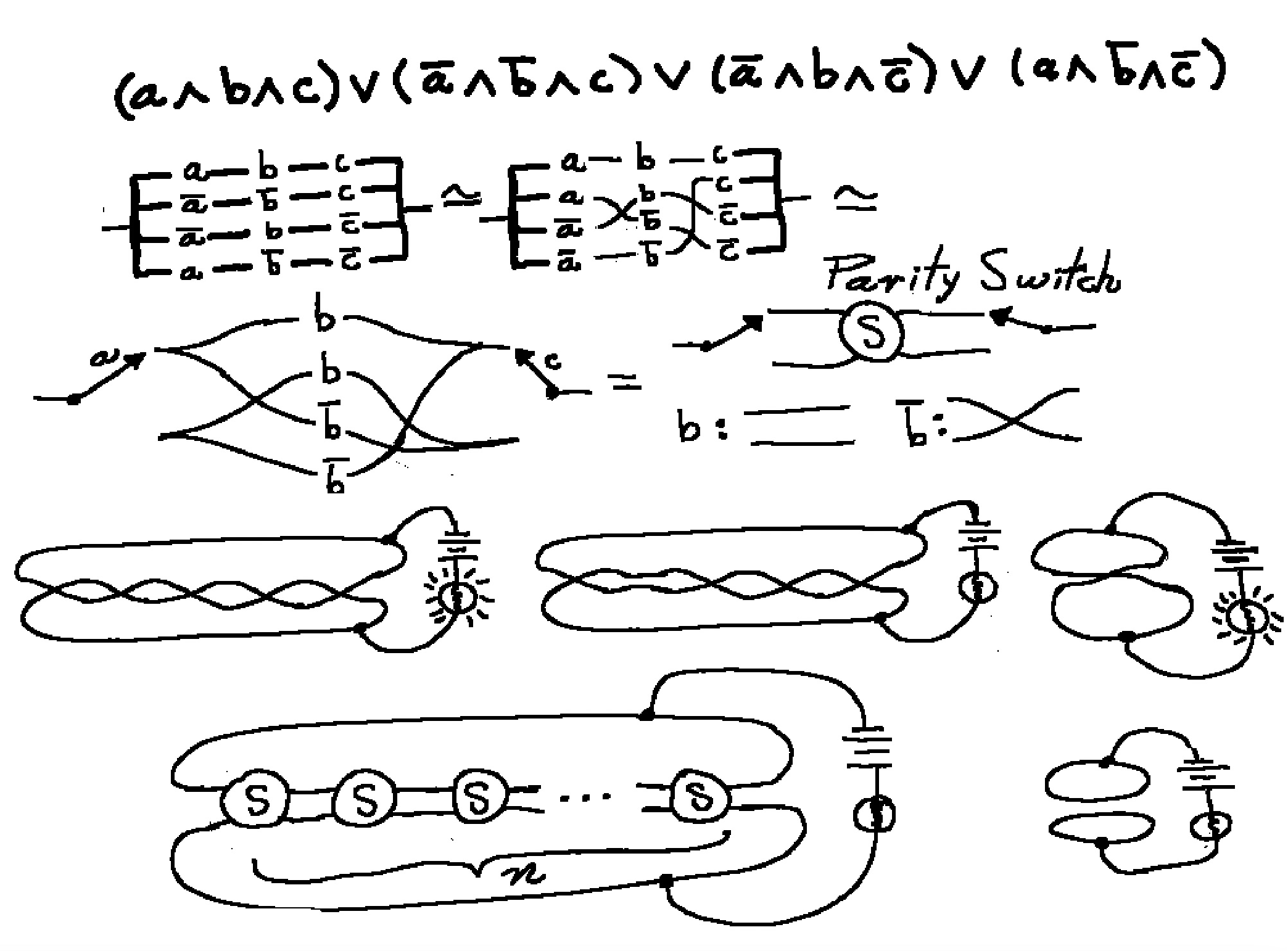}
     \end{tabular}
     \caption{\bf Lightbulb Control and the Cross Switch}
     \label{control}
\end{center}
\end{figure}

\begin{figure}[htb]
     \begin{center}
     \begin{tabular}{c}
     \includegraphics[width=10cm]{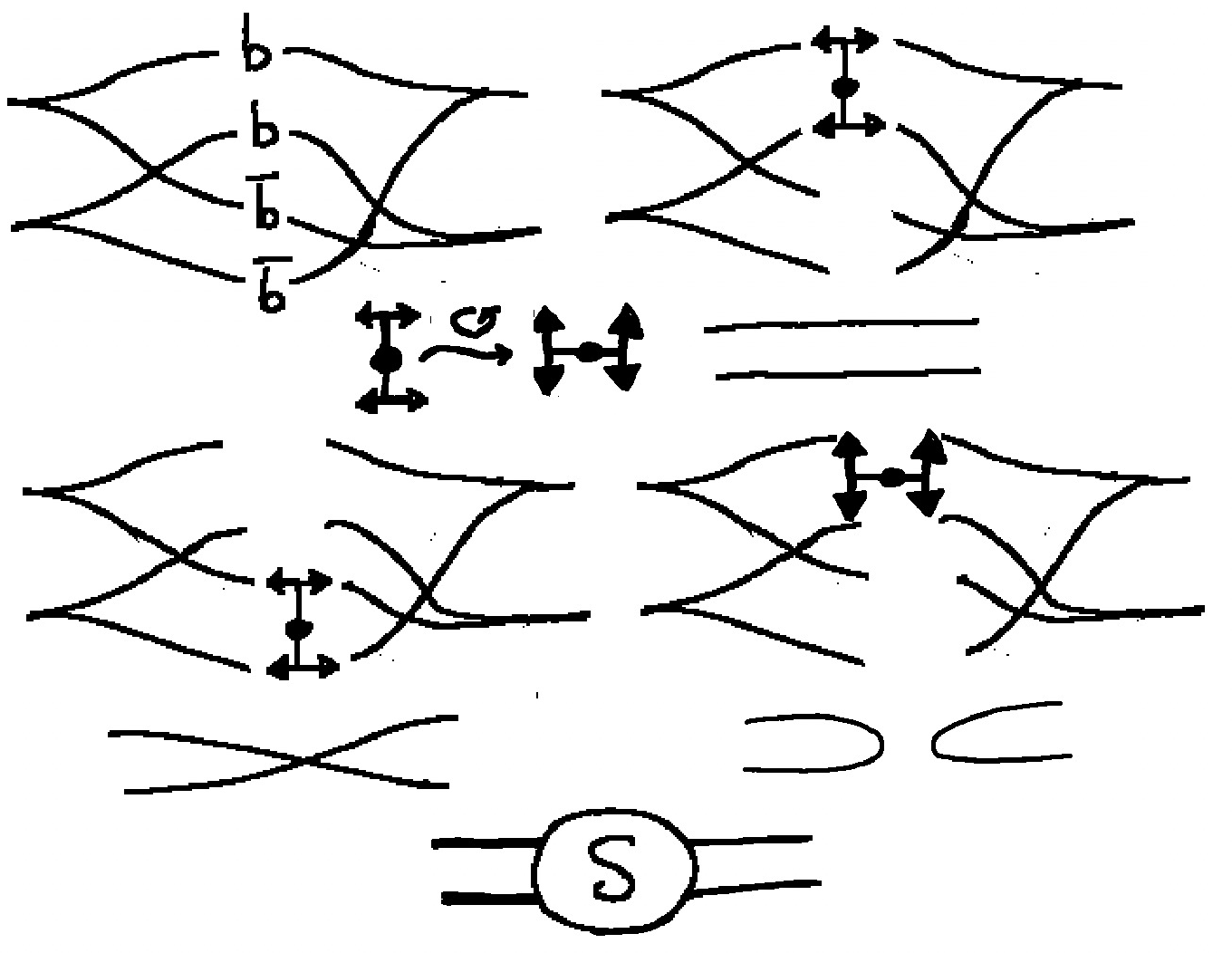}
     \end{tabular}
     \caption{\bf Three State Cross Switch - Slide and Rotate}
     \label{crosswitch}
\end{center}
\end{figure}

In Figure~\ref{control} we illustrate the corresponding translation for the case of three
binary switches controlling a single bulb. Now the symbolic logical condition is a disjunction of four conjunctions as there are four conditions that will light the bulb. The first part of the figure illustrates the translation to a switching device that will pass current from right to left exactly under these four switching conditions. We illustrate how a person with a sharp eye for the structure of this device can see it as consisting in two single pole, double throw switches at the ends, and a more complex switch ``b" in the middle. A careful look at the middle switch reveals that it can be seen as as ``crossing switch" $S$ where the switch $S$ has two input lines and two output lines, and its two states consist in these lines being either parallel or crossed over. Note that in our notation, two switching lines that cross through one another do not actually touch in a circuit realization. Thus the crossing is 
virtual. It is a bit of luck to find the crossing switch $S,$ as it allows one to design a circuit that can control one light bulb with any number $n$ of the crossing switches. The rest of the figure shows how this works. By connecting 
$n$ crossing switches in series, as shown, and connecting to form a closed loop arrangement as shown, we find that the circuit draws will light the bulb when and odd number of switches are in the crossed position. Figure~\ref{crosswitch} illustrates the cross switch in isolation. The idea for using a cross switch to solve the $n$-switch light bulb problem is due to the cyberneticist Ricardo Uribe \cite{Uribe}. Uribe liked to refer to this as a 
paradoxical solution and also he liked to point out that the basic structure of the solution when the light is on, is a Mobius band. Just so, if we put the terminals of the bulb/battery part of the circuit on nearby edges of a Mobius (using the edge of the band as wire), then the current can run between them because the band has only one edge. In an ordinary band, there are two edges and the current cannot get across. Thus does topology appear in the logic of circuit design.\\

\begin{figure}[htb]
     \begin{center}
     \begin{tabular}{c}
     \includegraphics[width=10cm]{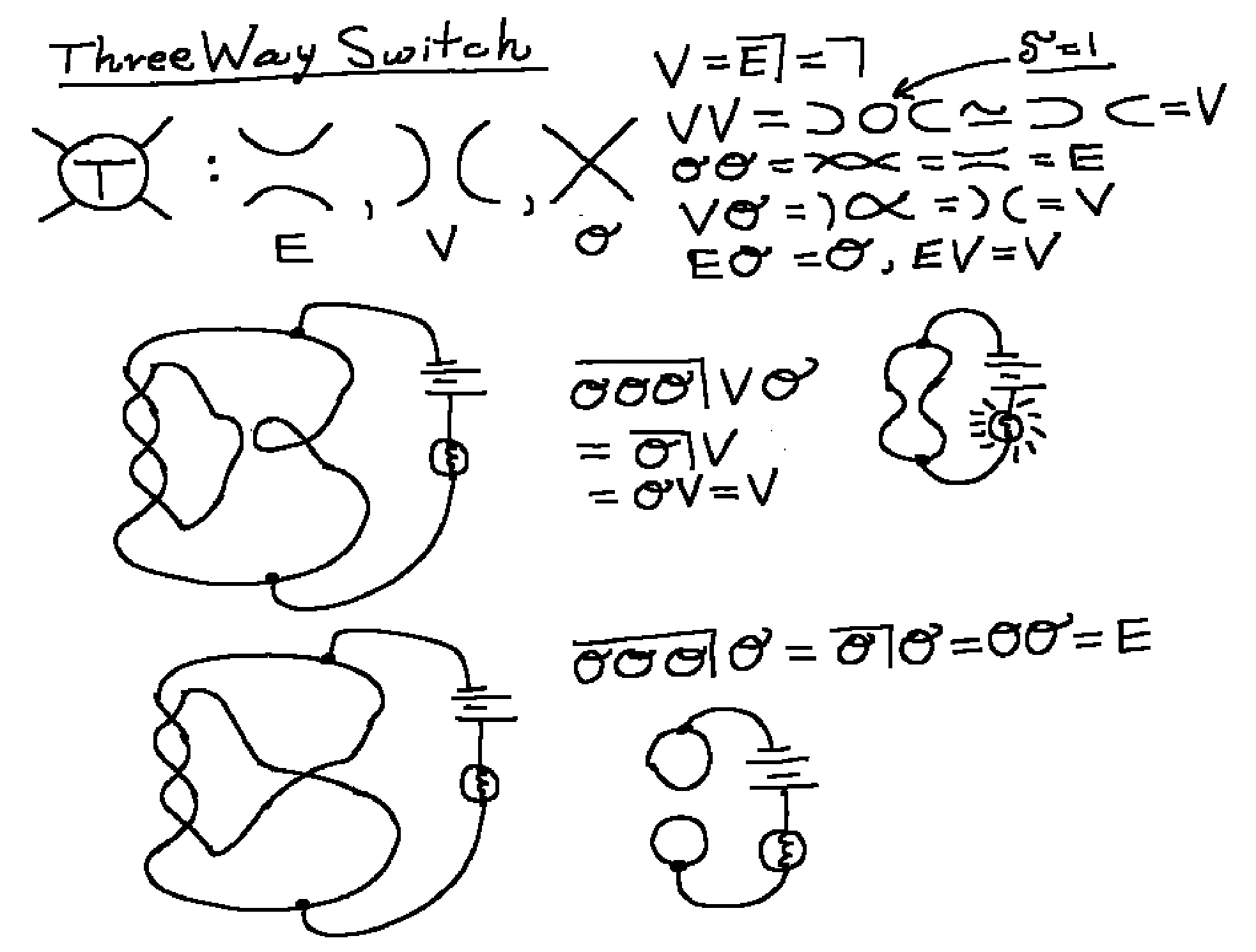}
     \end{tabular}
     \caption{\bf From Cross Switch to Three State Cross Switch}
     \label{three}
\end{center}
\end{figure}

In Figure~\ref{crosswitch} we show the details of the cross switch using a sliding contact to accomplish the states of the switch, with  $b$ operative when the sliding contact is in the up position and $\bar{b}$ 
operative when the sliding contact is in the down position. In this figure, we add the possibility of rotating the sliding contact by ninety degrees. This gives a {\it third} state for the switch as illustrated in the figure.
Thus we now have a three-state switch where the two lines are effectively either horizontally parallel (E),  vertically parallel (V), or crossed over one another (O). These three states correspond exactly to our crossing algebra states where $V = \cross{E} = \cross{\,\,}$ as discussed earlier in the paper. Figure~\ref{three} illustrates the correspondence. In Figure~\ref{three} we show how we can take a flat diagram for a rational or arborescent link and convert it to a switching circuit so that the light bulb will be on when there is (in the rational case) one component and the light will be off when there are two components. In the algebra for such switching circuits we use
$V V = \cross{\,\,} \cross{\,\,} = \cross{\,\,} = V$, taking $\delta = 1$ in the previous formulation. This means that local closed components produced by the $VV$ interaction will not be counted, but they do not contribute to the 
question of circuit connectivity as we have wired it up in these figures. The crossing algebra is correspondingly simpler and more like a multiple valued logic for these three way switching circuits. It is ironic that the logic of these circuits depends on the ``paradoxical" element $O$ with $\cross{O} = O,$ perhaps confirming Uribe's idea of the paradoxical nature of this design.\\

\begin{figure}[htb]
     \begin{center}
     \begin{tabular}{c}
     \includegraphics[width=10cm]{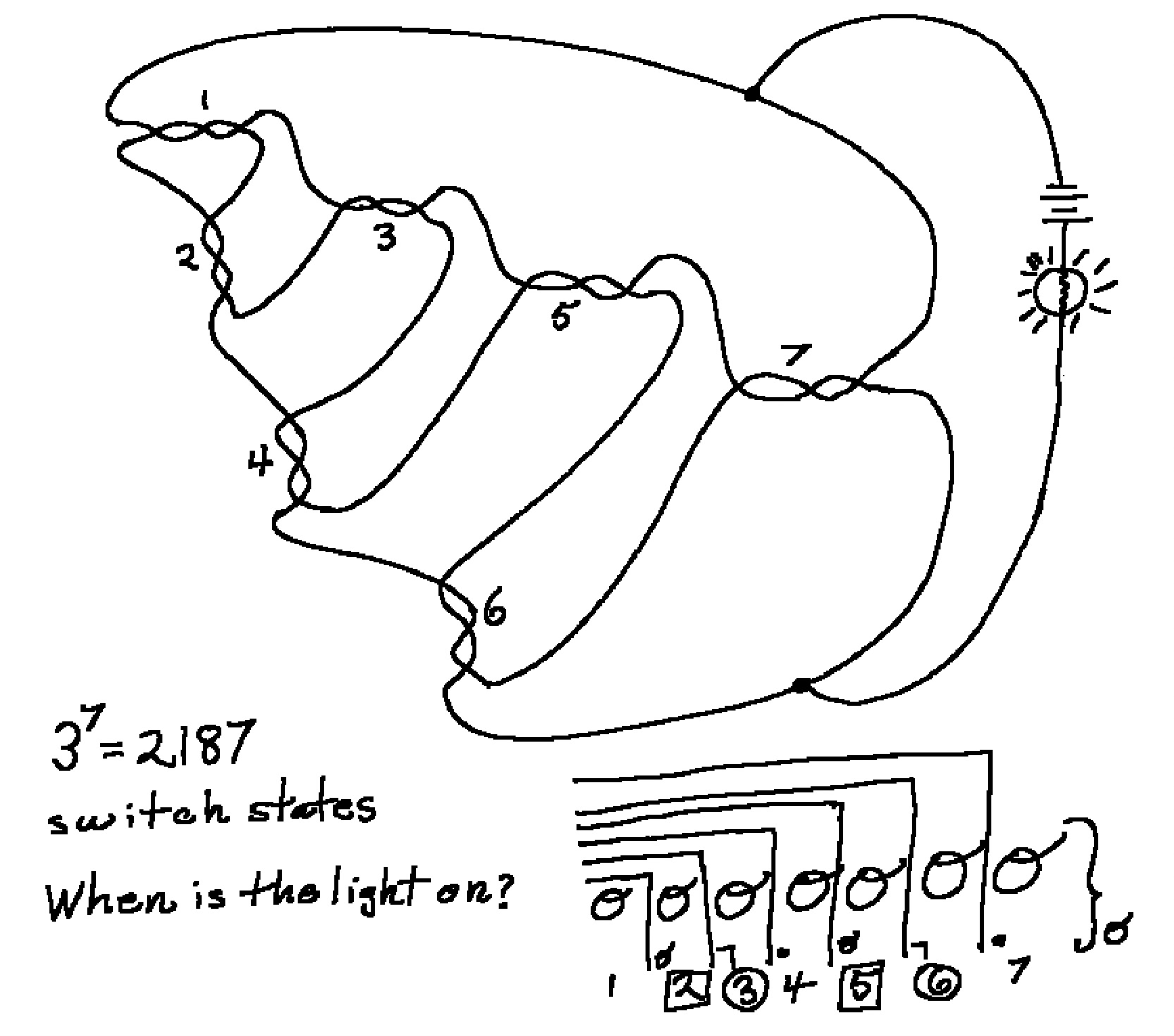}
     \end{tabular}
     \caption{\bf Circuit Design and Crossing Algebra}
     \label{circuit}
\end{center}
\end{figure}

In Figure~\ref{circuit} we illustrate a larger circuit and its crossing algebra expression. The reader can apply all the techniques of this paper and find out (without tracing the circuit) that the light will be on for the circuit state that is indicated in the figure. And the reader can use the crossing algebra to find those switches where one application of the switch between $E$ and $O$ will or will not change the state of the bulb.
It is of interest to note that in this way we revisit the concept of Claude Shannon that there should be an algebraic understanding of switching circuits, and in this case the algebra is non-boolean. Imaginary logical values such as our paradoxical $O$ are crucial for understanding the circuit behaviour.\\

\begin{figure}[htb]
     \begin{center}
     \begin{tabular}{c}
     \includegraphics[width=10cm]{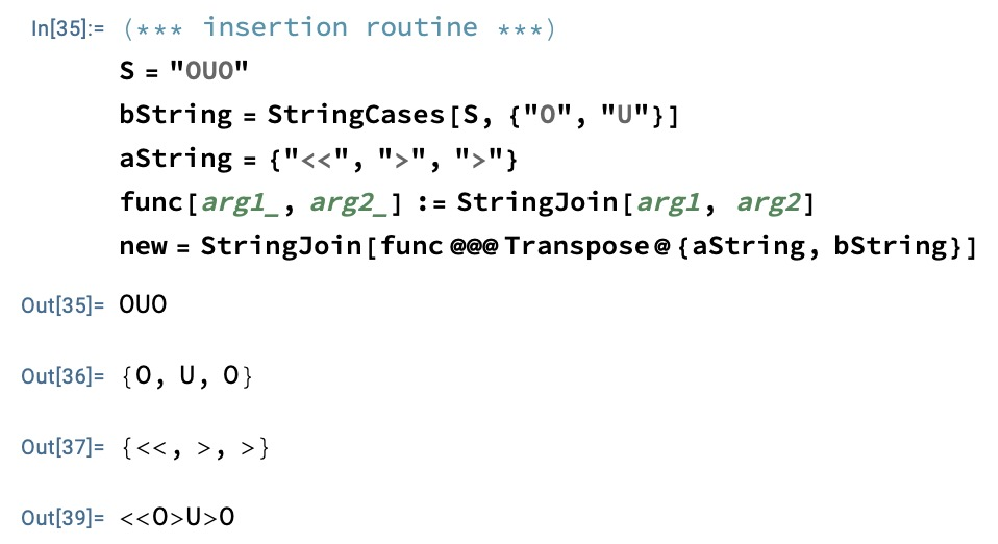}
     \end{tabular}
     \caption{\bf Insertion Code}
     \label{insert}
\end{center}
\end{figure}

\section{Logic and Foundations}
The crossing algebra is closely related to the calculus of indications of G. Spencer-Brown as it is explained in his work Laws of Form \cite{GSB}.
The calculus of indications of George Spencer-Brown (GSB) corresponds to crossing algebra with only the cross $\cross{\,\,}$ and evaluations of expressions involving the cross with $d=1.$
Thus the calculus of indications is given by the rules:\\

\noindent {\bf Calculus of Indications.}
 $$\cross{\cross{\,\,}} = \,\,\,$$
$$\cross{\,\,}\cross{\,\,} = \cross{\,\,}$$\\

Here we use the empty word for $\cross{\cross{\,\,}}.$ Every expression simplifies uniquely to a single mark or to the empty word.
Thus $$\cross{\cross{\,\,}\cross{\,\,} \cross{\cross{\cross{\,\,}} }}\cross{\,\,}= \cross{\cross{\,\,} \cross{\cross{\cross{\,\,}} }} \cross{\,\,}
= \cross{\cross{\,\,} \cross{\,\,}}\cross{\,\,} =\cross{\cross{\,\,}}\cross{\,\,} = \cross{\,\,}.$$

For the Calculus of Indications a natural interpretation is in terms of Boolean algebra, as explained in \cite{GSB} and \cite{MLogic,KnotLogic}. The first generalisation that we have used is to
add an element $O$ so that $$\cross{O} = O,$$ $$\cross{\,\,} O = \cross{\,\,},$$ $$OO = \,\, .$$ Call this the {\it Contracted Crossing Algebra} ($CCA$). The contracted crossing algebra is just sufficient for
finding the number of components in a rational knot of link in continued fraction form, as we have discussed in the previous sections of the paper. We call attention to it here because it bears a striking resemblance to a three valued extension of the calculus of indications that is essentially mapped to the three valued logic of Lukasiewicz \cite{L}. In this calculus for self-reference (CSR) \cite{CSR,FD,SRRF} one has an element $@$ in the calculus so that 
 $$\cross{@} = @,$$ $$\cross{\,\,} @ = \cross{\,\,},$$ $$@@ = @ .$$ 
Thus in the $CSR$ we have the same rules for $@$ as for $O$ in the contracted crossing algebra except that $@@=@$ while $OO= \,\,.$\\

The $CSR$ and its correlative three-valued Lukasiewicz logic is a natural extension of a two-valued logical calculus in the face of the paradoxical element $@$ with $\cross{@} = @.$
Such logics were originally designed in the light of adding to the usual values of True and False some intermediate values such as Possibly True or Possibly False. In the case of the simple three valued situation, the value $@$ is neither True, nor False, where by convention one can take the marked state $\cross{\,\,}$ as  True and the unmarked state $\cross{\cross{\,\,}} = \,\,$ as False. One interprets $XY$ as $X$ or $Y$ and the negation of $X$ as $~X = \cross{X}.$ In the algebra of this logic one no longer has the law of the excluded middle in the form
$X \cross{X} = \cross{\,\,}$ since $@ \cross{@} = @@ = @$ and $@$ is distinct from the marked state.  Similarly, in the crossing algebra we have $O \cross{O} = OO = \,\, ,$ so that 
again the law of the excluded middle is not satisfied, but in the crossing algebra the value of $P \cross{P}$ can be unmarked rather than indeterminate as in the $CSR.$\\

The $CSR$ has a natural algebra for its arithmetic and the axioms for this algebra are a generalisation of the axioms for the primary algebra of Laws of Form and directly related to Boolean algebra. A key difference from Boolean algebra in the primary algebra is the the cross is both and operator and a value in the algebra. The same holds in the $CSR$ and in the crossing algebra. The crossing algebra seems to present different problems for its axiomatization. For example, the following is an identity in contracted crossing algebra:
$$\cross{\cross{\cross{A}B}C}D = \cross{\cross{\cross{D}C}B}A$$ whenever $A,B,C,D$ are chosen from $\{O, E\}$ and one of them evaluates to $E.$ The identity is part of the larger family of identities
$$\cross{\cross{\cross{A_{1}}A_{2}} \cdots }A_{n} = \cross{\cross{\cross{A_{n}}\cdots}A_{2}}A_{1}$$  where $A_{1},A_{2},\cdots A_{n}$ are chosen from $\{O, E\}$ and one of them evaluates to $E.$ We understand a proof of these identities because the closures of rational tangles forming rational knots and links $(A_{1},A_{2}, \cdots A_{n})$ are identical topologically when the order of the terms is reversed.
Since the crossing algebra expressions count the number of components in the corresponding closure, they are equal when the number of components is two, since the only result that corresponds to two components is the unmarked evaluation $E.$ In the case of one component, it can happen that one expression equals $\cross{\,\,}$ and the other equals $O.$ For example 
$$\cross{\cross{O}O}E = \cross{OO} = \cross{\,\,}$$ while
$$\cross{\cross{E}O}O = \cross{\cross{\,\,}}O = O.$$ 
This result suggests that a full understanding of the axiomatic properties of the crossing algebra will require the use of the diagrammatics of the tangles or an equivalent structure.\\

The motivation for the crossing algebra is based on iconic representation of knots and links by diagrams, and it crosses over into iconic representations of logic and to closely related situations in topology and combinatorics. To see this more clearly, consider the  consequences of boundary interactions for curves in the plane. Here arcs of two curves that are shared cancel each other as in the arithmetic of mod 2 cycles or chains. Thus two circles that share their boundaries can cancel to become no circle at all. This can be taken to be one interpretation of $OO= \,\,$ $\,\,$ in the context where the boundary of a distinction is identified as the third value. Then one can regard $@$ and $O$ as coexisting in a larger multiple valued arithmetic where $@$ and $@$ represent disjoint circles in $@@$ while $O$ and $O$ represent interacting circles in $OO.$ These are motivating remarks at the foundations of these structures.\\

Note that at the base arithmetical level we have a generalisation of the Spencer-Brown calculus of indications with arithmetical (combinatorial) initials
$$\cross{\,\,} \cross{\,\,}  = \delta \cross{\,\,} $$
$$\cross{\cross{\,\,} }  = \,\,.$$
$$\cross{\delta} = \cross{\,\,} \delta = \delta \cross{\,\,} .$$
Here the algebraic element $\delta$ acts as a counter or memory for the number of adjacent crosses that appear in an expression.
In this way we have a dictionary between this arithmetic and the boolean nature of the original calculus where 
$$\cross{\,\,} \cross{\,\,} = \cross{\,\,}.$$
We can take $$\cross{\,\,} \cross{\,\,} \cross{\,\,} = \delta^{2} \cross{\,\,}$$ as a representative for the number $3,$ just as it survives our special evaluations as $\delta^{3}$ and is used to count components in a diagram. In this logical arithmetic, we do not reduce the $\cross{\,\,}$ to a $\delta$ since this is only appropriate for final evaluations.\\

The iconics we have used for tangles become iconics for these abstract combinatorial algebras.\\

In relating the algebraic constructions in this paper to other structures it is worth mentioning that the rule $\cross{\,\,}\cross{\,\,} = \delta \cross{\,\,}$ is, by our iconics, at the beginning of the diagrammatic Temperley-Lieb Algebra \cite{K,KL}. For in our iconics we have  $\cross{\,\,} = \cross{\Asmooth} = \Bsmooth$ and $\Bsmooth \Bsmooth = \delta \Bsmooth$ corresponds to the 
basic projector identity $U^2 = \delta U$ in the Temperley-Lieb algebra. In computing the bracket polynomial in the previous section we used the crossing algebra to count loops for the states.
In the Temperley-Lieb algebra models for the Kauffman bracket and the Jones polynomial the loop counting is accomplished by properties of the Temperley-Lieb algebra either diagrammatically or algebraically, depending on the context of the calculations. The loop structure is even more important for constructing the Khovanov complex \cite{Kho} for that link homology theory generalizing the bracket and Jones Polynomials. As we have already pointed out, it would be of interest to have an algebraic way, distinct from tensor networks,  to construct the Khovanov complex for arborescent links and their generalisations.\\

\begin{figure}[htb]
     \begin{center}
     \begin{tabular}{c}
     \includegraphics[width=10cm]{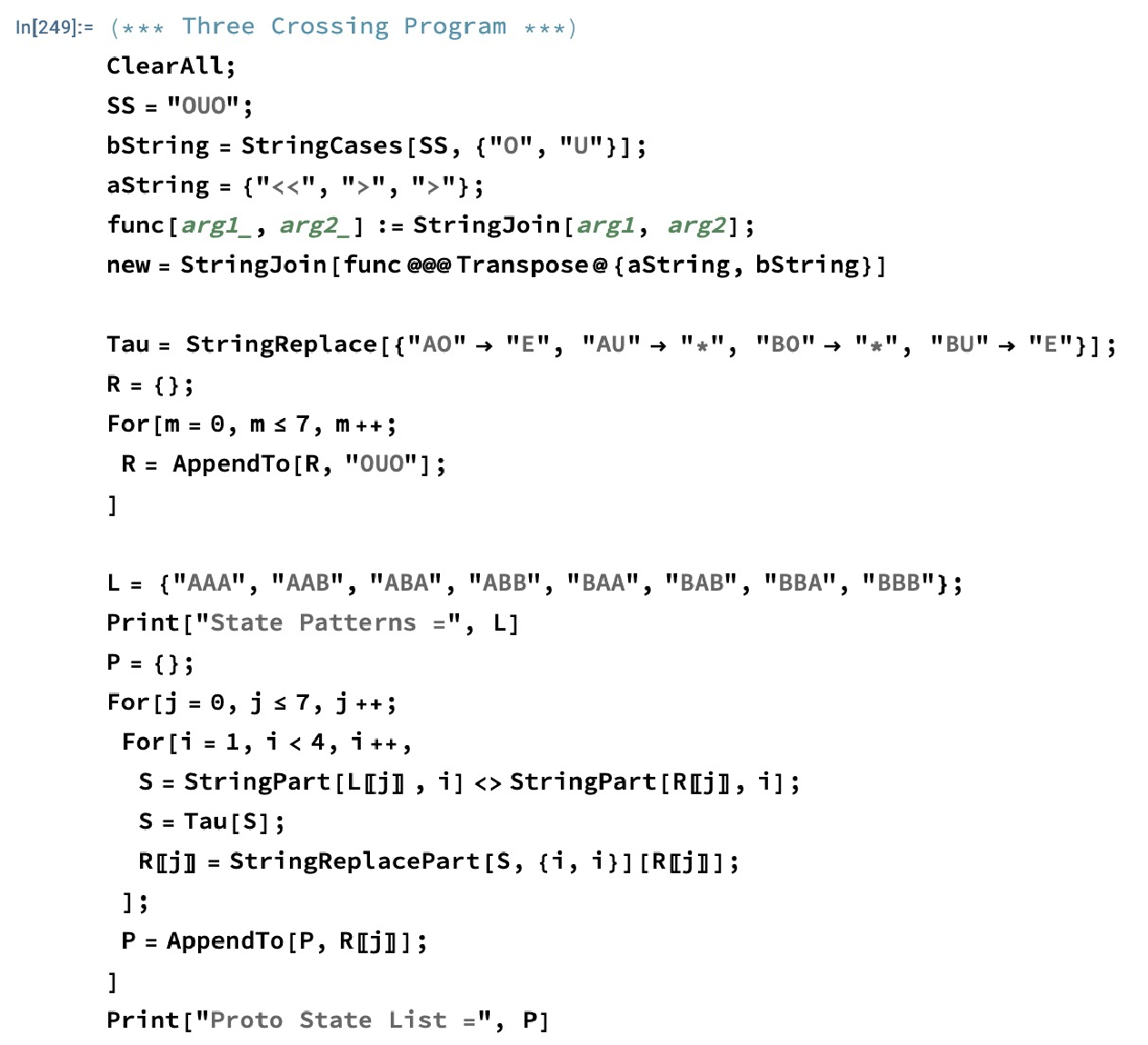}
     \end{tabular}
     \caption{\bf Program Part 1}
     \label{pg1}
\end{center}
\end{figure}

 \begin{figure}[htb]
     \begin{center}
     \begin{tabular}{c}
     \includegraphics[width=10cm]{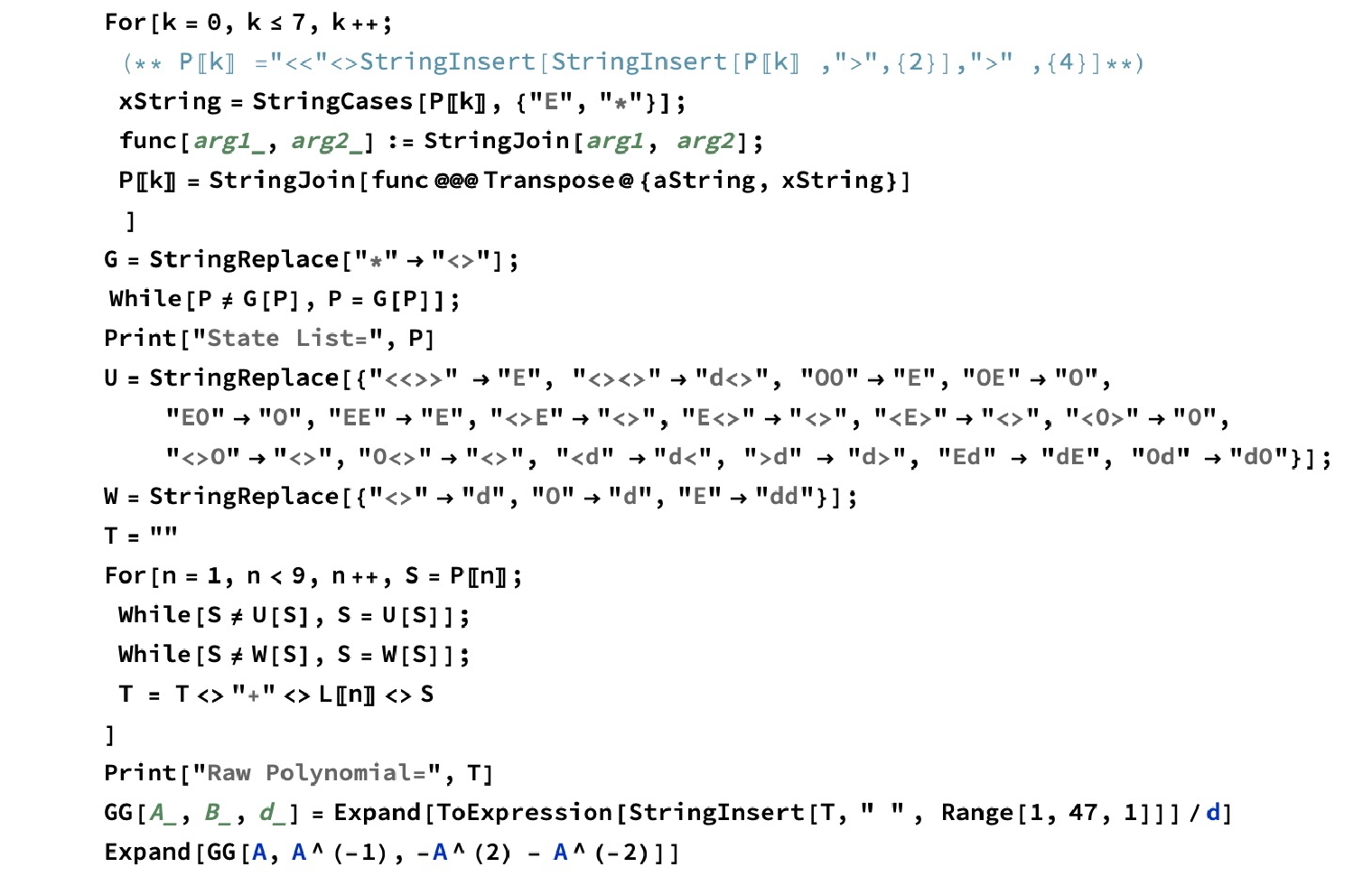}
     \end{tabular}
     \caption{\bf Program Part 2}
     \label{pg2}
\end{center}
\end{figure}

 \begin{figure}[htb]
     \begin{center}
     \begin{tabular}{c}
     \includegraphics[width=10cm]{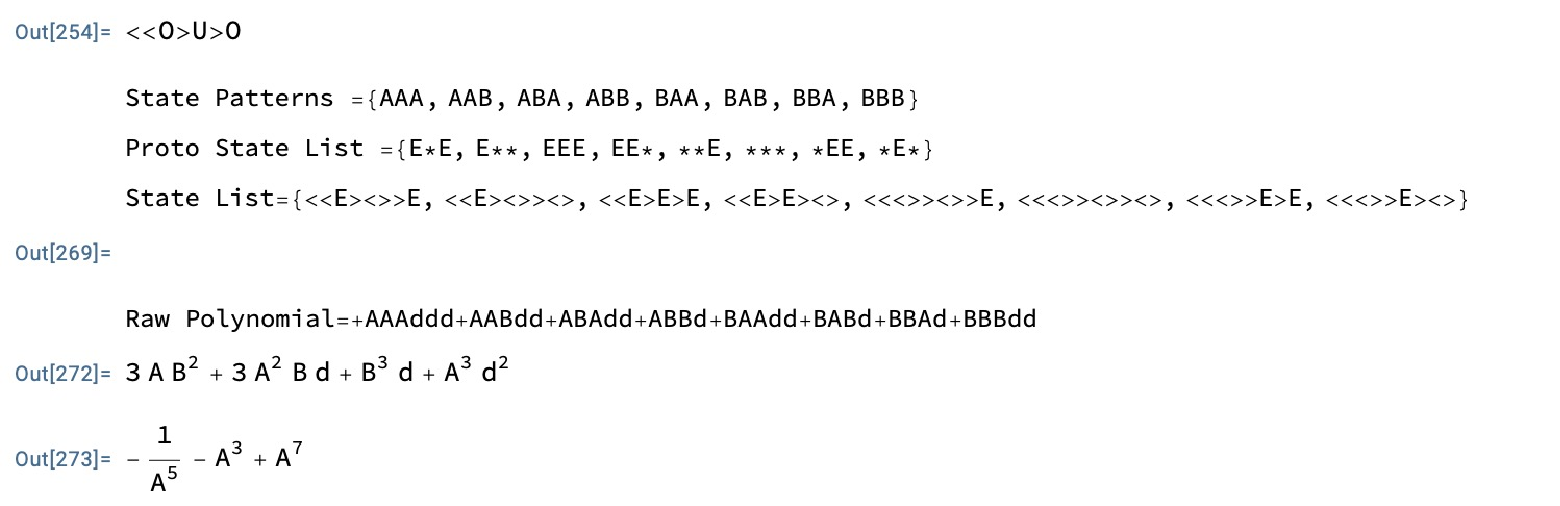}
     \end{tabular}
     \caption{\bf Output}
     \label{out}
\end{center}
\end{figure}

 \section{Appendix- More about computation.}
In Figures ~\ref{insert}, ~\ref{pg1}, ~\ref{pg2}, ~\ref{out} we illustrate a prototype program that can take as input a crossing algebra expression for an arborescent link and output the 
bracket polynomial of that link. The method is the same as our outline above, with the production of the states automated by using string manipulations. In Figure~\ref{insert} we show how a typical input code such as 
$$\cross{\cross{O}U}O = <<O>U>O$$ with $O$ standing for a crossing of type $[+1]$ and $U$ standing for a crossing of type $[-1].$ The input method uses the bracket framework separated from the
symbols $O,O,U$ and a separate string in the form $OOU.$ A program combines the symbol string and the bracket string. This same program acts multiple times in the full program illustrated in 
Figure~\ref{pg1} and Figure~\ref{pg2} to create all the strings representing the states for the bracket. Figure ~\ref{out} is the last stage of the program that produces a raw polynomial and then converts it
to the bracket polynomial by letting $B= A^{-1}$ and $dd = -A^{2} - A^{-2}.$ The reader will note that the input to this program can be made more friendly and that this particular program has an internal list of all the state forms as products of $A$'s and $B$'s. At the present time we use different lists for different crossing numbers. Better technology is in the offing.\\

\end{document}